\documentclass[final,3p,times]{elsarticle}
\usepackage{lineno,hyperref}
\usepackage{amssymb}
\usepackage{subfigure}
\usepackage{amsmath}
\usepackage{mathrsfs}
\newtheorem{theorem}{Theorem}[section]

\modulolinenumbers[5]

\journal{Wave Motion  (https://doi.org/10.1016/j.wavemoti.2018.02.001)}

\begin{document}

\begin{frontmatter}

\title{Method of Relaxed Streamline Upwinding for Hyperbolic Conservation Laws}


\author{Ameya D. Jagtap }
\cortext[mycorrespondingauthor]{Corresponding author}
\ead{ameyadjagtap@gmail.com, ameya@tifrbng.res.in}

\address{Centre for Applicable Mathematics, Tata Institute of Fundamental Research, Bangalore, Karnataka - 560065, India.}

\begin{abstract}
In this work a new finite element based Method of Relaxed Streamline Upwinding is proposed to solve hyperbolic conservation laws. Formulation of the proposed scheme is based on relaxation system which replaces hyperbolic conservation laws by semi-linear system with stiff source term also called as relaxation term. The advantage of the semi-linear system is that the nonlinearity in the convection term is pushed towards the source term on right hand side which can be handled with ease. Six symmetric discrete velocity models are introduced in two dimensions which symmetrically spread foot of the characteristics in all four quadrants thereby taking information symmetrically from all directions. Proposed scheme gives exact diffusion vectors which are very simple. Moreover, the formulation is easily extendable from scalar to vector conservation laws. Various test cases are solved for Burgers equation (with convex and non-convex flux functions), Euler equations and shallow water equations in one and two dimensions which demonstrate the robustness and accuracy of the proposed scheme. New test cases are proposed for Burgers equation, Euler and shallow water equations. Exact solution is given for two-dimensional Burgers test case which involves normal discontinuity and series of oblique discontinuities. Error analysis of the proposed scheme shows optimal convergence rate. Moreover, spectral stability analysis gives implicit expression of critical time step.
\end{abstract}


\end{frontmatter}

\linenumbers

\section{Introduction}
Many natural processes are governed by hyperbolic conservation laws like high speed flows governed by compressible Euler equations, shallow water flows like flow in a canal, river \textit{etc} are governed by shallow water equations, astrophysical flows or space weather governed by magnetohydrodynamic equations \textit{etc}. These equations describe the transport and propagation of waves (both linear and nonlinear) in space and time. Due to nonlinear nature of convection term such equations admit discontinuous solution even when the initial condition is sufficiently smooth. This precludes the possibility of finding closed form solution except in some simple cases. As a last resort, numerical methods are extensively used to solve these equations. Over the past few decades various numerical methods are proposed to solve hyperbolic conservation laws.
In the literature of finite volume and finite difference framework, upwind methods are used for solving hyperbolic conservation laws. 
Upwind methods use biased stencil in the direction determined by the sign of characteristic speeds. These methods are first proposed by Courant \textit{et al.} in \cite{CIR}.
Various upwind methods are available in the literature like flux splitting methods, flux difference splitting \textit{etc}. Upwind methods are also used in finite element framework where they are part of much larger group called as stabilized finite element methods. There are various stabilized finite element methods available in the literature like Taylor Galerkin method, Streamline-Upwind Petrov-Galerkin method, Discontinuous Galerkin method, Least-Square Galerkin method, Characteristic Galerkin \textit{etc}. A detailed discussion about these methods can be found in \cite{OCZT, JSHTW, PBGU, JJDB, TJC}. 

Relaxation schemes introduced by Jin \& Xin \cite{Jin} for hyperbolic conservation laws without source term are an attractive alternative to the upwind schemes. In the recent years, the simplicity of this scheme attracted many researchers around the world. Relaxation scheme is based on the relaxation system which replaces nonlinear convection term present in hyperbolic partial differential equation(s) by semi-linear system with stiff relaxation term on right hand side. This system is equivalent to original hyperbolic conservation law in the limit of vanishing relaxation parameter. Simple procedure present in Relaxation schemes to handle nonlinear convection term avoids more complex Riemann solvers. First and second order relaxation schemes are first introduced in \cite{Jin} while higher order relaxation schemes are discussed in \cite{HS, CYL, MSeaid}. Relaxation schemes for a hyperbolic conservation laws with source term are presented in \cite{CLL, ACh, TJF}. Natalini \cite{NaR} interpreted Jin \& Xin's relaxation system as a discrete velocity Boltzmann equation with BGK model for the collision term. In \cite{PAG} Aregba-Driollet and Natalini introduced numerical schemes based on discrete velocity Boltzmann equation which are called as discrete velocity kinetic schemes. Relaxation schemes are also employed in the lattice Boltzmann framework \cite{MBAK, GB1}. An alternative relaxation system for one-dimensional conservation law proposed by Murthy \cite{AVS} retain the semi-linear structure of original relaxation system but at the same time satisfies the integral constraint which is more consistent than the original relaxation system. For more details of relaxation schemes refer \cite{Jin, Nat, NaR, PAG, MB, BFLS, CRR, CF, HR, LMP} and the references there in. 

Relaxation scheme is also used in finite element framework for one dimensional scalar and vector (elastodynamics) problems \cite{CATM}. In this paper the relaxation system based Streamline-Upwind scheme, named as Method of Relaxed Streamline Upwinding (MRSU) is developed in finite element framework for hyperbolic conservation laws in one and two dimensions. The proposed scheme belongs to the class of stabilized finite element methods.
Some of the salient features of proposed formulation are
\begin{enumerate}
\item Weak formulation of governing hyperbolic conservation laws in conservation form is obtained from relaxation system assuming instantaneous relaxation to equilibrium.
\item Only test space of convection part of the governing equation is enriched to obtain the required stabilization term. 
\item Group discretization also called as group formulation which is shown to be more accurate is used for flux function \cite{Fle,Mor}.
\item Proposed scheme can be easily extended from scalar to vector conservation laws. 
\item Exact stabilization vectors can be obtained for both scalar as well as vector conservation laws. Importantly, computationally expensive Jacobian matrices are not involved in the stabilization terms.  
\item Six symmetric discrete velocity models are introduced in two dimensions which include four points along diagonal ($\mathscr{D}4$), nine points including rest particle ($\mathscr{AD}9$), eight points without rest particle ($\mathscr{AD}8$), four points along axis ($\mathscr{A}4$), four points along diagonal with one rest particle ($\mathscr{D}5$) and four points along axis with one rest particle ($\mathscr{A}5$) discrete velocity models.
\item To show the efficacy of the proposed scheme various test cases of Burgers equation, Euler equations and shallow water equations are solved in one and two dimensions. Moreover, some new test cases of Burgers equation, Euler and shallow water equations are also introduced. In case of two-dimensional Burgers equation a set of test cases involving normal and oblique discontinuity is proposed along with their exact solution.
\item Error analysis of the proposed scheme shows optimal rate of convergence. Furthermore, spectral stability analysis gives implicit expression of critical time step. 
\end{enumerate}

This paper is arranged as follows. After introduction in section 1, section 2 describes governing hyperbolic conservation laws like Burgers equation, Euler and shallow water equations. Section 3 gives relaxation system for hyperbolic conservation laws. In section 4, relaxed formulation of hyperbolic conservation laws is explained which will be used to develop MRSU scheme. Section 5 explains two and three discrete velocity models for one-dimensional problems whereas section 6 describes various symmetric discrete velocity models for two-dimensional problems. In section 7, Chapman-Enskog type expansion of relaxation system is performed which gives stability condition for such system. Section 8 gives weak MRSU formulation in detail. Section 9 describes temporal discretization of semi-discrete MRSU scheme followed by section 10 where simple gradient based shock capturing parameter is defined. In section 11, spectral stability analysis of the proposed scheme is carried out which gives expression of critical time step. In section 12, many numerical experiments are carried out for Burgers equation, Euler and shallow water equations which support author's claim of robustness and accuracy in the proposed numerical scheme. Finally, this paper is concluded in section 13.

\section{Governing Equations}
The governing hyperbolic conservation laws in two dimensions are
\begin{equation}\label{HCLori}
\frac{\partial U}{\partial t} + \frac{\partial G_1(U)}{\partial x}+ \frac{\partial G_2(U)}{\partial y}= 0, \,\,\, (\mathbf{x},t)  \in \Omega_T \coloneqq \Omega \times (0, \,T] \subset \mathbb{R}^2 \times \mathbb{R}_+
\end{equation}
with appropriate initial and boundary conditions. Fluxes $G_i(U), i = 1,2$ are functions of conserved variable $U$. In case of Burgers equation both $U$ and $G_i(U)$'s  are scalar quantities whereas in case of Euler and shallow water equations they are vector quantities. 

For two-dimensional Euler equations, $U$ and $G_i$'s are given as
$$ U =  \left\{ {\begin{array}{c}
\rho\\
\rho u_1\\
\rho u_2\\
\rho E\\
\end{array} } \right\}, \ \ G_i =\left\{ {\begin{array}{c}
\rho u_i \\
\delta_{i1}p+ \rho u_1u_i \\
\delta_{i2}p+ \rho u_2u_i\\
 p u_i+ \rho u_i E\\
\end{array} } \right\}, \ \ i =1, 2$$
where $\rho, u_1, u_2, E, p$ are density, velocity components in $x$ and $y$ directions, total energy and pressure respectively and $\delta_{ij}$ is Kronecker delta.  Total energy is given by
$$E = \frac{p}{\rho (\gamma-1)} + \frac{1}{2} || \mathbf{u}||_{\mathcal{L}_2}^2$$
The eigenvalues of flux Jacobian matrices $\frac{\partial G_i}{\partial U}, \,\, i =1, 2$ are $$u_i \pm a ,\,\, u_i,\,\,u_i$$
where $a = \sqrt{\gamma p/\rho}$ is acoustic speed.

For two-dimensional shallow water equations, $U$ and $G_i$'s  are given as
$$ U =  \left\{ {\begin{array}{c}
H\\
H u_1\\
H u_2\\
\end{array} } \right\}, \ \ G_i =\left\{ {\begin{array}{c}
H u_i \\
H u_1u_i+\frac{1}{2}\delta_{i1} \, g H^2\\
H u_2u_i+\frac{1}{2}\delta_{i2} \,g H^2\\
\end{array} } \right\}, \ \ i =1, 2$$
where $H, u_1, u_2,g$ are water column height, velocity components in $x$ and $y$ directions and acceleration due to gravity respectively.
Here, eigenvalues of flux Jacobian matrices $\frac{\partial G_i}{\partial U}, \,\, i =1, 2$ are 
$$u_i \pm \sqrt{gH},\,\, u_i$$

\section{Relaxation System for Hyperbolic Conservation Laws}
This section gives a brief introduction to the relaxation system for hyperbolic conservation laws. Consider two-dimensional system of hyperbolic conservation laws
\begin{equation}\label{HCL}
\frac{\partial U}{\partial t} +\frac{\partial G_1(U)}{\partial x}+ \frac{\partial G_2(U)}{\partial y} = 0, \,\,\, U(\mathbf{x},0) = U_0(\mathbf{x})
\end{equation}
where $\mathbf{x} \in \mathbb{R}^{2}$, $U \in \mathbb{R}^N$ and flux function $G_i(U) \in \mathbb{R}^N$ is nonlinear. As the above system of equations is hyperbolic, so the Jacobian $\frac{\partial G_i}{\partial U}$ is diagonalizable. Jin and Xin \cite{Jin} proposed following relaxation system for equation \eqref{HCL}
\begin{align*}
\frac{\partial U}{\partial t} + \frac{\partial W_1}{\partial x} + \frac{\partial W_2}{\partial y}  & = 0 
\\ \frac{\partial W_1}{\partial t} + \boldsymbol{\Theta}_1^2 \frac{\partial U}{\partial x} &= -\frac{1}{\epsilon} (W_1-G_1(U))
\\ \frac{\partial W_2}{\partial t} + \boldsymbol{\Theta}_2^2 \frac{\partial U}{\partial y} &= -\frac{1}{\epsilon} (W_2-G_2(U))
\end{align*}
where $W_i \in \mathbb{R}^N$, $\boldsymbol{\Theta}_i^2 \in \mathbb{R}^{N\times N}$ is a diagonal matrix with non-negative elements $\Theta_{ij}^2$ with $ i =1, 2, \cdots, \mathcal{D}$ ($\mathcal{D}$ represents number of dimensions),  $ j = 1, 2, \cdots, N$ ($N$ is number of discrete velocities) and $\epsilon$ is the relaxation time. 

The advantage of relaxation system is, convection term is linear and the nonlinearity is moved to the right hand side as source term. The solution of above relaxation system approaches solution of hyperbolic conservation law in the limit $\epsilon \rightarrow 0$ if following sub-characteristic condition is satisfied 
\begin{equation}\label{SCC}
 \frac{\sigma_1^2}{\Theta_{i1}^2} +\frac{\sigma_2^2}{\Theta_{i2}^2}+ \cdots +\frac{\sigma_N^2}{\Theta_{iN}^2} \leq 1
\end{equation}
where $\sigma_1, \sigma_2, \cdots \sigma_N$ are the eigenvalues of Jacobian matrix $\frac{\partial G_i(U)}{\partial U}$. 

In case of two-dimensional scalar conservation law, the relaxation system reads
\begin{align*}
\frac{\partial U}{\partial t} + \frac{\partial W_1}{\partial x} + \frac{\partial W_2}{\partial y}  & = 0 
\\ \frac{\partial W_1}{\partial t} + \lambda_1^2 \frac{\partial U}{\partial x} &= -\frac{1}{\epsilon} (W_1-G_1(U))
\\ \frac{\partial W_2}{\partial t} + \lambda_2^2 \frac{\partial U}{\partial y} &= -\frac{1}{\epsilon} (W_2-G_2(U))
\end{align*}
where $\lambda_1, \lambda_2$ are positive constants. It can be seen that, in the limit $\epsilon \rightarrow 0$ relaxation system for scalar and vector conservation laws leads to original conservation laws.

\section{Relaxed Formulation of Hyperbolic Conservation Laws}
This section introduces relaxed formulation of one-dimensional scalar hyperbolic conservation laws. Extension of this formulation to higher dimensions as well as for vector conservation laws would be straightforward. 

Consider one-dimensional scalar nonlinear hyperbolic conservation law
\begin{equation*}
\frac{\partial U}{\partial t} + \frac{\partial G(U)}{\partial x} = 0, \,\,\, (x,t)  \in \Omega_T \subset \mathbb{R} \times \mathbb{R}_+
\end{equation*}
with initial condition $ U(x,0) = U_0(x)$ and appropriate boundary condition. The nonlinear flux function is $G(U) = U^2/2$. Relaxation system for above equation is
\begin{align}\label{reSy}
\frac{\partial U}{\partial t} + \frac{\partial W}{\partial x} & = 0  \nonumber
\\ \frac{\partial W}{\partial t} + \lambda^2 \frac{\partial U}{\partial x} &= -\frac{1}{\epsilon} (W-G(U)) 
\end{align}
with initial condition $ U(x,0) = U_0(x), \ \ W(x,0) = G(U_0(x))$, where $\lambda$ is a positive constant. Above set of equations can be written in matrix form as

\begin{equation*}
\frac{\partial Q}{\partial t} + A\frac{\partial Q}{\partial x} = H 
\end{equation*}
where
\begin{equation*}
 Q =  \left\{ {\begin{array}{c}
U\\
W\\
\end{array} } \right\}, \ \ A =\left[ {\begin{array}{cc}
0&1\\
\lambda^2&0\\
\end{array} } \right], \ \  H =  \left\{ {\begin{array}{c}
0\\
\frac{1}{\epsilon} (G(U) - W) \\
\end{array} } \right\}
\end{equation*}
Matrix $A$ can be decomposed as
\begin{equation*}
A = R\Lambda R^{-1}
\end{equation*}
where $R$ is a model matrix whose columns are the eigvenvectors of $A$ and $\Lambda = \text{diag}\{-\lambda,\lambda \}$ is a spectral matrix with eigenvalues as diagonal entries. Introducing the characteristics variable vector $\mathbf{f} = R^{-1}Q$ above system can be written as
\begin{equation}\label{shce}
\frac{\partial \mathbf{f}}{\partial t} + \Lambda \frac{\partial \mathbf{f}}{\partial x} = R^{-1}H 
\end{equation}
where 
\begin{equation*}\mathbf{f} = R^{-1}Q  =\left\{ {\begin{array}{c}
f_1\\
f_2\\
\end{array} } \right\}= \left\{ {\begin{array}{c}
\frac{U}{2} - \frac{W}{2\lambda}\\
\frac{U}{2} + \frac{W}{2\lambda}\\
\end{array} } \right\}
\end{equation*}
Let 
\begin{equation}\label{F2v}
\mathbf{F} = \left\{ {\begin{array}{c}
F_1\\
F_2\\
\end{array} } \right\} =  \left\{ {\begin{array}{c}
\frac{U}{2} - \frac{G}{2\lambda}\\
\frac{U}{2} + \frac{G}{2\lambda}\\
\end{array} } \right\}
\end{equation}
using this, equation \eqref{shce} can be written as
\begin{equation}\label{shce1}
\frac{\partial \mathbf{f}}{\partial t} + \Lambda \frac{\partial \mathbf{f}}{\partial x} = -\frac{(\mathbf{f}-\mathbf{F})}{\epsilon} 
\end{equation}
with initial conditions $ U(x,0) = U_0(x), \ \ \mathbf{f}(0) = \mathbf{F}(U_0(x))$.

For small value of $\epsilon$ above equation looks similar to the classical Boltzmann equation with BGK model \cite{BGK} where $\epsilon$ is performing a role of relaxation time. Here $\mathbf{f}$ represents the distribution function whereas $\mathbf{F}$ represents local Maxwellian distribution or equilibrium distribution. In the above equation discrete velocities $-\lambda$ and $\lambda$ are involved whereas in a classical Boltzmann equation velocities are continuous. Due to this, Natalini \cite{NaR, PAG} interpreted relaxation system (equation \eqref{shce1}) as a discrete velocity Boltzmann equation.

Splitting method \cite{Jin} is employed to solve equation \eqref{shce1} which can be written in two steps as
\begin{align*}
\frac{\partial \mathbf{f}}{\partial t} + \Lambda \frac{\partial \mathbf{f}}{\partial x} = 0 & \ \ :\textrm{Convection}
\\ \frac{d\mathbf{f}}{dt} = -\frac{(\mathbf{f}-\mathbf{F})}{\epsilon}  & \ \ :\textrm{Relaxation}
\end{align*}
The solution of relaxation step is $\mathbf{f} = (\mathbf{f}(0)-\mathbf{F}) e^{-t/\epsilon}+\mathbf{F}$.
Assuming instantaneous relaxation to equilibrium \textit{i.e.}, $\epsilon = 0$ gives $\mathbf{f}  = \mathbf{F}$. Substituting this in convection step, we get
\begin{equation*}
\frac{\partial \mathbf{F}}{\partial t} + \Lambda \frac{\partial \mathbf{F}}{\partial x} = 0
\end{equation*}
which can be rewritten in conservation form as
\begin{equation}\label{shce21}
\frac{\partial \mathbf{F}}{\partial t} +  \frac{\partial \Lambda \mathbf{F}}{\partial x} = 0
\end{equation}
Hyperbolic conservation law is recovered from above equation by taking moments
\begin{equation*}
\frac{\partial P\mathbf{F}}{\partial t} +  \frac{\partial P (\Lambda \mathbf{F})}{\partial x} = 0
\end{equation*}
where 
\begin{subequations}\label{Mome111}
\begin{align}
P\mathbf{F} &= U \label{Mome1}
\\ P(\Lambda \mathbf{F}) & = G  \label{Mome2}
\end{align}
\end{subequations}
and moment vector $P = \underbrace{[1, 1, \cdots, 1]}_{N \,\text{times}}$. Equations \eqref{Mome1} and \eqref{Mome2} are called as \textit{Moment Relations}. Schemes based on above procedure are called as instantaneous Relaxation schemes or Relaxed schemes. In two dimensions, equation \eqref{shce21} becomes
\begin{equation}\label{shce2}
\frac{\partial \mathbf{F}}{\partial t} +  \frac{\partial \Lambda_1 \mathbf{F}}{\partial x}+\frac{\partial \Lambda_2 \mathbf{F}}{\partial y} = 0
\end{equation}
In the upcoming section, equation \eqref{shce2} will be used for deriving weak formulation of hyperbolic conservation laws.

In next two sections, we shall discuss various discrete velocity models in one and two dimensions, construct $\Lambda$ matrices for these models and obtain corresponding expressions of local Maxwellian distribution. 

\section{One-Dimensional Discrete Velocity Models}
\begin{figure} [htpb] 
\centering
\includegraphics[scale=0.63]{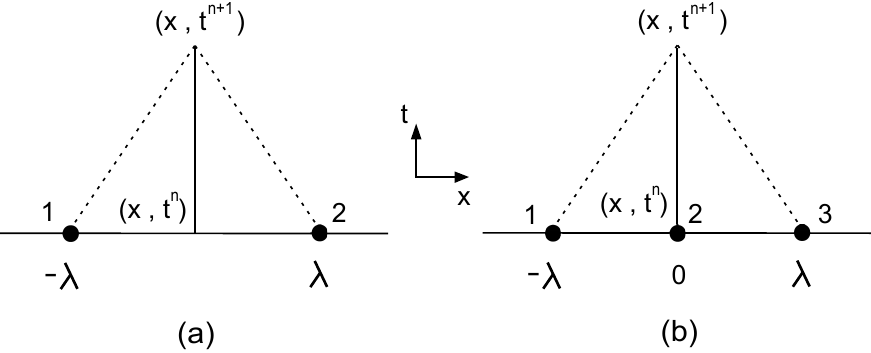}
\caption{(a) Two and (b) three discrete velocity models}
\label{fig:23dvm}
\end{figure}

Several choices of discrete velocity models are possible as long as the corresponding local Maxwellian distribution $\mathbf{F}$ satisfies the moment relations given by equation \eqref{Mome111}. In this section two models are considered for one-dimensional hyperbolic conservation laws, namely, two discrete velocity model with $N=2$ and three discrete velocity model with $N=3$. Two discrete velocity model is already introduced in the previous section where two moving particles are present. In this case $\Lambda$ is given as
$$\Lambda = \text{diag}\{-\lambda,\lambda \}$$
and the corresponding expression of $\mathbf{F}$ is given by equation \eqref{F2v} which satisfies the moment relations. Figure~\ref{fig:23dvm} (a) shows this model.

In case of three discrete velocity model there is a rest particle along with two moving particles (see figure~\ref{fig:23dvm} (b)). In this case
$\Lambda$ is given as
$$\Lambda = \text{diag}\{-\lambda,0, \lambda \}$$
and the corresponding expression of $\mathbf{F}$ is (see Appendix A for detailed derivation)
\begin{equation*}
\mathbf{F} = \left\{ {\begin{array}{c}
F_1\\
F_2\\
F_3 \\
\end{array} } \right\} =  \left\{ {\begin{array}{c}
\frac{U}{3} - \frac{G}{2\lambda}\\
\frac{U}{3} \\
\frac{U}{3} + \frac{G}{2\lambda}\\
\end{array} } \right\}
\end{equation*}
Again, $\mathbf{F}$ satisfies the moment relations. 

\section{Two-Dimensional Symmetric Discrete Velocity Models}
Asymmetric discrete velocity model for two-dimensional relaxation system is discussed in the literature \cite{PAG}. This section introduces the Symmetric Discrete Velocity (SDV) models in two dimensions which takes information symmetrically from all directions. Again, in two dimensions several choices of SDV models are possible. In general, $\Lambda$ matrices are constructed using $N$ diagonal blocks which corresponds to $N$ discrete velocities as
$$\Lambda_i = \textrm{diag}\{\lambda^1_i,\lambda^2_i, \cdots,\lambda^N_i \}$$
Here superscript indicates number of discrete velocity.

In order to admit kinetic entropy and satisfy the entropy inequality in the equilibrium limit $\epsilon \rightarrow 0$ by equation \eqref{shce1}, Bochut \cite{BF} has characterized the space of Maxwellians. The Maxwellians are written
as a linear combination of conserved variables and fluxes as
\begin{equation*}
 F_i(U) = k_{i0} U + \sum_{j=1}^{\mathcal{D}} k_{ij}G_j(U) \ \ \textrm{for} \,\, i = 1,2, \ldots N
\end{equation*}
where constants $k_{ij}$ are chosen in such a way that consistency conditions given by equation \eqref{Mome111} are satisfied. The $\Lambda$ matrices are constructed using orthogonal velocity method \cite{Jin}. Here, diagonal entries $\lambda^l_i$ 's are chosen such that they satisfy the following two conditions
\begin{align}\label{OVM2}
 & \sum_{l=1}^N \lambda^l_i = 0, \ \ \forall \,\, i = 1,2, \ldots \mathcal{D} \nonumber
 \\ & \sum_{l=1}^N \lambda^l_i\lambda^l_j  = 0, \,\,\, \textrm{where} \,\, i \neq j \,\,\, \textrm{and}  \ \ \forall \,\,i,j = 1,2, \ldots \mathcal{D} 
\end{align}
with this the Maxwellian is obtained as
\begin{equation*}
 F_k(U) = \frac{U}{N} + \sum_{j=1}^{\mathcal{D}} \frac{G_j(U)}{\left(\sum_{l=1}^N |\lambda^l_j|\right)} \textrm{sgn}(\lambda^k_j)
\end{equation*}
\begin{figure} [htpb] 
\centering
\includegraphics[ scale=0.7]{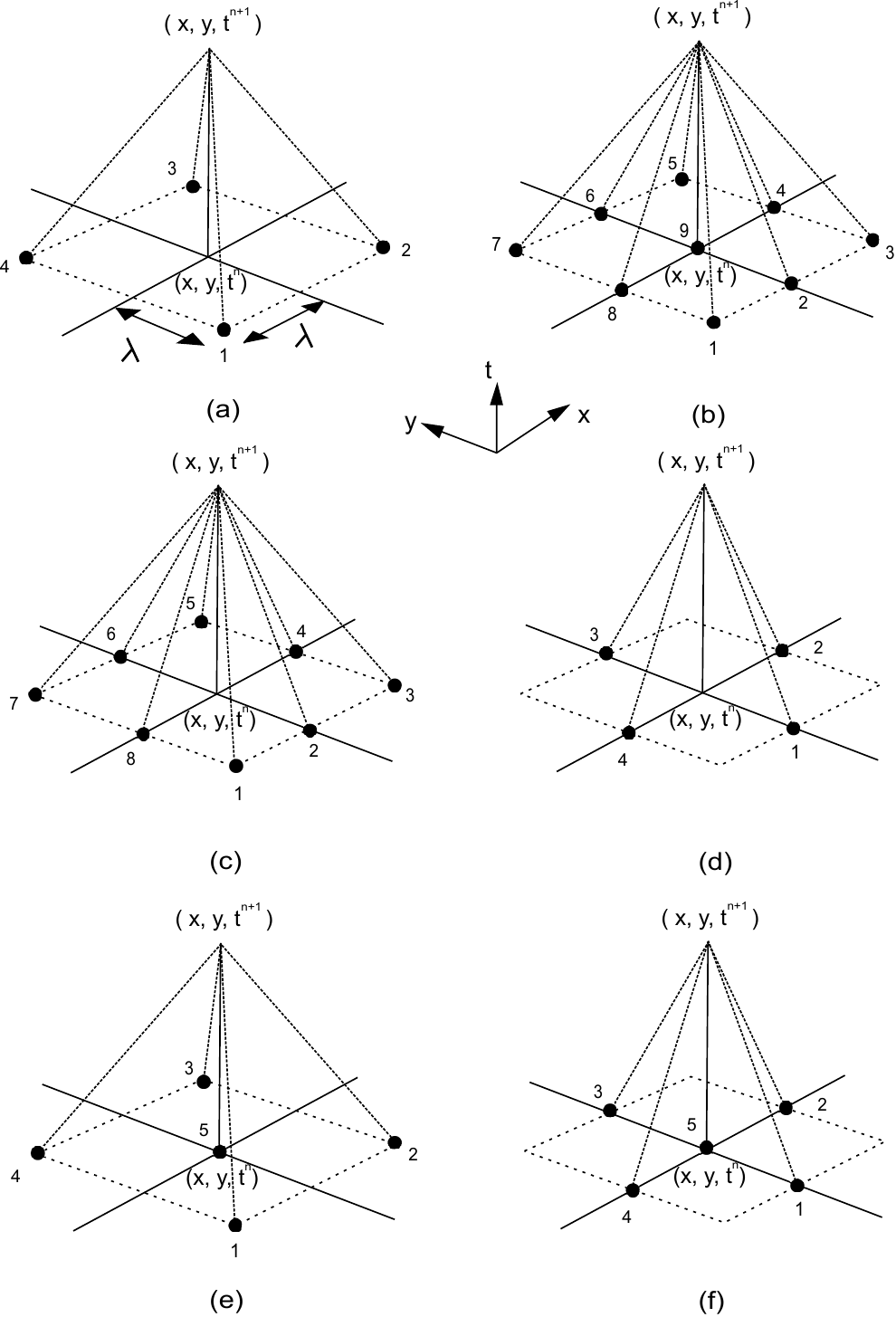}
\caption{(a) Four points along diagonal ($\mathscr{D}4$), (b) Nine points  ($\mathscr{AD}9$), (c) Eight points ($\mathscr{AD}8$), (d) Four points along axis ($\mathscr{A}4$), (e) Four points along diagonal with one rest particle ($\mathscr{D}5$) and (f) Four points along axis with one rest particle ($\mathscr{A}5$) SDV models.}
\label{fig:IsoRelx}
\end{figure}
Next, we consider different SDV models in two dimensions. Six SDV models are proposed, namely, four points along diagonal ($\mathscr{D}4$) with four discrete velocities \textit{i.e.}, $N=4$, nine points ($\mathscr{AD}9$) with $N=9$, eight points ($\mathscr{AD}8$) with $N=8$, four points along axis ($\mathscr{A}4$) with $N=4$, four points along diagonal with one rest particle ($\mathscr{D}5$) with $N=5$ and four points along axis with one rest particle ($\mathscr{A}5$) with $N=5$ SDV models as shown in figure ~\ref{fig:IsoRelx}. Hereafter for the sake of simplicity, superscript on discrete velocities are omitted. In all models, magnitude of discrete velocities (characteristic speeds) in both $x$ and $y$ direction is same which is $\lambda_1 = \lambda_2 = \lambda$ as shown in figure ~\ref{fig:IsoRelx} (a).

$\mathscr{D}4$ SDV model uses four characteristic speeds along diagonal from four quadrants as shown in figure \ref{fig:IsoRelx} (a). In this case $\Lambda_1$ and $\Lambda_2$ matrices are given as
\begin{align*}
\Lambda_1  & = \text{diag}\{-\lambda,\lambda,\lambda,-\lambda\}, \ \ \Lambda_2   = \text{diag}\{-\lambda,-\lambda,\lambda,\lambda\}  
\end{align*}
The characteristics variable and local Maxwellian distribution function are obtained as
\begin{equation*}
\mathbf{f} = \left\{ {\begin{array}{c}
f_1\\
f_2\\
f_3 \\
f_4 \\
\end{array} } \right\} =  \left\{ {\begin{array}{c}
\frac{U}{4} - \frac{W_1(U)}{4\lambda} - \frac{W_2(U)}{4\lambda}\\
\frac{U}{4} + \frac{W_1(U)}{4\lambda} - \frac{W_2(U)}{4\lambda}\\
\frac{U}{4} + \frac{W_1(U)}{4\lambda} + \frac{W_2(U)}{4\lambda}\\
\frac{U}{4} - \frac{W_1(U)}{4\lambda} + \frac{W_2(U)}{4\lambda}\\
\end{array} } \right\}, \ \ \mathbf{F} = \left\{ {\begin{array}{c}
F_1\\
F_2\\
F_3 \\
F_4 \\
\end{array} } \right\} =  \left\{ {\begin{array}{c}
\frac{U}{4} - \frac{G_1(U)}{4\lambda} - \frac{G_2(U)}{4\lambda}\\
\frac{U}{4} + \frac{G_1(U)}{4\lambda} - \frac{G_2(U)}{4\lambda}\\
\frac{U}{4} + \frac{G_1(U)}{4\lambda} + \frac{G_2(U)}{4\lambda}\\
\frac{U}{4} - \frac{G_1(U)}{4\lambda} + \frac{G_2(U)}{4\lambda}\\
\end{array} } \right\}
\end{equation*}

$\mathscr{AD}9$ SDV model uses nine characteristic speeds (including rest particle) from four quadrants as shown in figure \ref{fig:IsoRelx} (b). Here $\Lambda_1$ and $\Lambda_2$ matrices are given as
\begin{align*}
\Lambda_1  & = \text{diag}\{-\lambda,0, \lambda,\lambda,\lambda,0,-\lambda,-\lambda, 0\}, \ \ 
  \Lambda_2   =  \text{diag}\{-\lambda,- \lambda,-\lambda, 0,\lambda,\lambda,\lambda, 0, 0 \} 
\end{align*}
In this case the characteristics variable and local Maxwellian distribution function are obtained as
\begin{equation*}
\mathbf{f} = \left\{ {\begin{array}{c}
f_1\\
f_2\\
f_3 \\
f_4 \\
f_5\\
f_6\\
f_7 \\
f_8 \\
f_9\\
\end{array} } \right\} =  \left\{ {\begin{array}{c}
\frac{U}{9} - \frac{W_1(U)}{6\lambda} - \frac{W_2(U)}{6\lambda}\\
\frac{U}{9}  - \frac{W_2(U)}{6\lambda}\\
\frac{U}{9} + \frac{W_1(U)}{6\lambda} - \frac{W_2(U)}{6\lambda}\\
\frac{U}{9} + \frac{W_1(U)}{6\lambda} \\
\frac{U}{9} + \frac{W_1(U)}{6\lambda} + \frac{W_2(U)}{6\lambda}\\
\frac{U}{9} + \frac{W_2(U)}{6\lambda} \\
\frac{U}{9} - \frac{W_1(U)}{6\lambda} +\frac{W_2(U)}{6\lambda}\\
\frac{U}{9} - \frac{W_1(U)}{6\lambda}\\
\frac{U}{9} \\
\end{array} } \right\}, \ \ \mathbf{F} = \left\{ {\begin{array}{c}
F_1\\
F_2\\
F_3 \\
F_4 \\
F_5\\
F_6\\
F_7 \\
F_8 \\
F_9\\
\end{array} } \right\} =  \left\{ {\begin{array}{c}
\frac{U}{9} - \frac{G_1(U)}{6\lambda} - \frac{G_2(U)}{6\lambda}\\
\frac{U}{9}  - \frac{G_2(U)}{6\lambda}\\
\frac{U}{9} + \frac{G_1(U)}{6\lambda} - \frac{G_2(U)}{6\lambda}\\
\frac{U}{9} + \frac{G_1(U)}{6\lambda} \\
\frac{U}{9} + \frac{G_1(U)}{6\lambda} + \frac{G_2(U)}{6\lambda}\\
\frac{U}{9} + \frac{G_2(U)}{6\lambda} \\
\frac{U}{9} - \frac{G_1(U)}{6\lambda} +\frac{G_2(U)}{6\lambda}\\
\frac{U}{9} - \frac{G_1(U)}{6\lambda}\\
\frac{U}{9} \\
\end{array} } \right\}
\end{equation*}

$\mathscr{AD}8$ SDV model uses eight characteristic speeds from four quadrants as shown in figure \ref{fig:IsoRelx} (c) with
\begin{align*}
\Lambda_1  & = \text{diag}\{-\lambda,0, \lambda,\lambda,\lambda,0,-\lambda,-\lambda\} , 
\ \ 
  \Lambda_2   =  \text{diag}\{-\lambda,- \lambda,-\lambda, 0,\lambda,\lambda,\lambda, 0 \} 
\end{align*}
Here, the characteristics variable and local Maxwellian distribution function are obtained as
\begin{equation*}
\mathbf{f} = \left\{ {\begin{array}{c}
f_1\\
f_2\\
f_3 \\
f_4 \\
f_5\\
f_6\\
f_7 \\
f_8 \\
\end{array} } \right\} =  \left\{ {\begin{array}{c}
\frac{U}{8} - \frac{W_1(U)}{6\lambda} - \frac{W_2(U)}{6\lambda}\\
\frac{U}{8}  - \frac{W_2(U)}{6\lambda}\\
\frac{U}{8} + \frac{W_1(U)}{6\lambda} - \frac{W_2(U)}{6\lambda}\\
\frac{U}{8} + \frac{W_1(U)}{6\lambda} \\
\frac{U}{8} + \frac{W_1(U)}{6\lambda} + \frac{W_2(U)}{6\lambda}\\
\frac{U}{8} + \frac{W_2(U)}{6\lambda} \\
\frac{U}{8} - \frac{W_1(U)}{6\lambda} +\frac{W_2(U)}{6\lambda}\\
\frac{U}{8} - \frac{W_1(U)}{6\lambda}\\
\end{array} } \right\}, \ \ \mathbf{F} = \left\{ {\begin{array}{c}
F_1\\
F_2\\
F_3 \\
F_4 \\
F_5\\
F_6\\
F_7 \\
F_8 \\
\end{array} } \right\} =  \left\{ {\begin{array}{c}
\frac{U}{8} - \frac{G_1(U)}{6\lambda} - \frac{G_2(U)}{6\lambda}\\
\frac{U}{8}  - \frac{G_2(U)}{6\lambda}\\
\frac{U}{8} + \frac{G_1(U)}{6\lambda} - \frac{G_2(U)}{6\lambda}\\
\frac{U}{8} + \frac{G_1(U)}{6\lambda} \\
\frac{U}{8} + \frac{G_1(U)}{6\lambda} + \frac{G_2(U)}{6\lambda}\\
\frac{U}{8} + \frac{G_2(U)}{6\lambda} \\
\frac{U}{8} - \frac{G_1(U)}{6\lambda} +\frac{G_2(U)}{6\lambda}\\
\frac{U}{8} - \frac{G_1(U)}{6\lambda}\\
\end{array} } \right\}
\end{equation*}
$\mathscr{A}4$ SDV model uses four characteristic speeds along axis as shown in figure \ref{fig:IsoRelx} (d) with
\begin{align*}
\Lambda_1  & =  \text{diag}\{0,\lambda,0, -\lambda\} , \ \ \Lambda_2   = \text{diag}\{-\lambda,0, \lambda, 0\} 
\end{align*}
In this case characteristics variable and local Maxwellian distribution function are obtained as
\begin{equation*}
\mathbf{f} = \left\{ {\begin{array}{c}
f_1\\
f_2\\
f_3 \\
f_4 \\
\end{array} } \right\} =  \left\{ {\begin{array}{c}
\frac{U}{4}  - \frac{W_2(U)}{2\lambda}\\
\frac{U}{4} + \frac{W_1(U)}{2\lambda} \\
\frac{U}{4} + \frac{W_2(U)}{2\lambda}\\
\frac{U}{4} - \frac{W_1(U)}{2\lambda} \\
\end{array} } \right\}, \ \ \mathbf{F} = \left\{ {\begin{array}{c}
F_1\\
F_2\\
F_3 \\
F_4 \\
\end{array} } \right\} =  \left\{ {\begin{array}{c}
\frac{U}{4}  - \frac{G_2(U)}{2\lambda}\\
\frac{U}{4} + \frac{G_1(U)}{2\lambda} \\
\frac{U}{4} + \frac{G_2(U)}{2\lambda}\\
\frac{U}{4} - \frac{G_1(U)}{2\lambda} \\
\end{array} } \right\}
\end{equation*}

$\mathscr{D}5$ SDV model uses four characteristic speeds along diagonal from four quadrants and one rest particle as shown in figure \ref{fig:IsoRelx} (e) with
\begin{align*}
\Lambda_1  & = \text{diag}\{-\lambda,\lambda,\lambda,-\lambda, 0\}, \ \ \Lambda_2   = \text{diag}\{-\lambda,-\lambda,\lambda,\lambda, 0\}  
\end{align*}
The characteristics variable and local Maxwellian distribution function are obtained as
\begin{equation*}
\mathbf{f} = \left\{ {\begin{array}{c}
f_1\\
f_2\\
f_3 \\
f_4 \\
f_5 \\
\end{array} } \right\} =  \left\{ {\begin{array}{c}
\frac{U}{5} - \frac{W_1(U)}{4\lambda} - \frac{W_2(U)}{4\lambda}\\
\frac{U}{5} + \frac{W_1(U)}{4\lambda} - \frac{W_2(U)}{4\lambda}\\
\frac{U}{5} + \frac{W_1(U)}{4\lambda} + \frac{W_2(U)}{4\lambda}\\
\frac{U}{5} - \frac{W_1(U)}{4\lambda} + \frac{W_2(U)}{4\lambda}\\
\frac{U}{5}  \\
\end{array} } \right\}, \ \ \mathbf{F} = \left\{ {\begin{array}{c}
F_1\\
F_2\\
F_3 \\
F_4 \\
F_5 \\
\end{array} } \right\} =  \left\{ {\begin{array}{c}
\frac{U}{5} - \frac{G_1(U)}{4\lambda} - \frac{G_2(U)}{4\lambda}\\
\frac{U}{5} + \frac{G_1(U)}{4\lambda} - \frac{G_2(U)}{4\lambda}\\
\frac{U}{5} + \frac{G_1(U)}{4\lambda} + \frac{G_2(U)}{4\lambda}\\
\frac{U}{5} - \frac{G_1(U)}{4\lambda} + \frac{G_2(U)}{4\lambda}\\
\frac{U}{5} \\
\end{array} } \right\}
\end{equation*}

Finally, $\mathscr{A}5$ SDV model uses four characteristic speeds along axis and one rest particle as shown in figure \ref{fig:IsoRelx} (f) with
\begin{align*}
\Lambda_1  & =  \text{diag}\{0,\lambda,0, -\lambda, 0\} , \ \ \Lambda_2   = \text{diag}\{-\lambda,0, \lambda, 0, 0\} 
\end{align*}
In this case characteristics variable and local Maxwellian distribution function are obtained as
\begin{equation*}
\mathbf{f} = \left\{ {\begin{array}{c}
f_1\\
f_2\\
f_3 \\
f_4 \\
f_5 \\
\end{array} } \right\} =  \left\{ {\begin{array}{c}
\frac{U}{5}  - \frac{W_2(U)}{2\lambda}\\
\frac{U}{5} + \frac{W_1(U)}{2\lambda} \\
\frac{U}{5} + \frac{W_2(U)}{2\lambda}\\
\frac{U}{5} - \frac{W_1(U)}{2\lambda} \\
\frac{U}{5} \\
\end{array} } \right\}, \ \ \mathbf{F} = \left\{ {\begin{array}{c}
F_1\\
F_2\\
F_3 \\
F_4 \\
F_5 \\
\end{array} } \right\} =  \left\{ {\begin{array}{c}
\frac{U}{5}  - \frac{G_2(U)}{2\lambda}\\
\frac{U}{5} + \frac{G_1(U)}{2\lambda} \\
\frac{U}{5} + \frac{G_2(U)}{2\lambda}\\
\frac{U}{5} - \frac{G_1(U)}{2\lambda} \\
\frac{U}{5}\\
\end{array} } \right\}
\end{equation*}

The word '\textit{Symmetric}' in discrete velocity model needs further explanation. All SDV models are not only symmetric with respect to the $x$ and $y$ axes but they also satisfy following two conditions.
\begin{itemize}
\item All models satisfy the property of rotational invariance for rotation angle $\theta = \pm\frac{\xi \pi}{2}, \,\xi = 1,2,\cdots$.
\item The number of non-zero diagonal elements $N_z$ in $\Lambda_i, i=1,2$ must be equal, \textit{i.e.}, $\left.N_z\right|_{\Lambda_1} = \left.N_z\right|_{\Lambda_2}$.
\end{itemize}
In all SDV models the macroscopic variables are recovered by taking moments as $P\mathbf{F} = U, \ \ P(\Lambda_i \mathbf{F}) = G_i(U)$. In the upcoming sections, these models will be compared based on expression of diffusion term and the corresponding stability condition. Models with minimum numerical diffusion will be chosen for computation in both one and two dimensions.

\section{Chapman-Enskog Type Expansion of Relaxation System}
A Chapman-Enskog type expansion of relaxation system provides the condition under which relaxation system is a dissipative approximation to the given hyperbolic conservation laws \cite{CLL}. Chapman-Enskog type expansion of relaxation system of one-dimensional scalar conservation law with two discrete velocity model (equation \eqref{reSy}) is given as (see Appendix B)
\begin{align*}
\frac{\partial U}{\partial t} + \frac{\partial G(U)}{\partial x} &= \epsilon \frac{\partial }{\partial x} \left( \frac{\partial U}{\partial x} [\lambda^2 - (G'(U))^2] \right)+ \mathcal{O}(\epsilon^2)
\end{align*}
where $G'(U) = \frac{\partial G(U)}{\partial U}$. First term on right hand side represent the viscous dissipation term with coefficient of viscosity. Therefore, relaxation system provides a vanishing viscosity model for the original hyperbolic conservation laws. For stability, value of $\lambda$ should be chosen such that $ \lambda^2 \geq |(G'(U))|^2$. In case of three discrete velocity model the stability condition become $\frac{2}{3}\lambda^2 \geq |(G'(U))|^2 $.

Similarly, to obtain the stability condition for two-dimensional scalar conservation laws one can write the general Chapman-Enskog type expansion of two-dimensional relaxation system as
\begin{align*}
\frac{\partial U}{\partial t} + \frac{\partial G_1(U)}{\partial x}+& \frac{\partial G_2(U)}{\partial y} =  \epsilon \frac{\partial }{\partial x} \left( \frac{\partial U}{\partial x} [P(\Lambda_1^2 \mathbf{F}') - (G_1'(U))^2]+ \frac{\partial U}{\partial y}[P(\Lambda_1 \Lambda_2\mathbf{F}') - G_1'(U)G_2'(U)]\right) \nonumber
\\ & + \epsilon \frac{\partial }{\partial y} \left( \frac{\partial U}{\partial x} [P(\Lambda_2 \Lambda_1\mathbf{F}') - G_2'(U)G_1'(U)]+ \frac{\partial U}{\partial y} [P(\Lambda_2^2 \mathbf{F}') - (G_2'(U))^2]\right)  + \mathcal{O}(\epsilon^2)
\end{align*}
Table \ref{Table1}  gives expressions of diffusion coefficient terms $P(\Lambda_1^2 \mathbf{F}') $, $P(\Lambda_1 \Lambda_2 \mathbf{F}')$, $P(\Lambda_2 \Lambda_1 \mathbf{F}')$ and $P(\Lambda_2^2 \mathbf{F}')$ for various SDV models.
\begin{table}[h!]
\begin{center}
\small \begin{tabular}{ccccc} \hline
 SDV Model &  $P(\Lambda_1^2 \mathbf{F}') $& $P(\Lambda_1 \Lambda_2 \mathbf{F}')$& $P(\Lambda_2 \Lambda_1 \mathbf{F}')$&  $P(\Lambda_2^2 \mathbf{F}')$ \\  \hline
 $\mathscr{D}4$ & $\lambda^2$ & 0 &0 & $\lambda^2$   \\  
  $\mathscr{AD}9$& $\frac{2}{3}\lambda^2  $ & 0 &0 &$\frac{2}{3}\lambda^2$ \\
  $\mathscr{AD}8$& $\frac{3}{4}\lambda^2  $ & 0 &0 & $\frac{3}{4}\lambda^2  $\\
 $\mathscr{A}4$ & $\frac{1}{2}\lambda^2 $& 0 &0 &$\frac{1}{2}\lambda^2 $ \\
 $\mathscr{D}5$ & $\frac{4}{5}\lambda^2$ & 0 &0 & $\frac{4}{5}\lambda^2$   \\ 
 $\mathscr{A}5$ & $\frac{2}{5}\lambda^2$ & 0 &0 & $\frac{2}{5}\lambda^2$   \\  \hline 
 \end{tabular}
\caption{Diffusion coefficients for various SDV models.}\label{Table1}
\end{center}
\end{table}

Therefore, to make dissipation positive following condition must be satisfied
\begin{equation*}
\beta \lambda^2 \geq ( |(G_1'(U))|^2 +|(G_2'(U))|^2 ) 
\end{equation*}
where $\beta$ is the coefficient of numerical diffusion. $\beta = 1, 2/3, 3/4 ,1/2, 4/5, 2/5$ for $\mathscr{D}4$, $\mathscr{AD}9$, $ \mathscr{AD}8$, $\mathscr{A}4$, $\mathscr{D}5$, $\mathscr{A}5$ models respectively (see table \ref{Table1}).

For vector conservation laws the stability condition can be obtained following \cite{CLD, RN1,WQX}.
It can be shown that the solution of relaxation system approaches exact solution of original hyperbolic partial differential equation as $\epsilon \rightarrow 0$ if the sub-characteristic condition given by equation \eqref{SCC} is satisfied. In case of proposed SDV models, the unity value present on right hand side of the sub-characteristic condition is replaced by $\beta$. But, to balance stability and numerical diffusion, following values of $\Theta_{ij}$ based on the supremum eigenvalue of Jacobian matrices are chosen for Euler and shallow water equations.

For one-dimensional Euler equations 
\begin{equation*}
\Theta_{11} = \Theta_{12} = \cdots = \Theta_{1N} = \textrm{max}(\textrm{sup} |u +a|, \textrm{sup} |u|,\textrm{sup} |u -a|)
\end{equation*}
For one-dimensional shallow water equations
\begin{equation*}
\Theta_{11} = \Theta_{12} = \cdots = \Theta_{1N}  = \textrm{max}(\textrm{sup} |u +\sqrt{gH}|,\textrm{sup} |u -\sqrt{gH}|)
\end{equation*}
For two-dimensional Euler equations
\begin{align*}
\Theta_{11} = \Theta_{12} = \cdots = \Theta_{1N}  &= \textrm{max}(\textrm{sup} |u_1 +a|, \textrm{sup} |u_1|,\textrm{sup} |u_1 -a|)
\\ \Theta_{21} = \Theta_{22} = \cdots = \Theta_{2N}  &= \textrm{max}(\textrm{sup} |u_2 +a|, \textrm{sup} |u_2|,\textrm{sup} |u_2 -a|)
\end{align*}
For two-dimensional shallow water equations
\begin{align*}
\Theta_{11} = \Theta_{12} = \cdots = \Theta_{1N}  &= \textrm{max}(\textrm{sup} |u_1 +\sqrt{gH}|, \textrm{sup} |u_1|,\textrm{sup} |u_1 -\sqrt{gH}|)
\\ \Theta_{21} = \Theta_{22} = \cdots = \Theta_{2N}  &= \textrm{max}(\textrm{sup} |u_2 +\sqrt{gH}|, \textrm{sup} |u_2|,\textrm{sup} |u_2 -\sqrt{gH}|)
\end{align*} 

\section{Weak MRSU formulation of Hyperbolic Conservation Laws}
Weak MRSU formulation starts with equation \eqref{shce2} in conservation form. The standard Galerkin finite element approximation for equilibrium distribution function $\mathbf{F}$ and flux function $\Lambda_j\mathbf{F}$ are
\begin{equation*}
\mathbf{F} \approx \mathbf{F}^h = \sum_{\forall i} N_i^h \mathbf{F}_i^h, \ \ \ \  \Lambda_j\mathbf{F} \approx (\Lambda_j\mathbf{F})^h = \sum_{\forall i} N_i^h (\Lambda_j\mathbf{F})_i^h
\end{equation*}
where group formulation is used for flux function. Interpolation function $N^h \in C_c^0(\Omega^h)$ and computational domain $\Omega^h$ is divided into \textit{Nel} number of elements as
$ \Omega^h = \bigcup_{e= 1}^{\textit{Nel}} \Omega_e^h \ \ \textrm{such that} \ \ \Omega_i^h  \cap \Omega_j^h = \O , \,\,\, \forall i \neq j $.

Defining test and trial functions as 
\begin{align*}\mathcal{V}^h&=\{N^h~\in~\mathcal{H}^{1}( \Omega^h) \linebreak \, \textrm{and} \, \ N^h  = 0 \,\,\text{on}\,\, \Gamma_D \}
\\ \mathcal{S}^h  &=  \{  \mathbf{F}^h \in \mathcal{H}^{1}( \Omega^h) \, \textrm{and} \, \mathbf{F}^h = \mathbf{F}^h_D \,\,\text{on}\,\, \Gamma_D \}  
\end{align*}
($\Gamma_D$ is the Dirichlet boundary) the weak formulation is, find $\mathbf{F}^h \in \mathcal{S}^h $ such that $\forall \, N^h \in \mathcal{V}^h$
\begin{align*}
\int_{\Omega^h} N^h \cdotp   \left(\frac{\partial \mathbf{F}^h}{\partial t}  +    \frac{\partial (\Lambda_1 \mathbf{F})^h}{\partial x} +\frac{\partial (\Lambda_2 \mathbf{F})^h}{\partial y}\right)\, d\Omega^h + \sum_{e=1}^{\textit{Nel}} \int_{\Omega_e^h} & \left[ \tau_1 \Lambda_1^h \frac{\partial N^h}{\partial x} +\tau_2 \Lambda_2^h \frac{\partial N^h}{\partial y}\right] \cdotp \left(    \frac{\partial (\Lambda_1 \mathbf{F})^h}{\partial x} + \frac{\partial (\Lambda_2 \mathbf{F})^h}{\partial y} \right)\, d\Omega_e^h 
\\ & + \sum_{e=1}^{\textit{Nel}}  \int_{\Omega_e^h}  \delta^e \left(  \frac{\partial N^h}{\partial x}\cdotp\frac{\partial \mathbf{F}^h}{\partial x} + \frac{\partial N^h}{\partial y}\cdotp\frac{\partial \mathbf{F}^h}{\partial y}\right) \, d\Omega_e^h  = 0
\end{align*}

where $\tau_i = \frac{1}{\mathcal{D}|\Lambda_i^h|} \frac{1}{\sqrt{\mathcal{D}}} \left| \frac{\partial \mathbf{x}}{\partial \mathbf{r}}\right|, i = 1,2$ 
with $\mathcal{D},\,\,\mathbf{x} $ and $\mathbf{r}$ being dimension, physical coordinates and natural coordinates respectively.
$\delta^e$ is the shock capturing parameter which will be defined later. First term is a standard Galerkin approximation. In the second expression, term with square brackets is the enriched part of the test space which gives required numerical diffusion. The last term is a shock capturing term which is used only for two-dimensional vector conservation laws and will be activated near high gradient region.

By taking moments one can obtain weak MRSU formulation of hyperbolic conservation laws
\begin{align}\label{meq1}
\int_{\Omega^h} N^h \cdotp   \left(\frac{\partial P\mathbf{F}^h}{\partial t}  +    \frac{\partial P(\Lambda_1 \mathbf{F})^h}{\partial x} +\frac{\partial P(\Lambda_2 \mathbf{F})^h}{\partial y}\right)\, d\Omega^h + \sum_{e=1}^{\textit{Nel}} \int_{\Omega_e^h} & P\left[ \tau_1 \Lambda_1^h \frac{\partial N^h}{\partial x} +\tau_2 \Lambda_2^h \frac{\partial N^h}{\partial y}\right] \cdotp \left(    \frac{\partial (\Lambda_1 \mathbf{F})^h}{\partial x} + \frac{\partial (\Lambda_2 \mathbf{F})^h}{\partial y} \right)\, d\Omega_e^h \nonumber
\\ & + \sum_{e=1}^{\textit{Nel}}  \int_{\Omega_e^h}  \delta^e \left(  \frac{\partial N^h}{\partial x}\cdotp\frac{\partial P\mathbf{F}^h}{\partial x} + \frac{\partial N^h}{\partial y}\cdotp\frac{\partial P\mathbf{F}^h}{\partial y}\right) \, d\Omega_e^h  = 0
\end{align}
The moments are given as
$$P\mathbf{F}^h = U^h, \ \ P(\Lambda_i \mathbf{F})^h = G_i^h $$

As discussed in the introduction part, only test space of convection term of the governing equation is enriched which gives required stabilization. Exact stabilization matrices are obtained for both scalar as well as vector conservation laws. Relaxation system based approach to obtain weak formulation of hyperbolic conservation laws circumvent the process of evaluation of Jacobian matrices at each time step. Thus, proposed scheme is computationally less expensive especially, for system of hyperbolic conservation laws. Appendix C gives weak formulation of system of hyperbolic conservation laws in conservation form (without using relaxation system) which involves these Jacobian matrices.  

Next, we evaluate coefficients of diffusion given by second term of equation \eqref{meq1}. In one dimension the second term can be written as
$$ \sum_{e=1}^{\textit{Nel}} \int_{\Omega_e^h}  P\left[  \tau \Lambda^h \frac{\partial N^h}{\partial x} \right] \cdotp \left(    \frac{\partial (\Lambda \mathbf{F})^h}{\partial x} \right)\, d\Omega_e^h =  \frac{1}{\mathcal{D}\sqrt{\mathcal{D}}} \left| \frac{\partial \mathbf{x}}{\partial \mathbf{r}}\right| \,D \{ P\,(\textrm{sgn}(\Lambda^h)(\Lambda \mathbf{F})^h) \}$$
where diffusion matrix $D = \int_{\Omega^h}  \left(  \frac{\partial N^h}{\partial x} \right)^T   \frac{\partial N^h}{\partial x} \, d\Omega^h$. Here we used the definition of $\tau$ and term $P\,(\textrm{sgn}(\Lambda^h)(\Lambda \mathbf{F})^h)$ is the coefficient of diffusion. Similarly, one can find diffusion coefficients in two dimensions.

\subsection{One-Dimensional Equation(s)} The expression of diffusion coefficient using two discrete velocity model is obtained as
$$ P\,(\textrm{sgn}(\Lambda^h)(\Lambda \mathbf{F})^h)  =  \begin{cases} \lambda^h U^h & \text{for scalar conservation law} \\  \Theta_{1j}^h U^h, \,\,\,  j = 1, 2 & \text{for vector  conservation laws} \end{cases}$$
whereas for three discrete velocity model
$$ P\,(\textrm{sgn}(\Lambda^h)(\Lambda \mathbf{F})^h)  = \begin{cases} \frac{2}{3}\lambda^h U^h & \text{for scalar conservation law} \\  \frac{2}{3}\Theta_{1j}^h U^h, \,\,\,  j = 1, 2, 3 & \text{for vector  conservation laws} \end{cases}
$$
For one-dimensional relaxation system three discrete velocity model is less diffusive for both scalar as well as vector conservation laws hence, hereafter this model will be used for all one-dimensional computations.

\subsection{Two-Dimensional Equation(s)}
Table \ref{Table2} shows the expressions of diffusion as well as cross-diffusion terms involved in MRSU scheme for scalar and vector conservation laws using various SDV models. In case of vector conservation laws $j = 1, 2, \cdots, N$.
\begin{table}[htpb]
\begin{center}
\small \begin{tabular}{ccccc} \hline 
 SDV Model  &  $P\,(\textrm{sgn}(\Lambda_1^h)(\Lambda_1 \mathbf{F})^h) $& $P\,(\textrm{sgn}(\Lambda_1^h)(\Lambda_2 \mathbf{F})^h) $& $P\,(\textrm{sgn}(\Lambda_2^h)(\Lambda_1 \mathbf{F})^h)$&  $P\,(\textrm{sgn}(\Lambda_2^h)(\Lambda_2 \mathbf{F})^h) $ \\  \hline
 & Scalar $|$ Vector & Scalar $|$ Vector & Scalar $|$ Vector  & Scalar $|$ Vector \\ 
 $\mathscr{D}4$ & $\lambda^h U^h$ \,\,\,\, $\Theta_{1j}^h U^h$ & 0 \,\,\,\,\, 0 &0 \,\,\,\,\, 0 & $\lambda^h U^h$ \,\,\,\, $\Theta_{2j}^h U^h$   \\  
  $\mathscr{AD}9$& $\frac{2}{3}\lambda^h U^h $ \,\,\,\, $\frac{2}{3}\Theta_{1j}^h U^h$& 0 \,\,\,\,\, 0 &0 \,\,\,\,\, 0 &$\frac{2}{3}\lambda^h U^h$ \,\,\,\, $\frac{2}{3}\Theta_{2j}^h U^h$ \\
  $\mathscr{AD}8$& $\frac{3}{4}\lambda^h U^h $ \,\,\,\, $\frac{3}{4}\Theta_{1j}^h U^h$ & 0 \,\,\,\,\, 0 &0 \,\,\,\,\, 0 & $\frac{3}{4}\lambda^h U^h $ \,\,\,\, $\frac{3}{4}\Theta_{2j}^h U^h$\\
 $\mathscr{A}4$ &    $\frac{1}{2}\lambda^h U^h $\,\,\,\,  $\frac{1}{2}\Theta_{1j}^h U^h$& 0 \,\,\,\,\, 0 &0 \,\,\,\,\, 0 &$\frac{1}{2} \lambda^h U^h $\,\,\,\,  $\frac{1}{2}\Theta_{2j}^h U^h$ \\
 $\mathscr{D}5$ &    $\frac{4}{5}\lambda^h U^h $\,\,\,\,  $\frac{4}{5}\Theta_{1j}^h U^h$& 0 \,\,\,\,\, 0 &0 \,\,\,\,\, 0 &$\frac{4}{5} \lambda^h U^h $\,\,\,\,  $\frac{4}{5}\Theta_{2j}^h U^h$ \\
 $\mathscr{A}5$ &    $\frac{2}{5}\lambda^h U^h $\,\,\,\,  $\frac{2}{5}\Theta_{1j}^h U^h$& 0 \,\,\,\,\, 0 &0 \,\,\,\,\, 0 &$\frac{2}{5} \lambda^h U^h $\,\,\,\,  $\frac{2}{5}\Theta_{2j}^h U^h$ \\ \hline 
 \end{tabular}
\caption{Diffusion and cross-diffusion terms for various SDV models.}\label{Table2}
\end{center}
\end{table}
It can be seen that the diffusion vectors are proportional to the supremum eigenvalues of the Jacobian matrices. 
One can use global supremum eigenvalue but, it produces excessive numerical diffusion especially, when the supremum eigenvalue (characteristics speed) is large compared to the average speed in the given domain which in turn smears discontinuity. Also, large global supremum eigenvalue reduces the stable time step which increases computational cost. Due to these reasons a local supremum eigenvalue of the Jacobian matrix is used for solving problems. Second observation from table \ref{Table2} is, among all two-dimensional SDV models $\mathscr{A}5$ model gives minimum diffusion for both scalar as well as vector conservation laws hence, hereafter this model will be used for all computations. 

Values of $\lambda^h$, $\Theta_{1j}^h$ and $ \Theta_{2j}^h $ for one and two-dimensional scalar as well as vector conservation laws are chosen according to the stability conditions discussed previously in section 7.

\begin{theorem}
In MRSU scheme using SDV models, the constant factor involved in coefficient of numerical diffusion terms (excluding cross-diffusion terms) is given by $\frac{N_z}{N}$ and it satisfies following inequality $$\frac{N_z}{N} \leq 1$$ where $N_z$ is the number of non-zero diagonal elements of $\Lambda$ matrix and $N$ is the number of discrete velocities.

Proof: In general the diffusion term in MRSU scheme is written as
\begin{align*}
P\,(\textrm{sgn}(\Lambda_i^h)(\Lambda_i \mathbf{F})^h)   & = \begin{cases} \frac{N_z}{N} \lambda^h U^h & \text{for scalar conservation law} \\ \frac{N_z}{N} \Theta_{ij}^h U^h, j = 1, 2, \cdots, N & \text{for vector conservation laws} \end{cases}
\end{align*}
where we used expression of $\mathbf{F}^h$ for a discrete velocity model. Since $N_z \leq N$ always, this implies $\frac{N_z}{N} \leq 1$.
\end{theorem}

Above mentioned inequality becomes equality only for two discrete velocity model in one dimension and $\mathscr{D}4$ model in two dimensions where $N_z = N$.


\begin{theorem}
Cross-diffusion vectors $P\,(\textrm{sgn}(\Lambda_1^h)(\Lambda_2 \mathbf{F})^h)  $ and $P\,(\textrm{sgn}(\Lambda_2^h)(\Lambda_1 \mathbf{F})^h) $ vanishes in all symmetric models. 

Proof: It is the direct consequence of condition given by equation \eqref{OVM2}
\end{theorem}

\section{Temporal Discretization}
For temporal discretization of semi-discrete scheme a simple forward Euler time discretization is used.
$$ \frac{\partial U^h}{\partial t} = \frac{(U^h)^{n+1}- (U^h)^n}{\Delta t}  + \mathcal{O}(\Delta t) $$
The fully discretized system of equations is solved using Generalized Minimal Residual (GMRES) method.

Main aim of this paper is to present new relaxation based stabilized finite element method for hyperbolic conservation laws. Here, first order scheme is presented but, proposed scheme is extendable to any high order using higher order time integration like strong stability preserving Runge-Kutta methods along with higher order interpolation function in space.

\section{Shock Capturing Parameter}
For two-dimensional vector conservation laws, diffusion along streamline direction is not sufficient to suppress the oscillations near high gradient regions. Hence, additional diffusion term with a shock capturing parameter is required which can sense these high gradient regions and add necessary amount of diffusion. There are many shock capturing parameters available in the literature \cite{TH, TS}.  In this work a simple gradient based shock capturing parameter $\delta^{e}$ is presented as follows.  

\begin{figure} [h!] 
\centering
\includegraphics[scale=0.65]{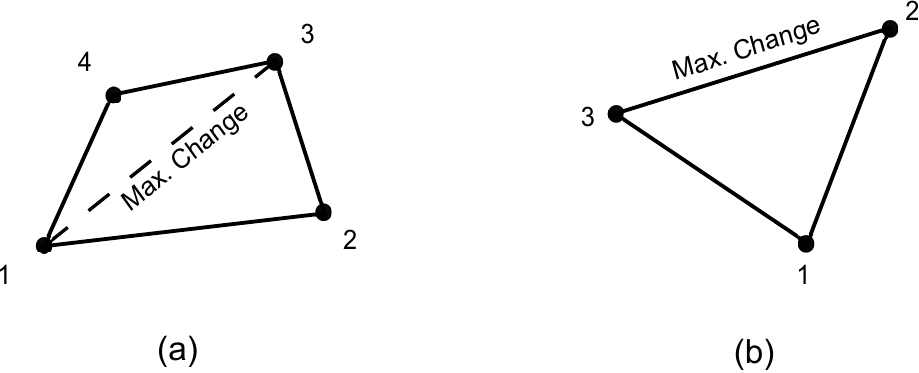}
\caption{Four node quadrilateral and three node triangular elements in physical domain}
\label{fig:4nq}
\end{figure}

Figure ~\ref{fig:4nq} (a) shows a typical four node quadrilateral element. As shown in figure, the maximum change in $\Phi$ occur across node 1 and 3 (where $\Phi$ could be density, temperature, pressure or even water column height; in present work, density is used for Euler equations and water column height is used for shallow water equations). The element based shock capturing parameter is then defined for node 1 and 3 as
\begin{equation*}
\delta_{1}^{e} = \delta_{3}^{e} = \frac{1}{\mathcal{D}\sqrt{\mathcal{D}}} \left| \frac{\partial \mathbf{x}}{\partial \mathbf{r}}\right|\left(\frac{ || \nabla \Phi ||_{\infty}^{e} }{ ||  \Phi ||_{\infty}^{e}  }\right)
\end{equation*}
where subscripts 1 and 3 represent node numbers.  For nodes 2 and 4, it is defined as
\begin{equation*}
\delta_i^{e} = \frac{1}{\mathcal{D}\sqrt{\mathcal{D}}} \left| \frac{\partial \mathbf{x}}{\partial \mathbf{r}}\right| \left(\frac{ \Phi_{\text{Max}}^{e} - \Phi_i^{e} }{ ||\Phi ||_{\infty}^{e}  }\right) \ \  \text{where} \,\, i =2,4.  
\end{equation*}
where $\Phi_{\text{Max}}^{e}$ and $ \Phi_i^{e} $ are maximum value of $\Phi$ and value of $\Phi$ at $i^{th}$ node point within the element respectively. The norm $||\Phi ||_{\infty}^{e}$ is defined as
$$ ||\Phi ||_{\infty}^{e} \coloneqq \text{max} (|\Phi_1|,\cdots, |\Phi_r|)$$
where $r$ represent total number of nodes in the element.

At element level matrix form, the shock capturing parameter is given by
\begin{equation*}
\delta^{e} = \text{diag}\,\{ \delta_{1}^{e},\delta_{2}^{e},\delta_{3}^{e},\delta_{4}^{e}\} 
\end{equation*}
The upper and lower bounds on the value $|| \nabla\Phi ||_{\infty}^{e}$ are  
\begin{equation*}
0 \leq || \nabla \Phi ||_{\infty}^{e} \leq ||\Phi ||_{\infty}^{e}
\end{equation*}
It is important to note that, addition of extra shock capturing term in weak formulation makes the formulation inconsistent with the original equation. Hence $\delta^{e}$ is defined in such a way that, as the element size $\partial \mathbf{x} \rightarrow 0$, $\delta^{e}$ should disappear. This condition is achieved by including $\partial \mathbf{x}$ in the numerator, which vanishes as we refine the mesh.  Similarly, one can define such delta parameter for triangular elements shown in figure ~\ref{fig:4nq} (b).

Calculation of global $\delta$ from element level $\delta^e$ at every time step can be computationally expensive. Moreover, it can be shown that it stagnates the convergence of the scheme \cite{LCC}. The remedy to this problem is to freeze $\delta$ if
$$ |\textrm{RES(current) - RES(previous)}| \,\,<\,\, \textrm{Tol}$$
where RES and Tol are residue and desired tolerance value respectively. This procedure results in the stagnation free convergence. Further, this procedure can be used for both steady-state and transient problems.

\section{Spectral Stability Analysis}
Explicit numerical schemes are conditionally stable due to restriction given by CFL criteria \cite{CFL}. From computational point of view it is important to find the maximum stable time step also called as critical time step $\Delta t_{\textrm{cr}}$ within which the scheme is stable. Various methods of stability analysis are available in the literature like von-Neumann, spectral stability analysis \textit{etc}. For further details about stability analysis, refer \cite{RM}. Unlike von-Neumann stability analysis, spectral stability analysis also includes boundary points. In this section a spectral stability analysis of MRSU scheme is performed by stating following theorem.

\begin{theorem}
The critical time step $\Delta t_{\textrm{cr}}$ of explicit MRSU scheme for one-dimensional linear convection equation with periodic boundary conditions satisfies following inequality
 \begin{equation}\label{Dcr}
\Delta t_{\textrm{cr}} \leq \frac{1}{\varrho\left(\frac{I}{\Delta t} -M^{-1} \left[c C  + \frac{h}{2} \lambda^h D  \right] \right) }
\end{equation}
where $M ,C$ and $D$ are mass, convection and diffusion matrices given as
$$M = \int_{\Omega^h} (N^h)^T N^h \, d\Omega^h, \,\,C = \int_{\Omega^h} (N^h)^T \left(\frac{dN^h}{dx}\right) \, d\Omega^h, \,\,D = \int_{\Omega^h} \left(\frac{dN^h}{dx}\right)^T\left(\frac{dN^h}{dx}\right) \, d\Omega^h $$
and $\varrho$, $c$, $h$ represent spectral radius of a matrix, constant wave speed, elemental length respectively. 

Proof:  Fully discrete MRSU scheme for one-dimensional linear convection equation with periodic boundary conditions and appropriate initial condition is
\begin{equation}\label{sSa}
(U^h)^{n+1}= \left[I -M^{-1}\Delta t \left(c C  + \frac{h}{2}  \lambda ^h D\right) \right] (U^h)^n
\end{equation}

Let $\mathcal{A} \coloneqq \left[I -M^{-1}\Delta t \left(c C + \frac{h}{2}  \lambda^h D\right) \right] $
be an amplification matrix, then above equation becomes 
\begin{equation*}
(U^h)^{n+1}= \mathcal{A} (U^h)^n
\end{equation*}
Using recurrence relation one can write
\begin{equation*}
(U^h)^{n+1}= \mathcal{A}^{n+1} (U^h)^0
\end{equation*} 
where $(U^h)^0$ is the initial solution. For bounded solution, amplification matrix must satisfy
$$||\mathcal{A}^{n+1}|| \leq 1, \ \  \forall n   \in \mathbb{N} $$
Thus, we get $ ||\mathcal{A}|| \leq 1 \Rightarrow  \varrho(\mathcal{A}) \leq 1$
where $\varrho(\mathcal{A})$ is the spectral radius of amplification matrix.
After simplification the critical time step is obtained as
 \begin{equation*}
\Delta t_{\textrm{cr}} \leq \frac{1}{\varrho\left(\frac{I}{\Delta t} -M^{-1} \left[c C + \frac{h}{2}  \lambda^h D \right] \right) }
\end{equation*}
where we used the property of spectral radius of a matrix $\varrho(\alpha\mathcal{A}) = \alpha \varrho(\mathcal{A}), \forall \alpha \in \mathbb{R}$.
\end{theorem}
This result can be easily extended for higher-dimensional linear convection equation.
\begin{figure} [h!] 
\centering
\includegraphics[scale=0.46]{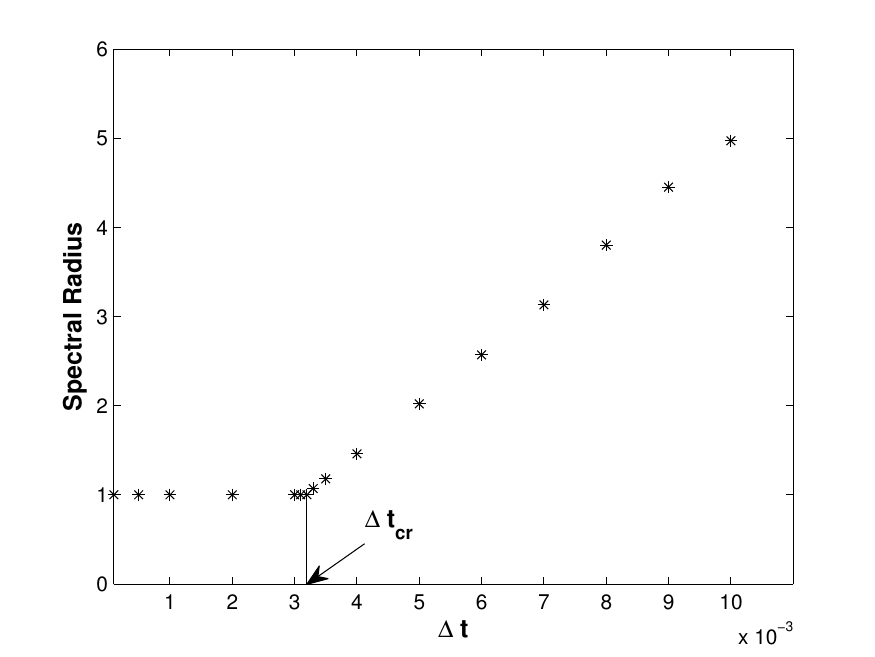}
\caption{Variation of spectral radius with $\Delta t$.}
\label{fig:Stability}
\end{figure}

As an example, lets consider one-dimensional linear convection equation with unity wave speed. The computational domain is $[0,\,1]$. Initial condition is given as
\begin{equation}\label{Cpul}
U(x,0) = \begin{cases} \frac{1}{2} \left(1+\cos \left[\frac{\pi(x-0.2)}{0.12}\right]\right) & \text{for} \,\, |x-0.2| \leq 0.12 \\   0 & \text{Otherwise}  \end{cases}
\end{equation}
and boundary conditions are periodic. 
The spectral radius of amplification matrix is computed numerically with $50$ node points. Figure \ref{fig:Stability} shows variation of spectral radius with respect to $\Delta t$ which gives the value of $\Delta t_{\textrm{cr}}$ around $3.2 \times 10^{-3}.$

\section{Numerical Experiments}
In this section various one and two-dimensional test cases are solved for Burgers equation, Euler and shallow water equations which shows the accuracy and robustness of the proposed scheme. These are the standard test cases which are chosen based on the complexity of the solution. Moreover, some new test cases are also introduced. 

\textbf{Remark}: Residue plots are given for some steady-state test cases where residue is calculated using relative $\mathcal{L}_2$ error as
\begin{equation*} \text{Residue} = \frac{||(U^h)^{n+1}- (U^h)^n||_{\mathcal{L}_2}}{||(U^h)^{n+1}||_{\mathcal{L}_2}} 
\end{equation*}

\subsection{Error Analysis using One-Dimensional Convection Equation}
Linear Lagrange interpolation function is used as a basis function in space. Experimental Order of Convergence (EOC) is calculated for one-dimensional convection equation with initial condition as a cosine wave (equation \eqref{Cpul}) convecting with unity wave speed in $\mathcal{L}_2$ and $\mathcal{H}_1$ norm. Table ~\ref{Table3} shows the EOC which is optimal for linear interpolation function $N^h \in C_c^0(\Omega^h)$.
\begin{table}[h!]
\begin{center}
\small \begin{tabular}{ccccc} \hline
No. of Nodes &  $\mathcal{L}_2$& EOC & $\mathcal{H}_1$&  EOC \\ \hline 
40 &0.1348 & - & 3.2934 &  - \\  
80& 0.0809 & 0.7366 &2.1803 &0.5951 \\
160 & 0.0451& 0.8430 &1.3226 &0.7211 \\
320 &0.0240   &0.9101&0.7521 & 0.8144  \\  
640& 0.0123  &0.9644 & 0.4134 &0.8634 \\
1280 & 0.00615 &1.0000 & 0.2195&0.9133\\
2560 & 0.00302 &1.0260 &0.1154 &0.9276\\ \hline 
 \end{tabular}
\caption{Convergence analysis of MRSU scheme.}\label{Table3}
\end{center}
\end{table}

\subsection{1-D Inviscid Burgers Test Cases with Smooth and Discontinuous Initial Data}
\begin{figure} [h!] 
\centering
\includegraphics[scale=0.37]{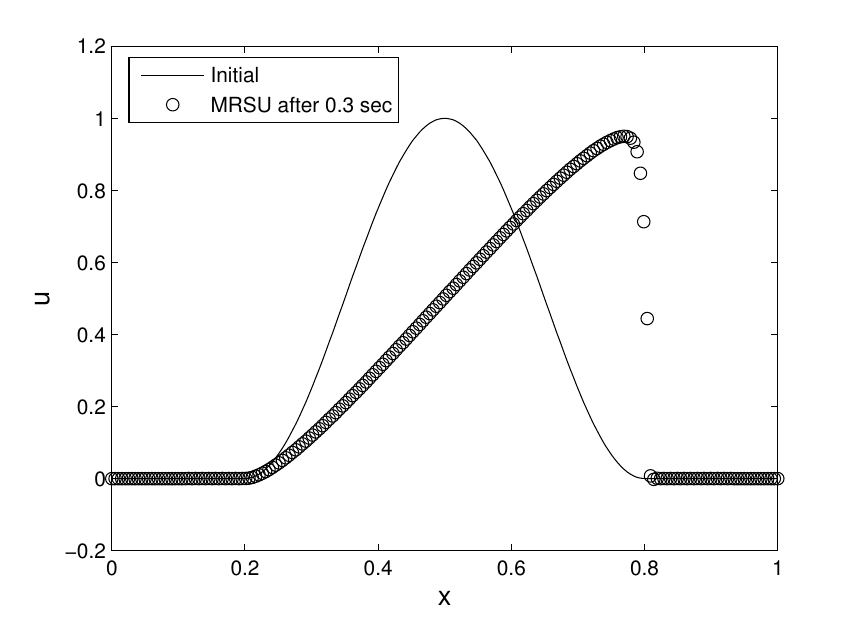}
\includegraphics[scale=0.37]{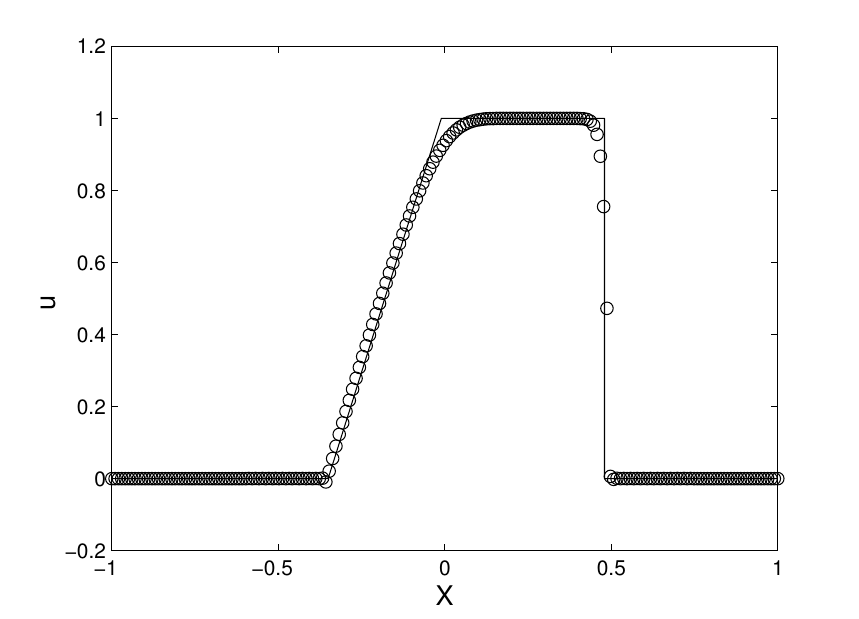}
\caption{Solution of one-dimensional Burgers equation using MRSU scheme.}
\label{fig:krgE11}
\end{figure} 
One-dimensional inviscid Burgers equation is given by
$$ \frac{\partial U}{\partial t} +  \frac{\partial U^2/2}{\partial x} = 0$$ The domain is $[0,\,\,1]$ for smooth initial condition and $[-1,\,\,1]$ for discontinuous initial condition.
Smooth initial condition is given by cosine pulse
\begin{equation*}
U(x,0) = \begin{cases} \frac{1}{2} \left(1+\cos \left[\frac{\pi(x-0.5)}{0.3}\right]\right) & \text{for} \,\, |x-0.5| \leq 0.3 \\   0 & \text{Otherwise}  \end{cases}
\end{equation*}
whereas discontinuous initial condition is represented by a square wave as
\begin{equation*}
U(x,0) = \begin{cases} 1 & \text{for} \,\, |x| < 1/3 \\   0 & \text{for}  \,\, 1/3 < |x| \leq 1\end{cases}
\end{equation*}
In both cases, node points are 200, final time is 0.3 and CFL number is 0.5. Figure ~\ref{fig:krgE11} (left) shows smooth initial profile (solid line) and MRSU solution after 0.3 seconds with circles. Figure ~\ref{fig:krgE11} (right) shows discontinuous exact solution (solid line) and MRSU solution (circles).

\subsection{Sod's Shock Tube Problem}
\begin{figure} [h!] 
\centering
\includegraphics[scale=0.35]{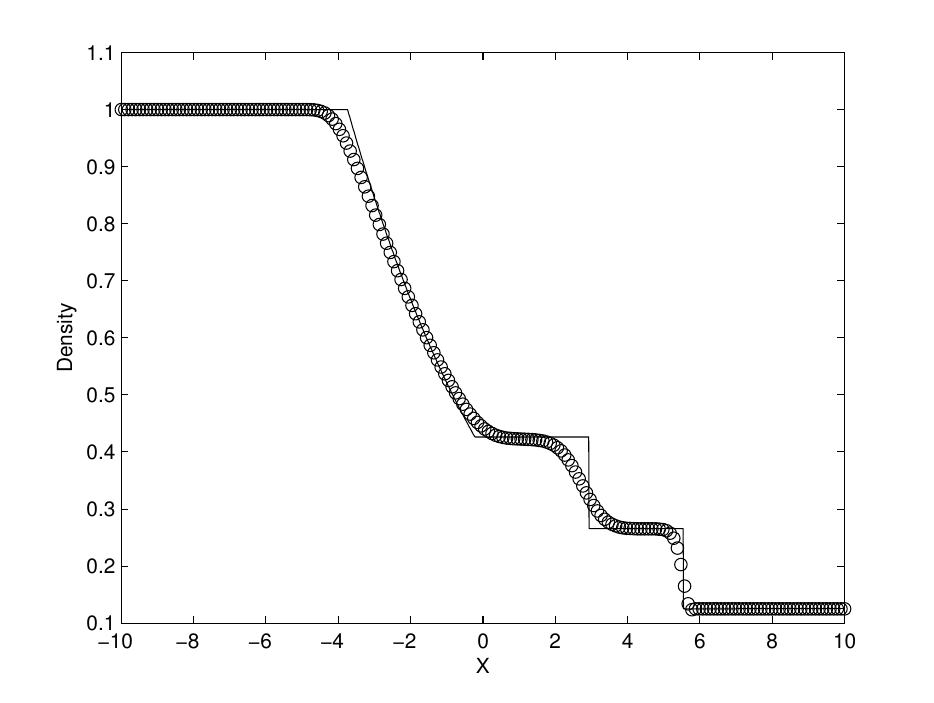}
\includegraphics[scale=0.35]{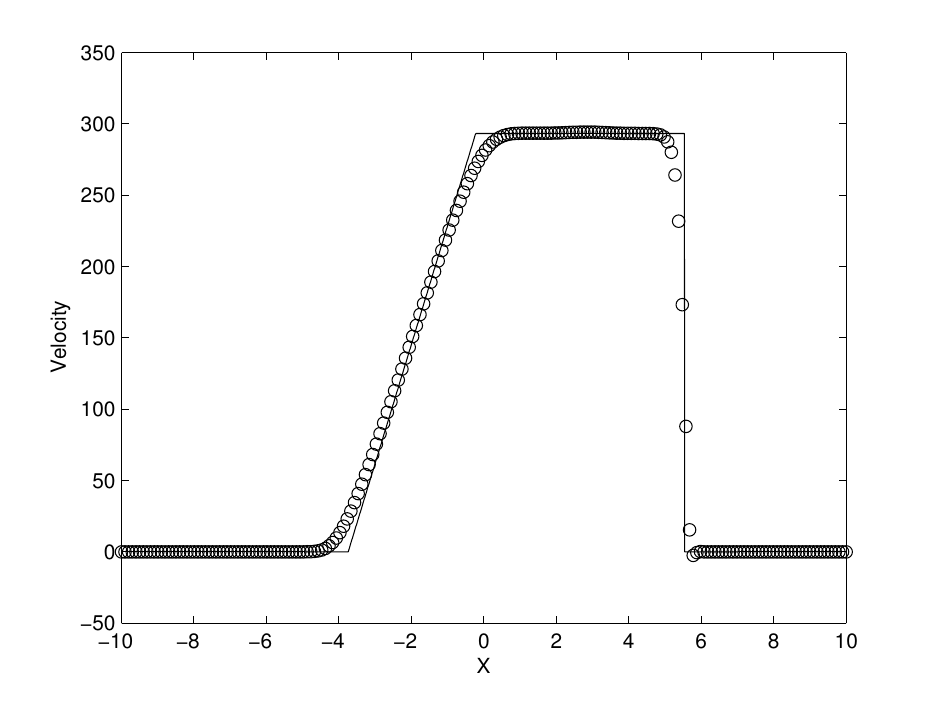}
\linebreak
\includegraphics[scale=0.35]{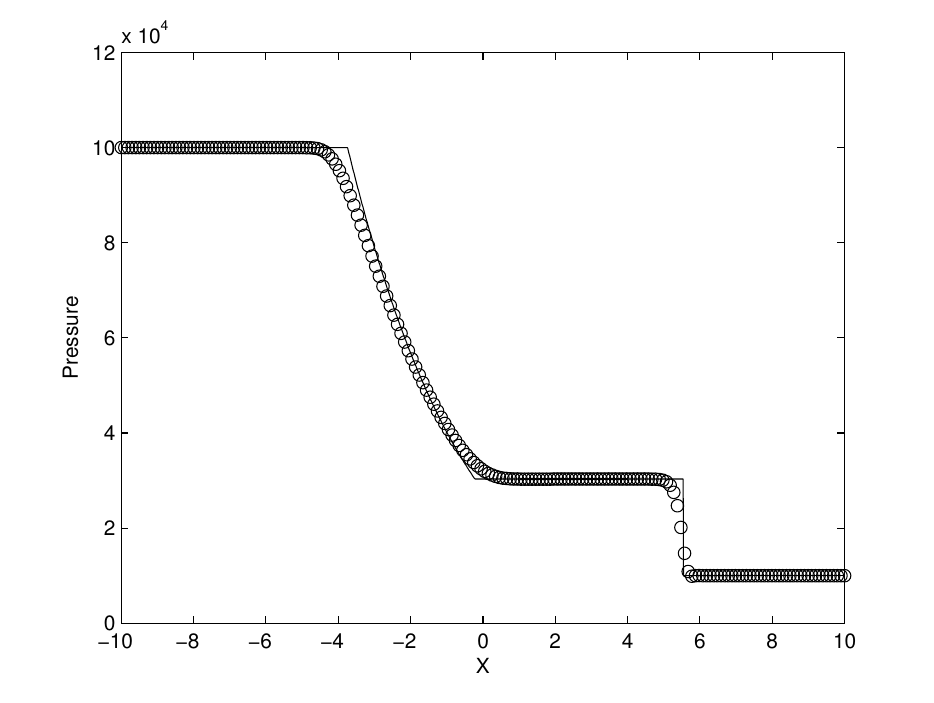}
\includegraphics[scale=0.35]{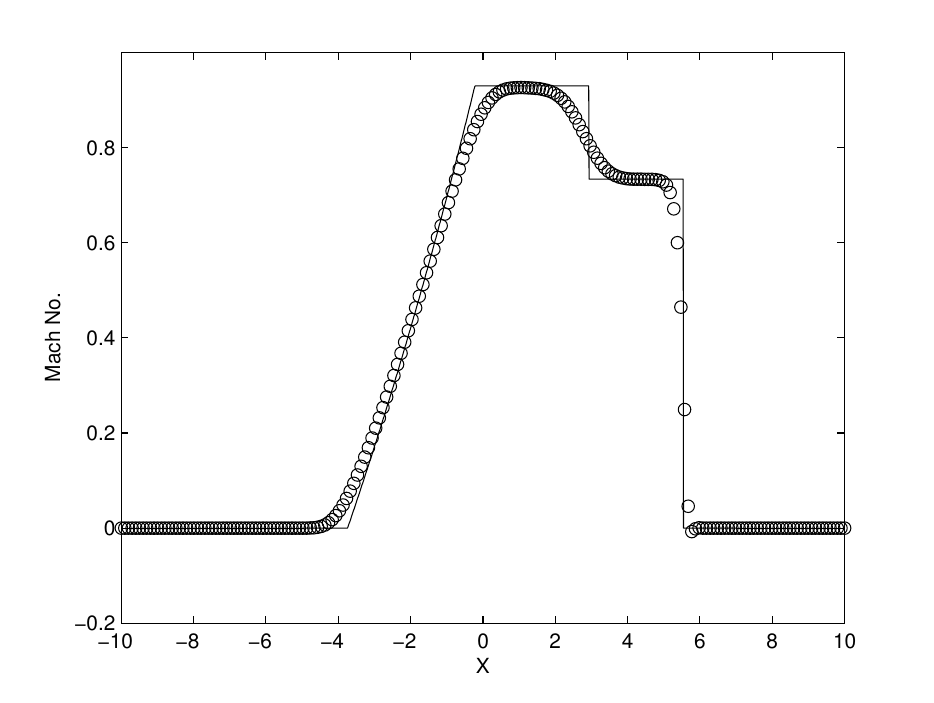}
\caption{Density, velocity, pressure and Mach number plots.}
\label{fig:krgE12}
\end{figure}
This test case is taken from Laney \cite{CL}. Here, the domain is [$-10,\,\,10$]. Sod's shock tube problem consists of a left rarefaction, a right shock wave and a contact discontinuity which separates the rarefaction and shock wave. 
The initial conditions are given by 
\begin{equation*}
(\rho, u, p)(x,0) = \begin{cases} 1,0,100000 & \text{If}\,\,\, -10<x<0 \\ 0.125,0,10000 & \text{If}\,\,\, 0\leq x<10\end{cases}
\end{equation*}
The number of node points is 200, CFL number is 0.2 and final time is 0.01. Figure ~\ref{fig:krgE12} shows density, velocity, pressure and Mach number plots. The solid line represents exact solution while numerical solution is given by circles.  Here, all the essential flow features like expansion wave, contact discontinuity and shock wave are captured. 

\subsection{Shock-Entropy Wave Interaction \cite{CWSo}}
In this test case a moving shock wave with Mach number 3 interacts with sinusoidal density profile. The domain is [$0,\,\,10$]. The initial conditions are given as
\begin{figure} [h!] 
\centering
\includegraphics[scale=0.33]{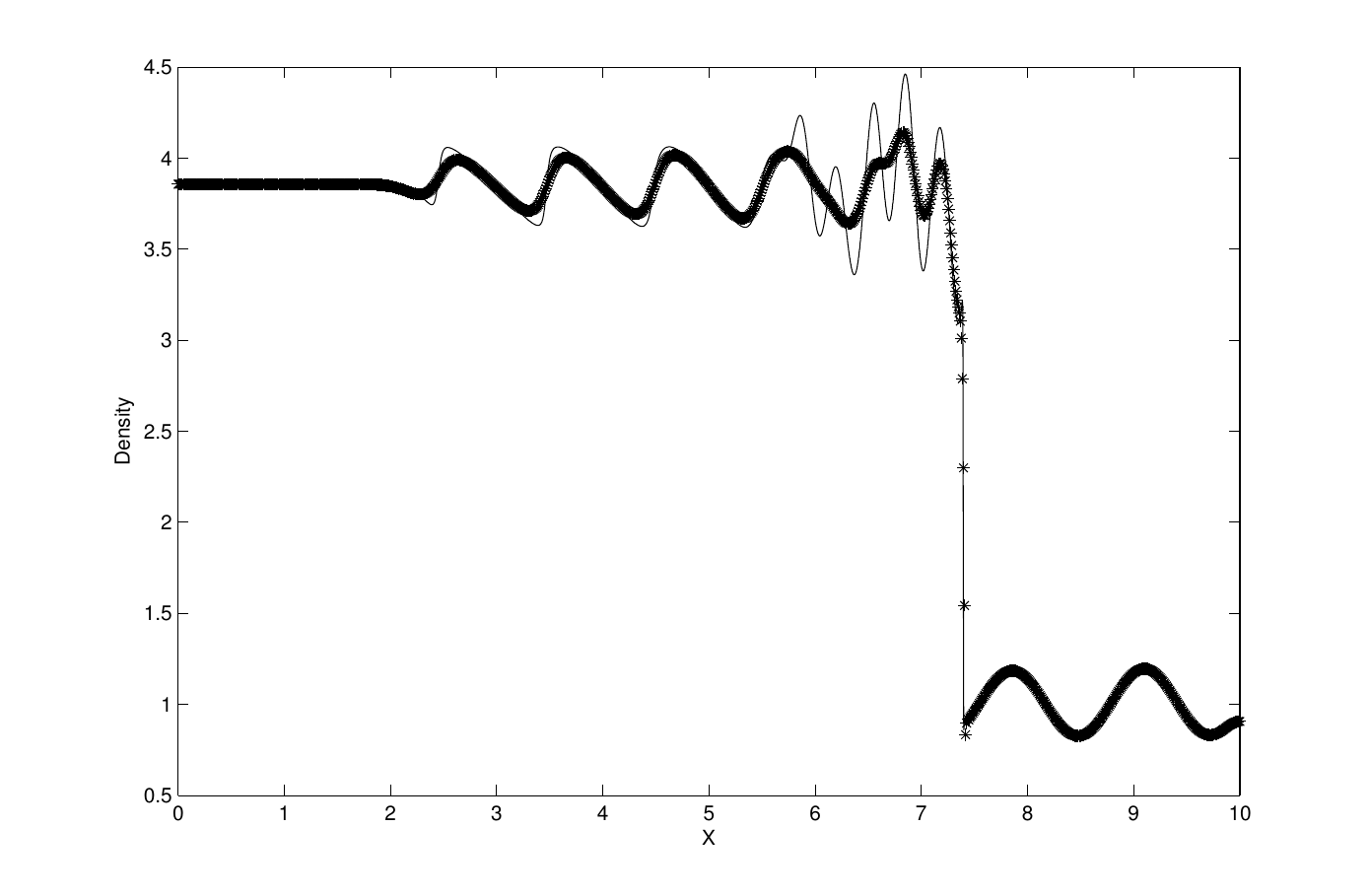}
\caption{Density plot for shock-entropy wave interaction.}
\label{fig:sew12}
\end{figure}
\begin{equation*}
(\rho, u, p ) (x,0) = \begin{cases} (3.857,2.629,10.3333) & \text{for} \,\, x < 1 \\   (1+0.2\sin(5x),0,1) & \text{for}  \,\, x \geq 1\end{cases}
\end{equation*}
Final time is 1.8. Number of nodes used is 1000 and CFL number is 0.4. This test case involves both shock wave and smooth profile. MRSU solution (represented by *) is compared with reference solution (solid line). 


\subsection{Woodward and Colella Blastwave Problem \cite{WCPP}}
It is one of the severe test case used to test the robustness and accuracy of the numerical scheme. The domain is [$0,\,\,1$]. 
\begin{figure} [h!] 
\centering
\includegraphics[scale=0.34]{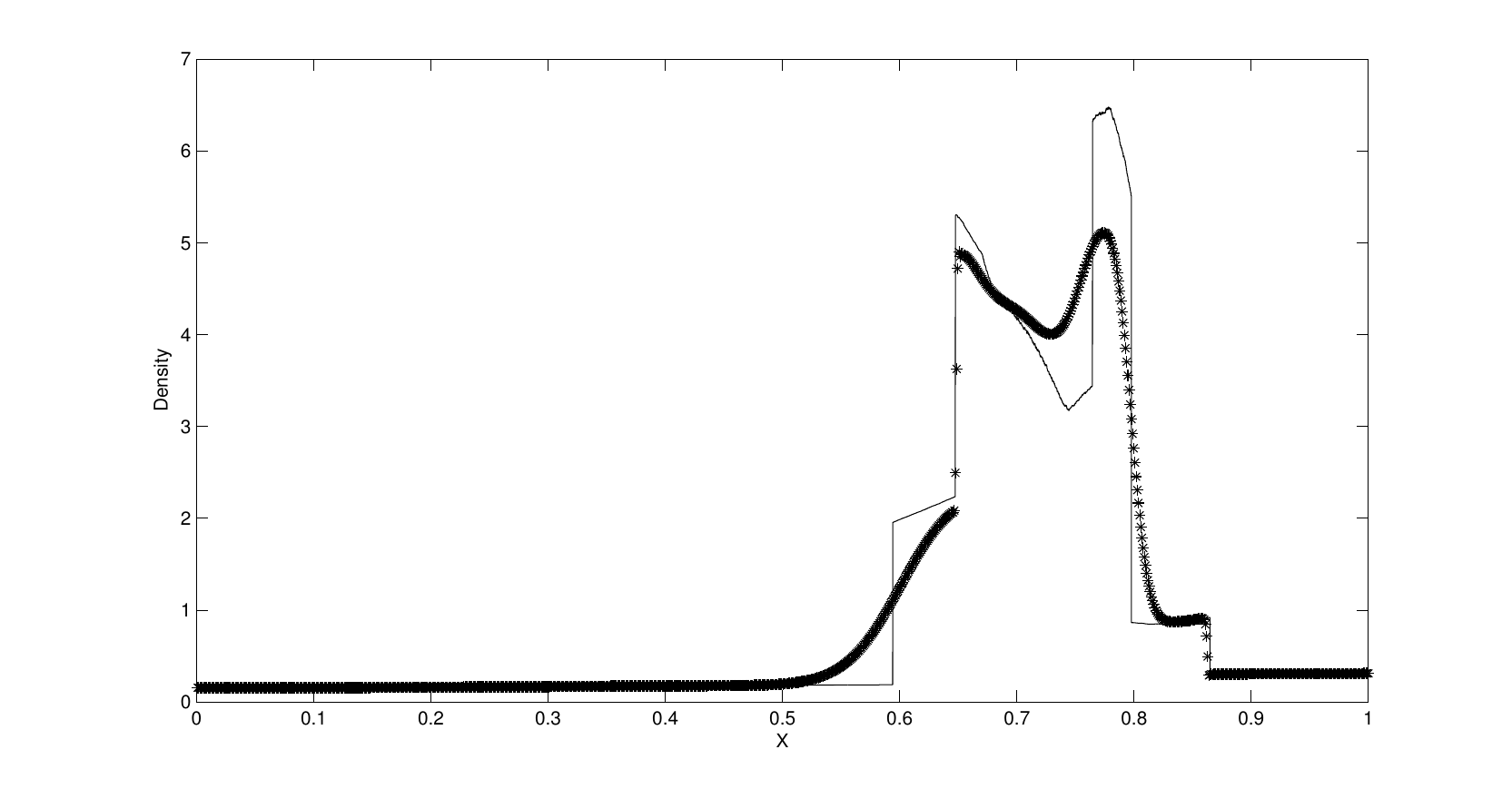}
\caption{Density plot for Blastwave problem with 1000 node points}
\label{fig:BWmrsu}
\end{figure}
Initial conditions are given as $\rho = 1$ and $u =0$ everywhere in the domain. Pressure is given as
\begin{equation*}
p(x,0) = \begin{cases} 1000 & \text{for} \,\, x \in [0,0.1] \\  0.01 & \text{for}  \,\, x \in [0.1,0.9]  \\  100 & \text{for}  \,\, x \in [0.9,1] \end{cases}
\end{equation*}
The final time is 0.038. The solution consists of interaction of strong expansion with strong shock and contact waves. MRSU scheme (represented by *) resolves all the flow features with just 1000 nodes as shown in figure ~\ref{fig:BWmrsu}. The reference solution is shown in solid line. 

\subsection{Dam Break Problem \cite{ZO}}
\begin{figure} [h!] 
\centering
\includegraphics[scale=0.35]{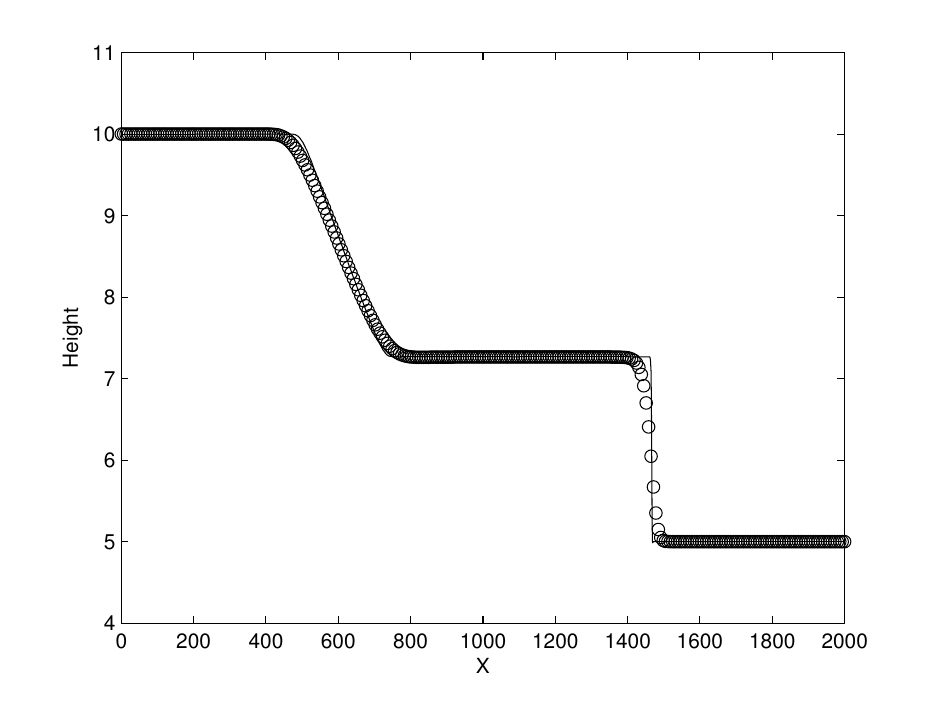}
\includegraphics[scale=0.35]{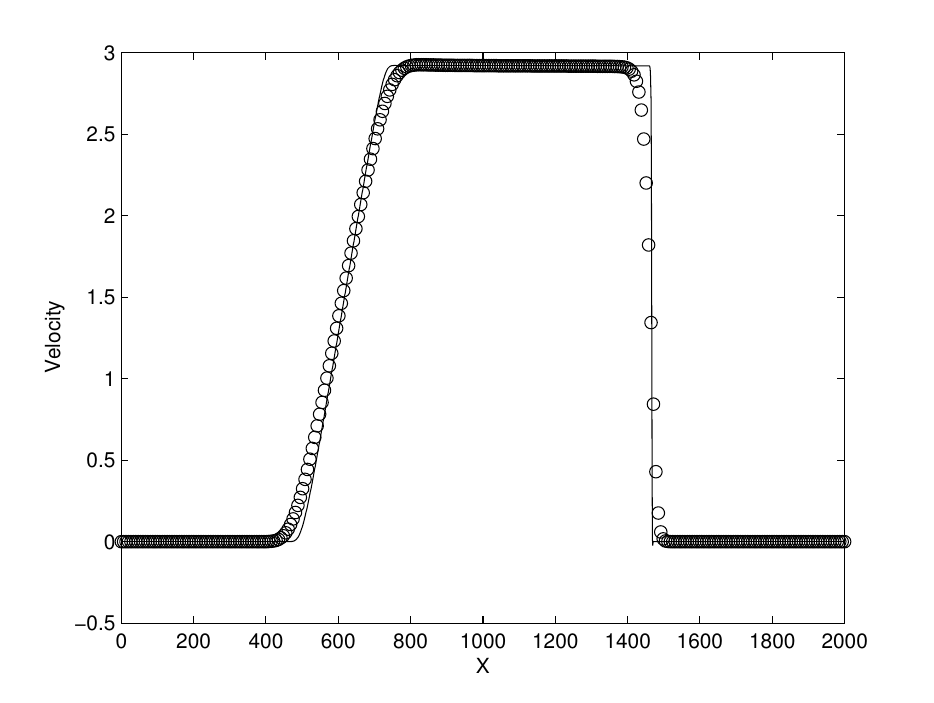}
\caption{Height and velocity plots for dam break problem with 200 nodes.}
\label{fig:DamB1}
\end{figure}

Figure ~\ref{fig:DamB1} shows result of dam break problem. CFL number is 0.25, number of nodes is 200 and final time is 50 seconds. MRSU scheme captures all the flow features like expansion region and hydraulic jump as shown by circles. The exact solution is given by solid line.


\subsection{Strong Shock Problem \cite{ZGM}}
\begin{figure} [h!] 
\centering
\includegraphics[scale=0.35]{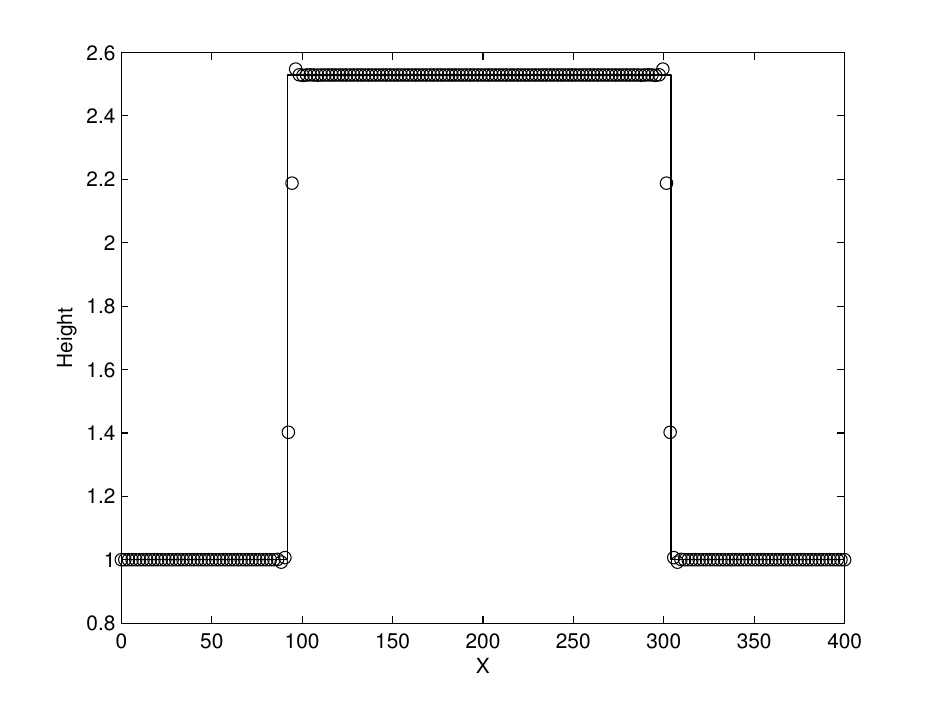}
\includegraphics[scale=0.35]{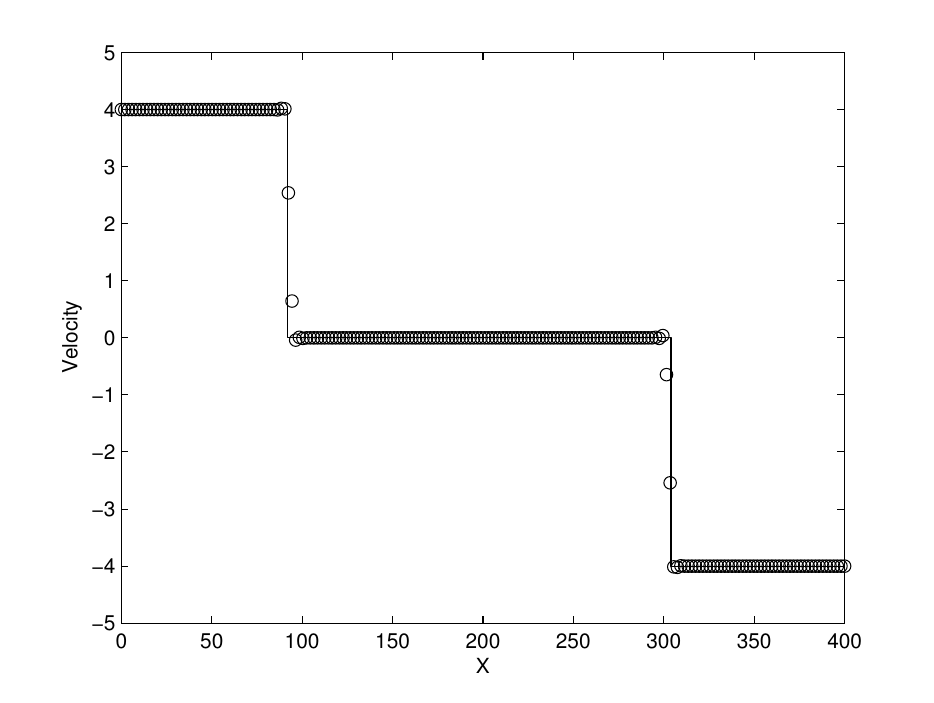}
\caption{Height and velocity plots for strong shock problem with 200 nodes.}
\label{fig:SSwP}
\end{figure}

Figure ~\ref{fig:SSwP} shows the results of strong shock problem. In this case CFL number is 0.4, number of nodes is 200 and final time is 40 seconds. MRSU scheme captures discontinuous wavefronts reasonably well.

\subsection{2D Burgers Steady State Test Cases: Normal and Oblique Discontinuities}
New set of test cases of two-dimensional Burgers equation along with their exact solution is given. The proposed set of test cases can be used to test the accuracy and robustness of numerical algorithms. The domain is $[-0.2,\,1] \times [0,\,\,1]$. Two-dimensional Burgers equation is given by equation 
\begin{equation}\label{2Dbur}
\frac{\partial U}{\partial t} + \frac{\partial \left(\frac{U^2}{2}\right)}{\partial x} + \frac{\partial U}{\partial y} = 0
\end{equation} with
following boundary conditions
$$ U(x,0)  = \textrm{\textbf{a}}  \,\,\, \textrm{for} \, \,-0.2<x<0; \,\, U(x,0)  = \textrm{\textbf{b}}  \,\,\, \textrm{for} \, \,0<x<1; \,\, U(-0.2,y)  = \textrm{\textbf{a}} \,\, \text{and}\,\, U(1,y)  = \textrm{\textbf{b}} \,\, \text{for}\,\, y\in  [0,\,1]$$
\begin{figure} [htpb]
\centering
\includegraphics[scale=1.1]{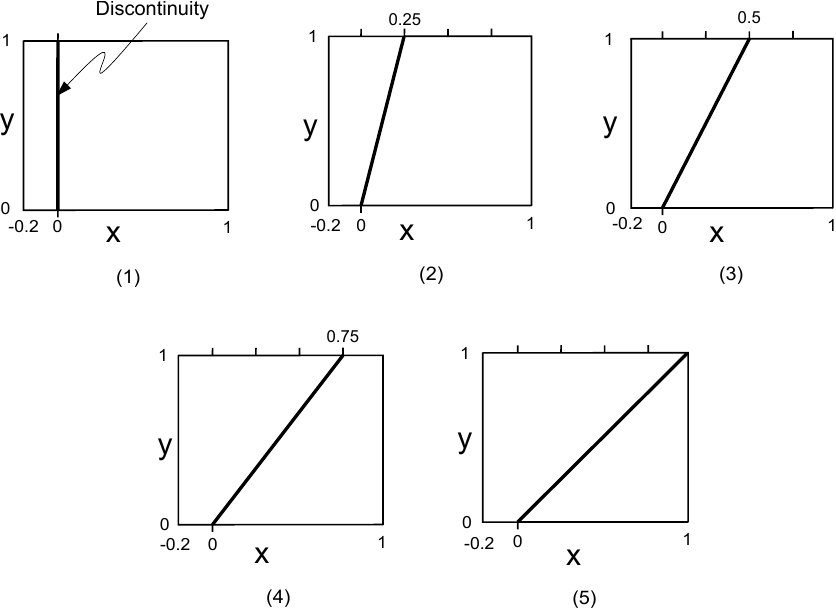}
\caption{2D Burgers exact solution involving normal and oblique discontinuities.}
\label{fig:2dbur2New}
\end{figure}
Now, lets consider different values of $\textbf{a}$ and $\textbf{b}$ given in table ~\ref{Table4} which corresponds to different solution shown in figure \ref{fig:2dbur2New}.
\begin{table}[htpb]
\begin{center}
\small \begin{tabular}{cccccc} \hline
 \textbf{Case} & (1)& (2) &(3) &(4)& (5) \\  \hline 
\textbf{a}  & 1 & 0.5& 1& 1.5 &2  \\
\textbf{b}  & -1 & 0 & 0& 0 & 0   \\
 \hline
 \end{tabular}
\caption{Values of \textbf{a} and \textbf{b} for various cases.}\label{Table4}
\end{center}
\end{table}
\begin{figure} [htpb]
\centering
\includegraphics[scale=0.39]{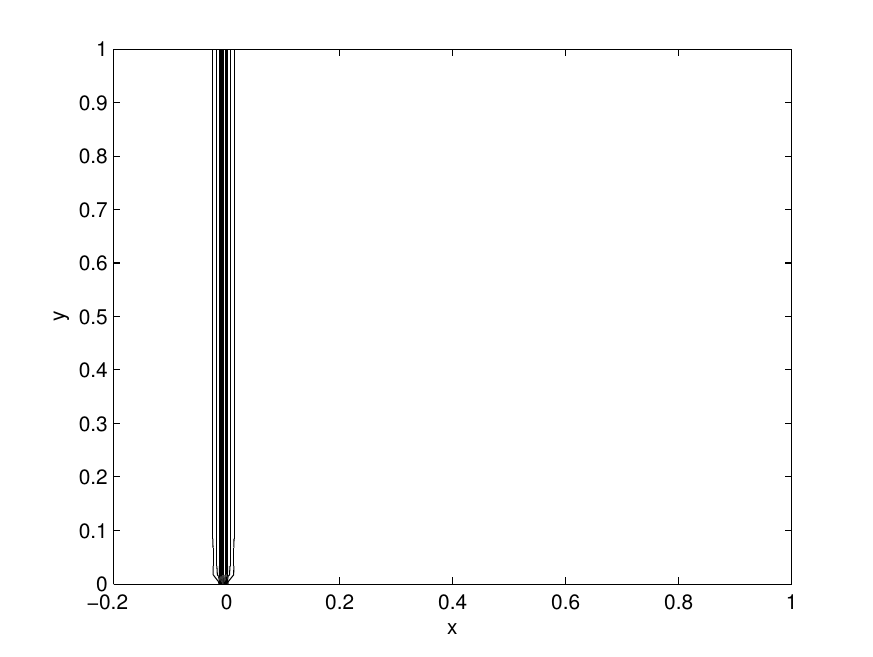}
\includegraphics[scale=0.39]{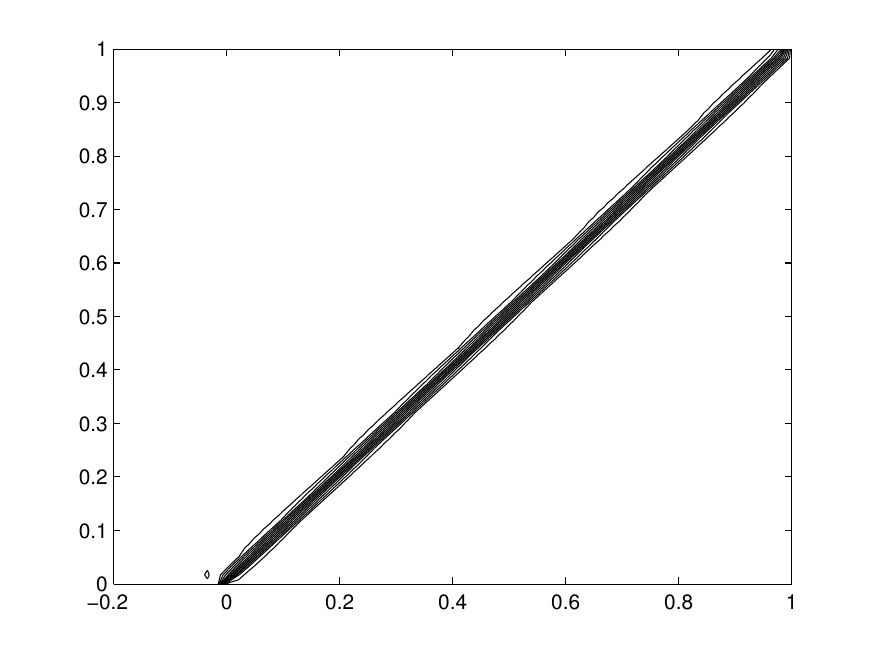}
\caption{Contour plots for normal discontinuity (case (1), left) and oblique discontinuity (case (5), right).}
\label{fig:2dbur2NewC}
\end{figure}
From the given figures it is clear that cases (2)-(5) produce oblique discontinuity with decreasing discontinuity angle with respect to horizontal line whereas case (1) produces normal discontinuity. The exact solution is 
\begin{equation*}
U(x,y) = \begin{cases} \textrm{\textbf{a}} & \textrm{If} \, x < \frac{\textrm{\textbf{a}}+\textrm{\textbf{b}}}{2}\, y \\ \textrm{\textbf{b}}  & \textrm{If} \,x \geq \frac{\textrm{\textbf{a}}+\textrm{\textbf{b}}}{2}\, y  \end{cases}
\end{equation*}
$72 \times 60$ quadrilateral mesh is used to solve case (1) and (5) as shown in figure \ref{fig:2dbur2NewC}. Both normal and oblique discontinuities are captured well by the proposed scheme.


\subsection{2D Burgers Test Case: Oblique Discontinuity}
\begin{figure} [htpb] 
\centering
\includegraphics[scale=0.37]{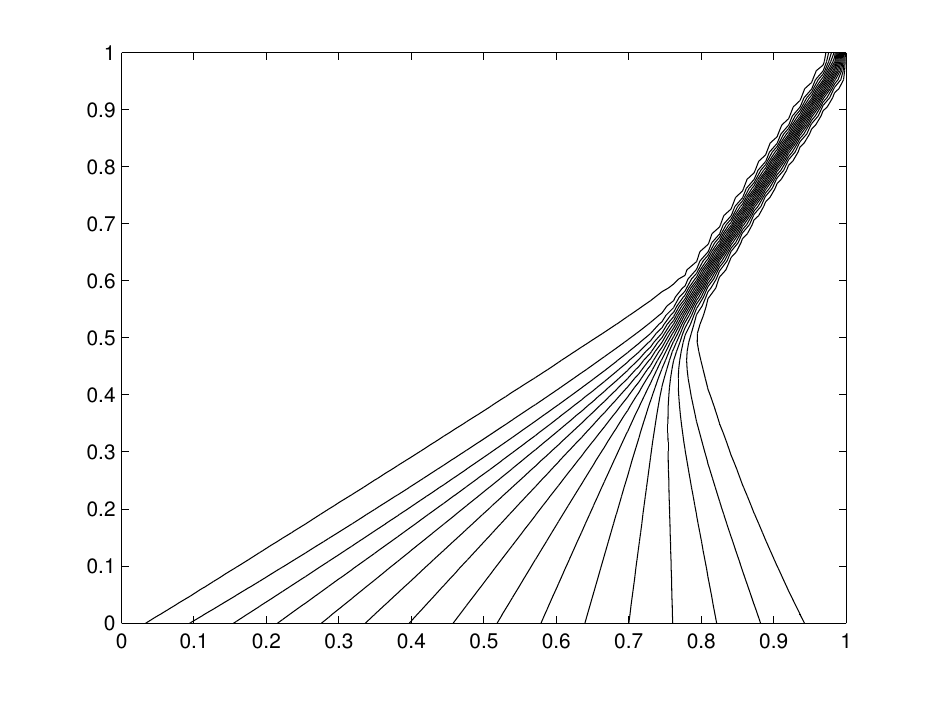}
\includegraphics[scale=0.37]{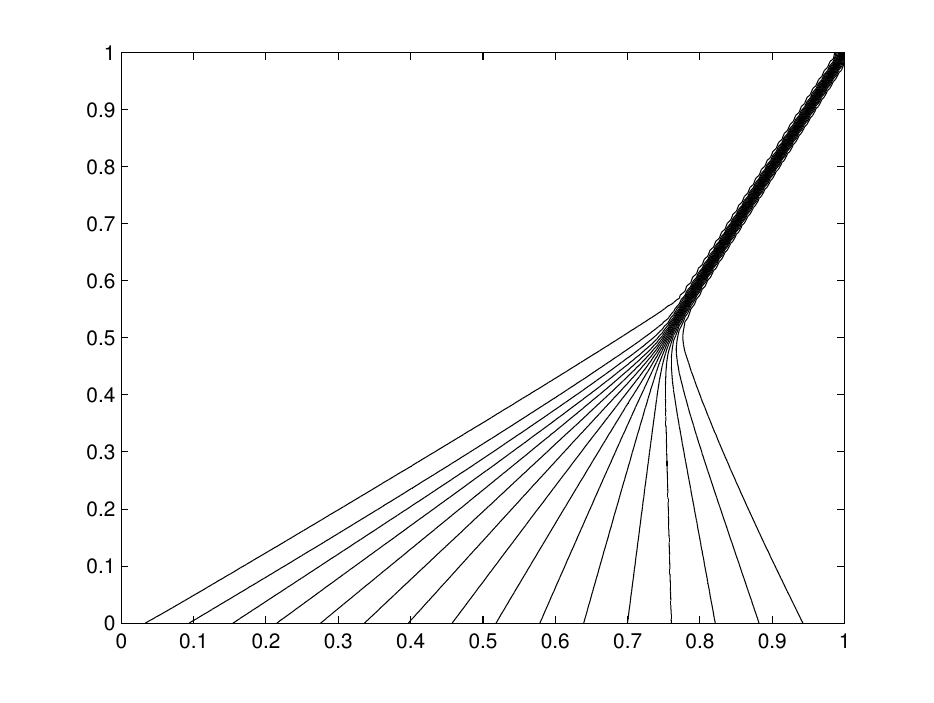}
\caption{$64 \times 64$ (left) and $128 \times 128$ (right) quadrilateral mesh.}
\label{fig:BT2}
\end{figure}

Two-dimensional Burgers equation is given by equation \eqref{2Dbur}. The domain is $[0,\,1]^2$ and boundary conditions are
$$U(0,y) = 1.5 \,\,\text{and}\,\, U(1,y) = -0.5\,\, \text{for} \,\,0<y<1; \,\,U(x,0) = 1.5-2x \ \ \text{for} \,\, 0<x<1$$
Exact solution is given in \cite{Spekreijse}. The oblique discontinuity is captured quite accurately as shown in figure ~\ref{fig:BT2} with different grid size.

\subsection{2D Burgers Test Case with Non-convex Flux Function - KPP Rotating Wave \cite{AKe}}
The domain is $[-2,\,2]\times [-2.5,\,1.5]$. Two-dimensional scalar conservation law with non-convex flux function is 
\begin{equation*}
\frac{\partial U}{\partial t} + \frac{\partial \sin U}{\partial x} + \frac{\partial \cos U}{\partial y} = 0
\end{equation*}
with initial conditions as
\begin{figure} [h!] 
\centering
\includegraphics[scale=0.4]{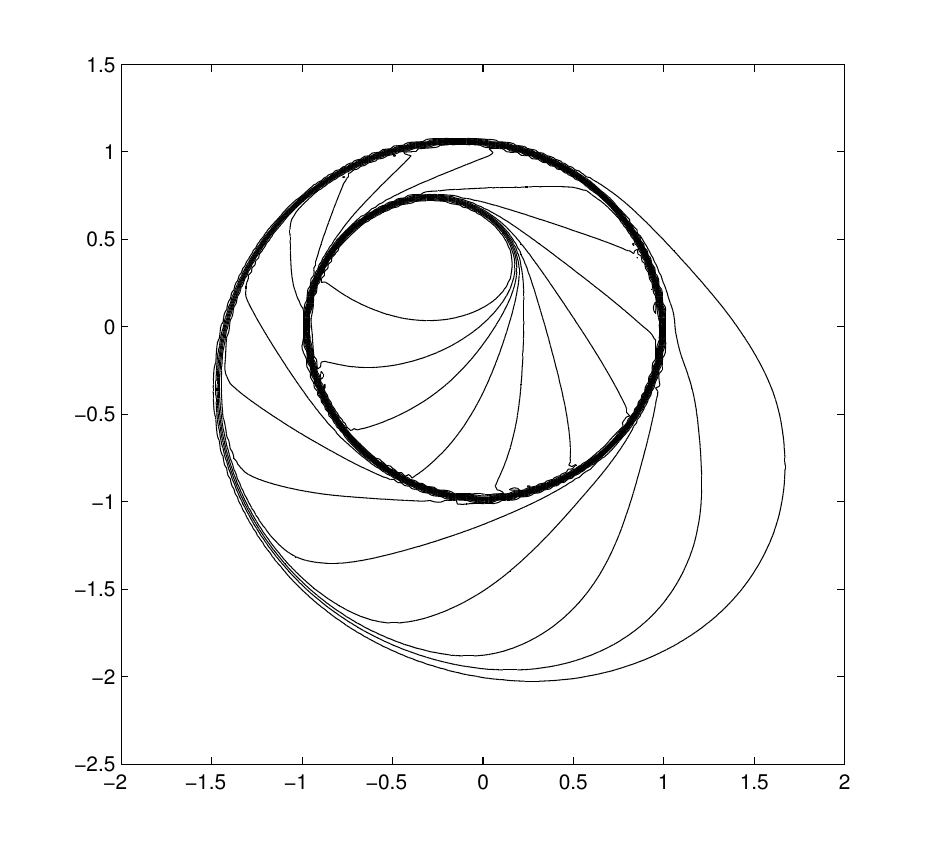}
\caption{KPP rotating wave}
\label{fig:KPPu4}
\end{figure}
\begin{equation*}
U(x,y,0) = \begin{cases} 3.5\pi & \textrm{If} \,\, x^2+y^2 < 1 \\   \frac{\pi}{4} &  \textrm{otherwise}\end{cases}
\end{equation*}
Figure ~\ref{fig:KPPu4} shows the contour plots on mesh $\Delta x = \Delta y = 1/50$. The rotating discontinuity is captured accurately by the proposed scheme.

\subsection{Shock Reflection Test Case \cite{Yee_Warming_Harten}}
In this test case the domain is rectangular $[0 ,\,\,3] \times [0,\,\,1]$. The boundary conditions are, inflow (left boundary) : $\rho =1, u_1=2.9, u_2=0, p=1/1.4$. Post shock condition (top boundary) : $\rho =1.69997, u_1=2.61934,u_2=-0.50633, p=1.52819$. Bottom boundary is a solid wall where slip boundary condition is applied, \textit{i.e.}, $\mathbf{u} \cdot n = 0$. At right boundary where flow is supersonic all primitive variables $\rho$, $u_1 $, $u_2$ and $p$ are extrapolated.
\begin{figure} [h!] 
\centering
\includegraphics[trim=1cm 0.9cm 1cm 12cm, clip=true, scale=0.35, angle = 0]{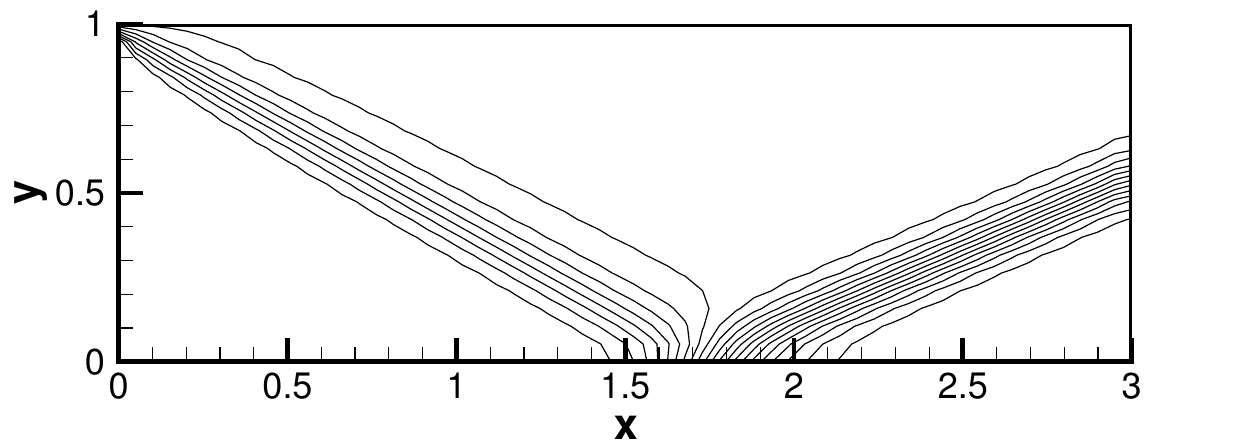}
\includegraphics[trim=1cm 0.9cm 1cm 12cm, clip=true, scale=0.35, angle = 0]{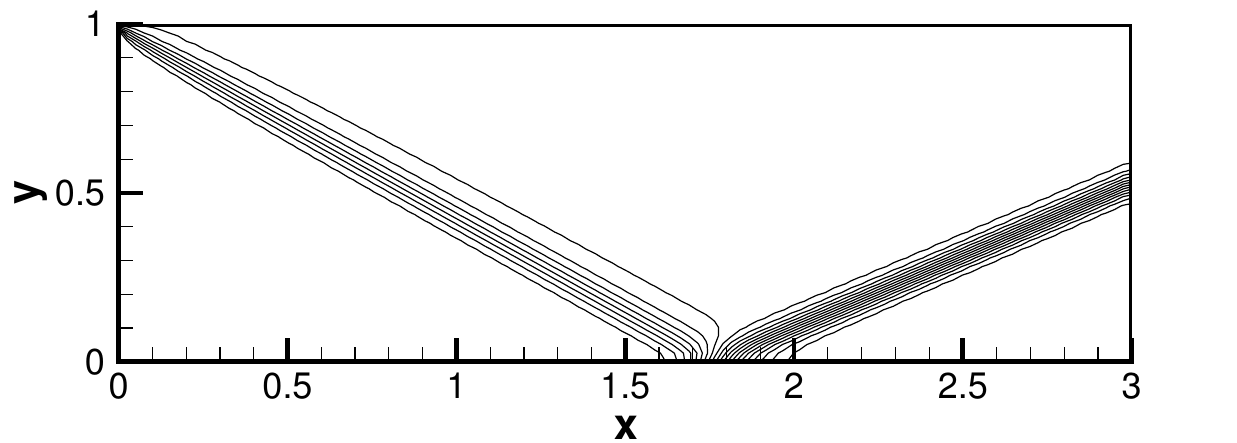}
\includegraphics[trim=1cm 0.9cm 1cm 12cm, clip=true, scale=0.35, angle = 0]{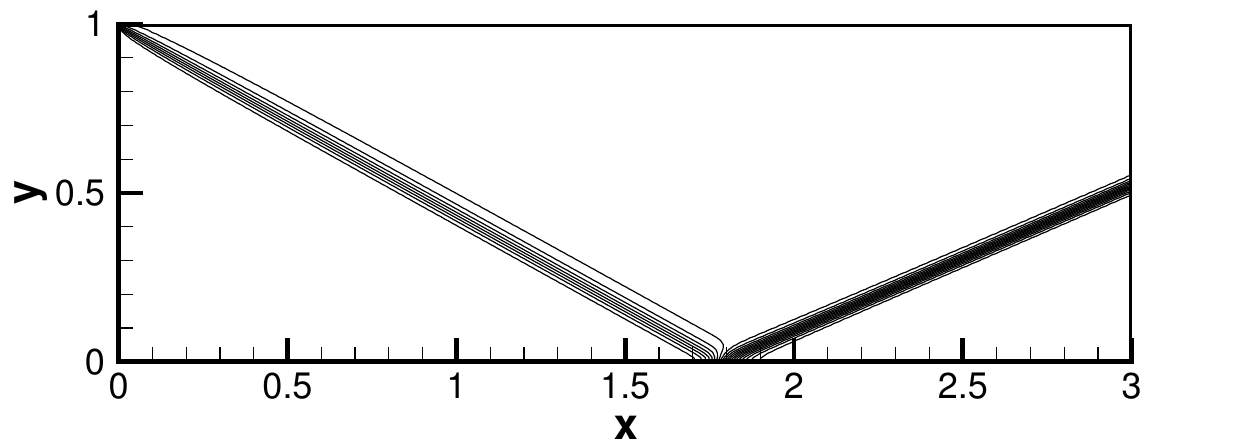}
\includegraphics[scale=0.4]{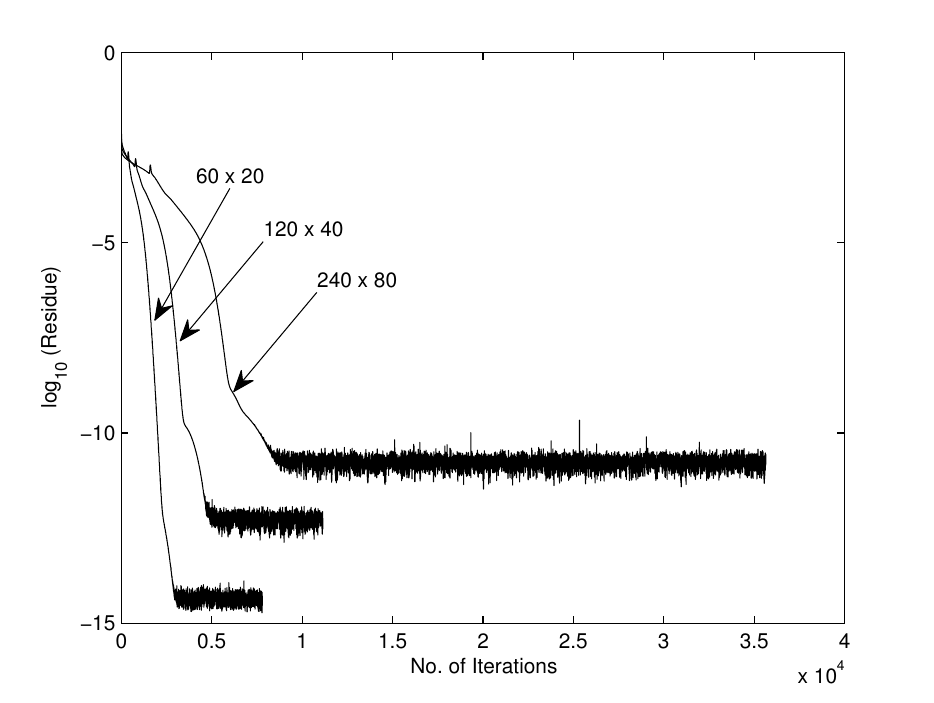}
\caption{Pressure contours (0.8:0.1:2.8) for $60 \times 20,\,\,\, 120 \times 40$  and $240 \times 80$ quadrilateral mesh along with residue plots.}
\label{fig:MRSu4}
\end{figure}
Pressure plots for  $60 \times 20$, $120 \times 40$ and $240 \times 80$ quadrilateral mesh along with residue plots are given in figure \ref{fig:MRSu4}.
\begin{table}[htpb]
\begin{center}
\small \begin{tabular}{cccc} \hline
  & 60 $\times$ 20 & 120 $\times$ 40 & 240 $\times$ 80 \\  \hline 
MRSU scheme  & 186  & 624 & 3032   \\
MSU scheme & 1335 & 5019 &  24129 \\
 \hline
 \end{tabular}
\caption{Computational cost in seconds.}\label{Table5}
\end{center}
\end{table}
Table \ref{Table5} gives computational cost comparison for relaxation system based MRSU scheme and MSU scheme which does not use relaxation system (see Appendix C) over different grid size. All the simulations are run on 3.10 GHz desktop machine till steady-state is achieved. From these values one can see that, MSU scheme is more expensive than that of MRSU scheme due to evaluation of Jacobian matrices at every time step.


\subsection{Half Cylinder Test Case}
\begin{figure} [h!] 
\centering
\includegraphics[trim=1cm 0.9cm 11cm 1cm, clip=true, scale=0.36, angle = 0]{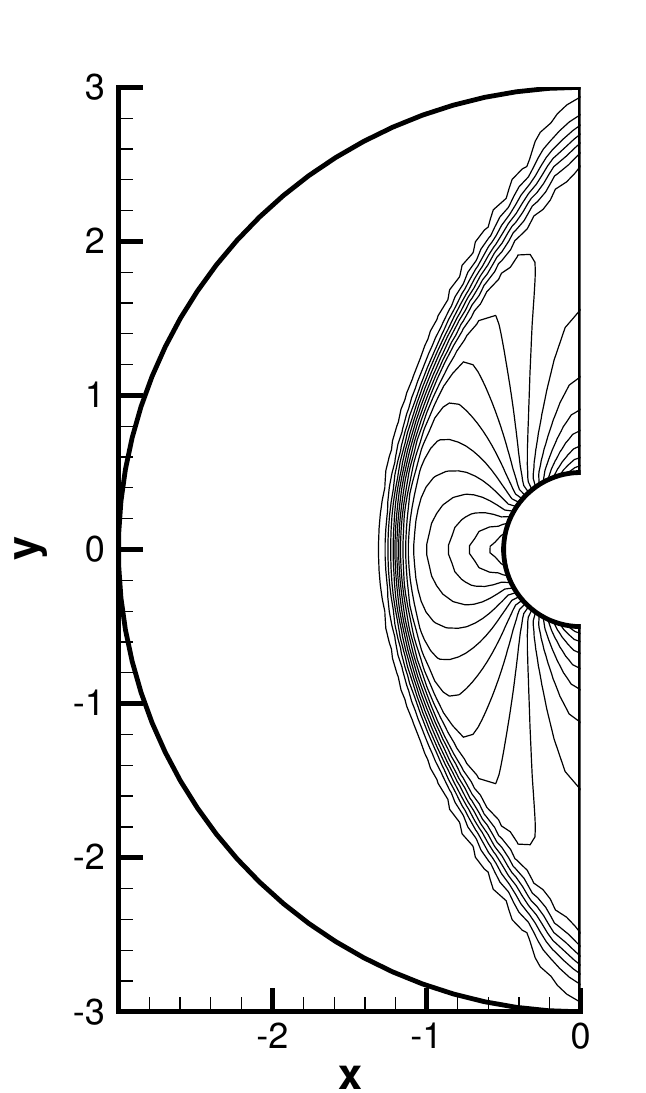}
\includegraphics[trim=1cm 0.9cm 11cm 1cm, clip=true, scale=0.36, angle = 0]{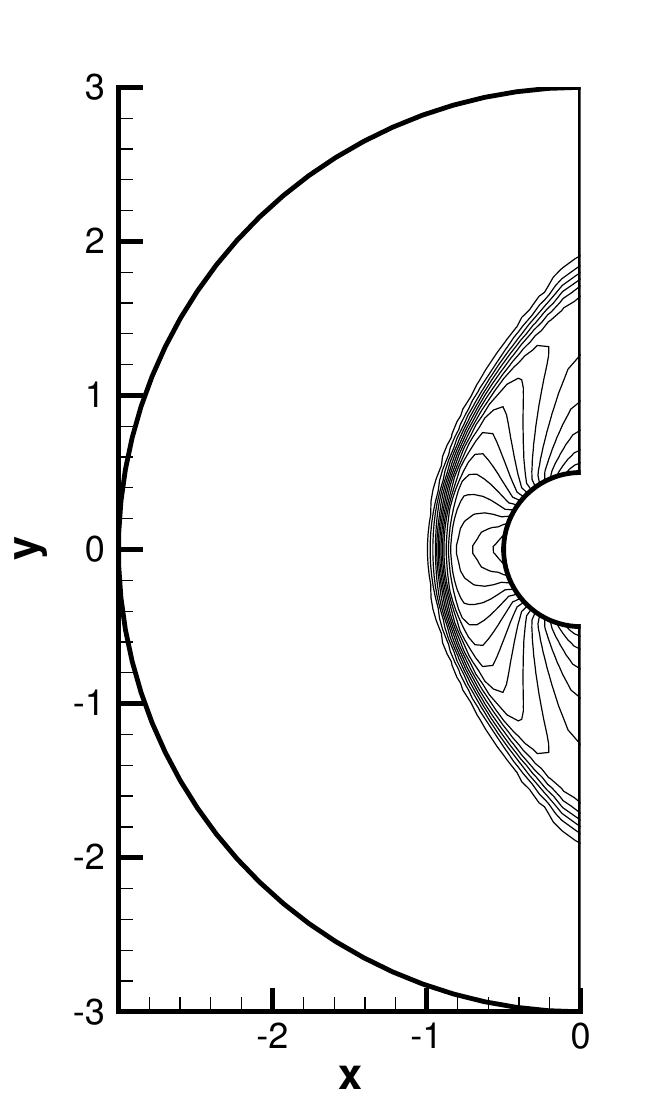}
\includegraphics[trim=1cm 0.9cm 11cm 1cm, clip=true, scale=0.36, angle = 0]{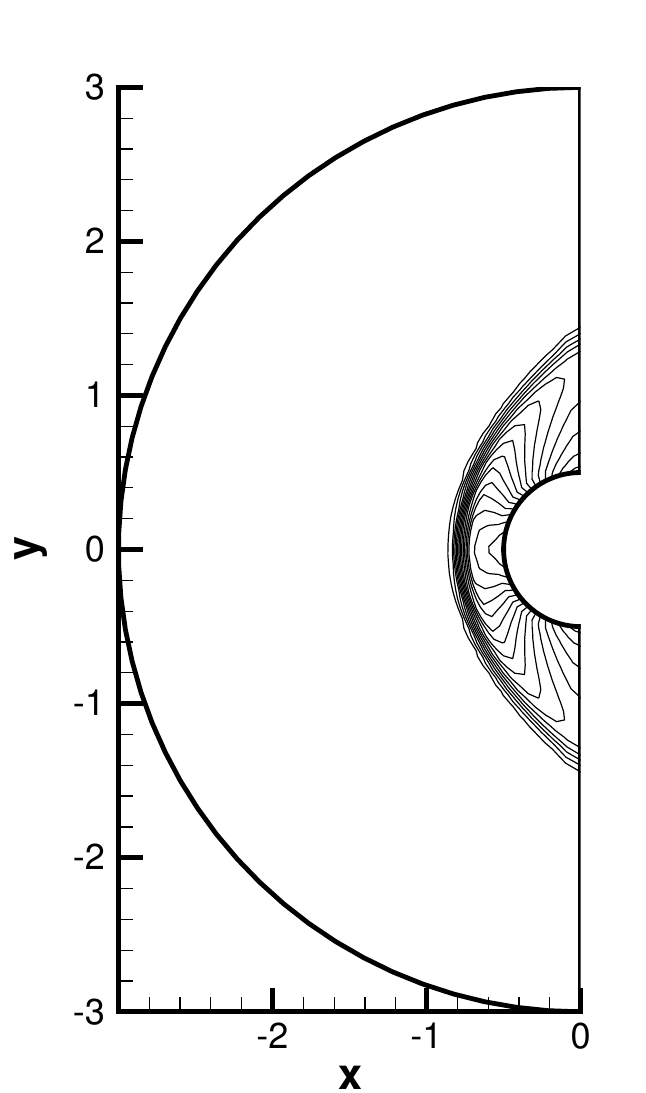}
\includegraphics[trim=1cm 0.9cm 11cm 1cm, clip=true, scale=0.36, angle = 0]{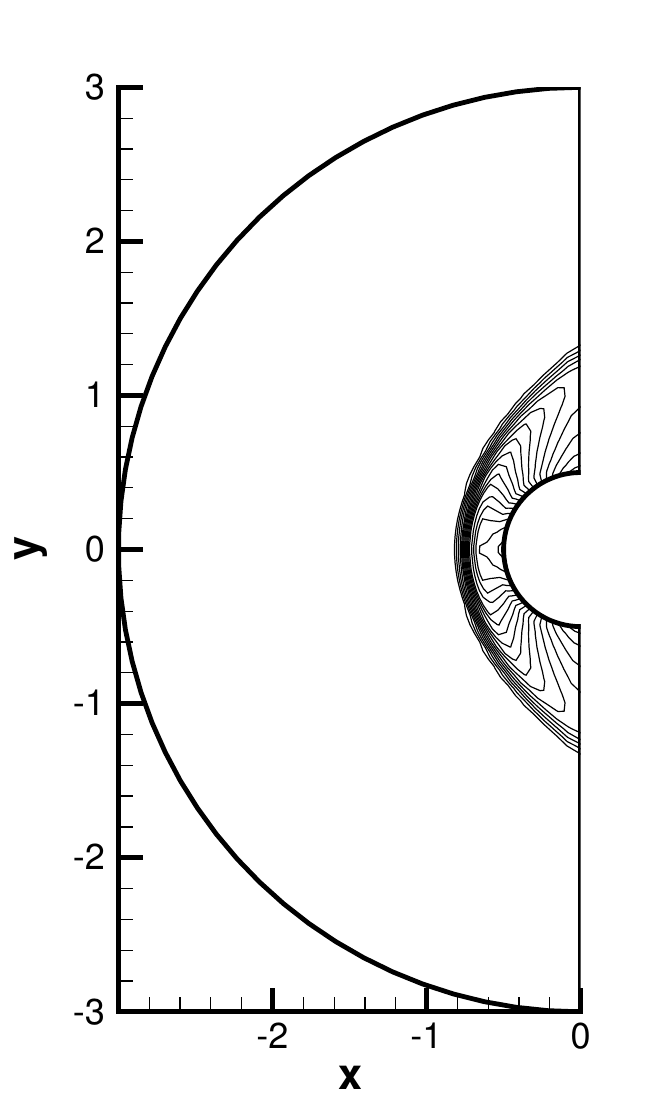}
\caption{Pressure contours for Mach 2 (0.8:0.2:4), Mach 3 (1:0.5:16), Mach 6 (2:2:32) and Mach 20 (20:20:360) using $45\times 45$ quadrilateral mesh.}%
\label{fig:KHcy4}
\end{figure}

Four test cases of inflow Mach number 2, 3, 6 and 20 are tested on a half cylinder \cite{HVi}. The domain is half circular. Left outer circle is inflow boundary, small circle inside the domain is a cylinder wall and the straight edges on right side are supersonic outflow boundaries (see figure \ref{fig:KHcy4}).
Pressure plots show bow shock in front of the half cylinder which is captured accurately at the right position in each case. These results are compared with existing results \cite{FSBillig}.

\subsection{Double Mach Reflection Test Case}
In the initial condition, Mach 10 shock wave makes an angle of $60^o$ with the reflecting wall. The undisturbed air in front of shock has density 1.4 and pressure 1. Initial conditions and boundary conditions are given in \cite{WCPP}.
\begin{figure} [h!] 
\centering
\includegraphics[trim=1cm 0.9cm 1cm 12cm, clip=true, scale=0.5, angle = 0]{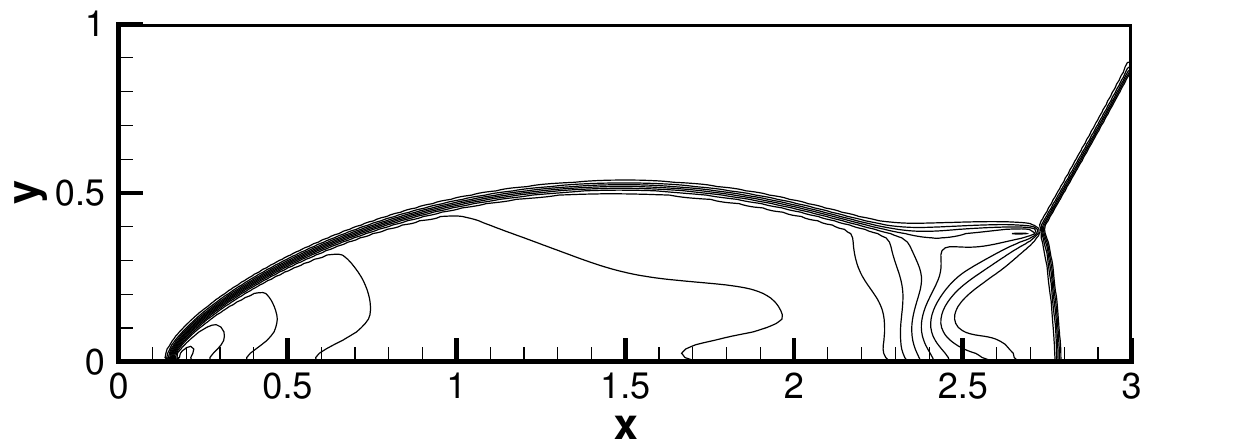}
\caption{Density contours (3:0.92:26) using $240 \times 80$ quadrilateral mesh.}
\label{fig:DMRmrsu}
\end{figure}
As the shock moves on the ramp, a self similar shock structure with two triple points evolve. In the primary triple point, Mach stem, incident shock and reflected shock meet as shown in figure \ref{fig:DMRmrsu}. An analytical work on double Mach reflection can be found in \cite{HLi, GaB} and the references therein, whereas experimental results are given in \cite{LGg}.


\subsection{Parallel Jet Flow \cite{DHP}}
The domain is $[0,\,1]^2$ and the initial conditions are given as
\begin{align*}
& M = 4, \,\,  \rho = 0.5 \,\, p = 0.25   \,\,\, \textrm{if}\,\, y >0.5 
\\ & M = 2.4, \,\,  \rho =1 \,\, p =1   \,\,\, \textrm{if} \,\,y <0.5 
\end{align*}
\begin{figure} [h!] 
\centering
\includegraphics[scale=0.25]{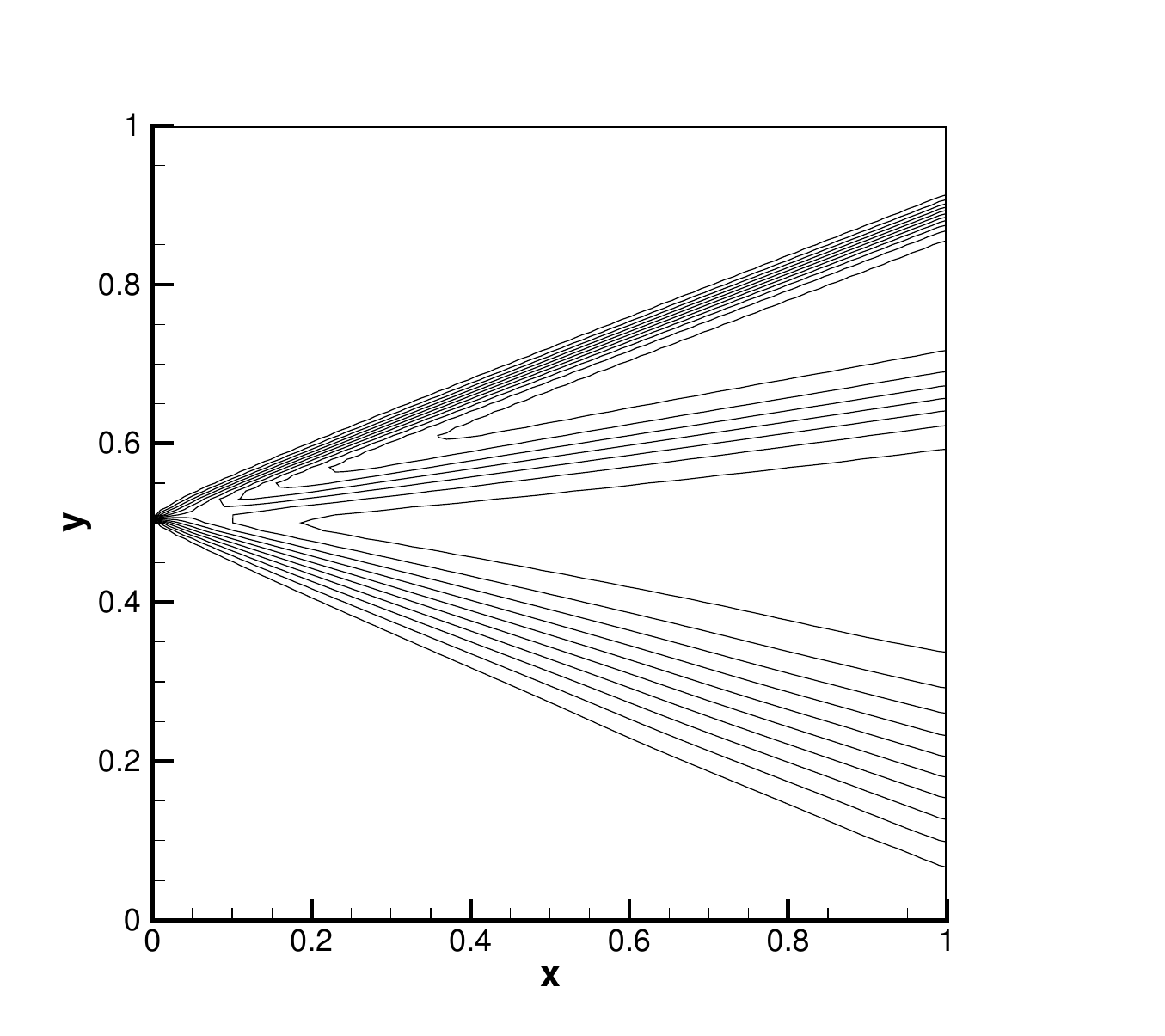}
\includegraphics[scale=0.25]{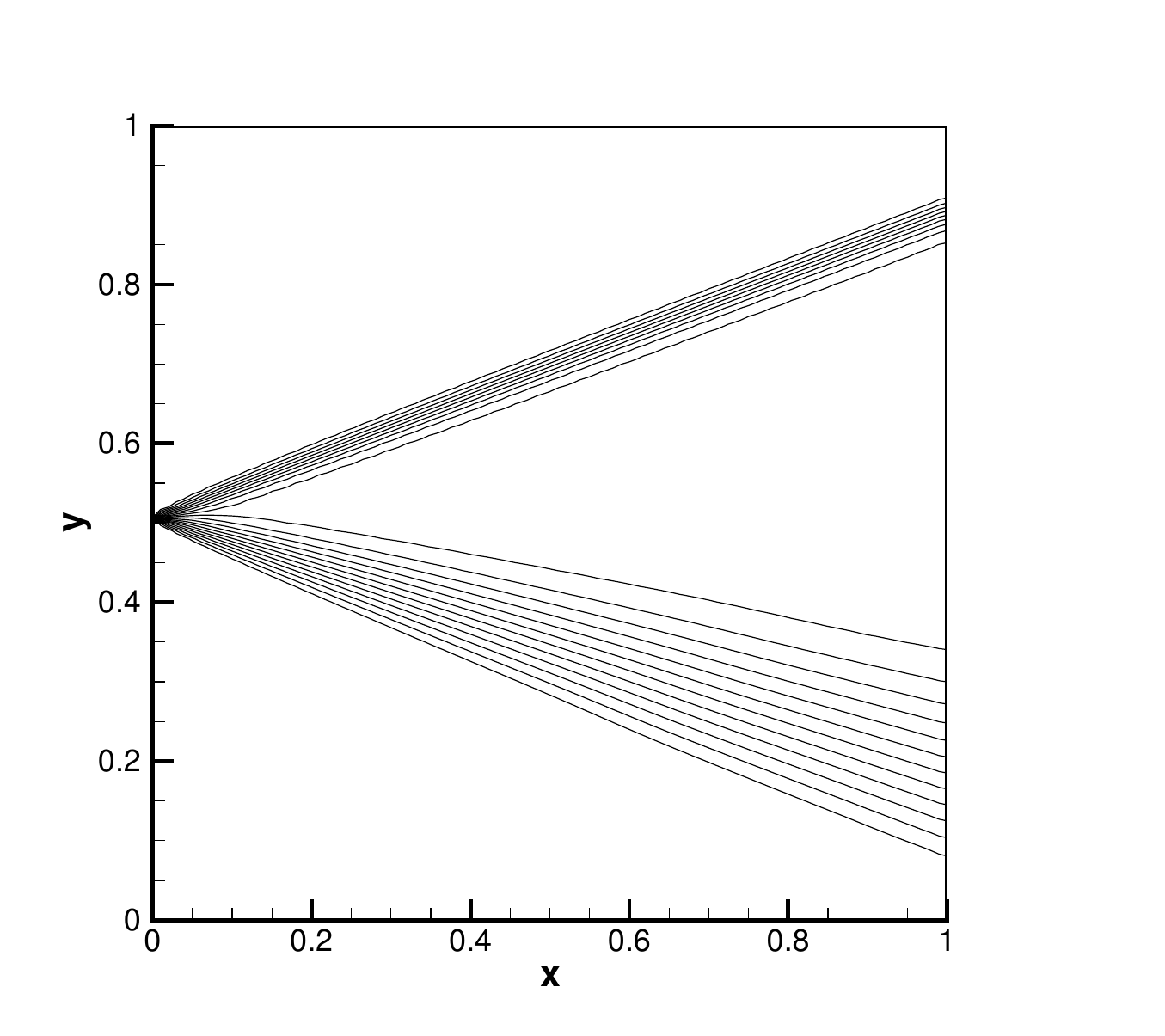}
\caption{Density contours (0.55:0.03:0.94) (left) and pressure contours (0.3:0.03:0.9) (right) for $100 \times 100$ quadrilateral mesh.}
\label{fig:JetMRSU}
\end{figure}
where $M$ is Mach number. Left boundary is the supersonic inflow and top, right and bottom boundaries are supersonic outflows where all the variables are extrapolated. Figure  \ref{fig:JetMRSU} shows the density and pressure plots $100 \times  100$ mesh. All the flow features like shock, expansion and contact are captured using MRSU scheme.



\subsection{Circular Explosion Problem\cite{TOro}}
In this test case the domain is $[-1,\,1]^2$. The initial condition represents two regions, first region is inside the circle with radius 0.4 and second region is outside the circle. The initial conditions are 
\begin{equation*}
(\rho, u_1,u_2,p)(x,y,0) = 
\begin{cases}
(1,\,\, 0,\,\,0,\,\,1) & \,\, \textrm{If}\,\, |r|\leq 0.4 \\ (0.125,\,\, 0,\,\,0,\,\,0.1) & \,\, \textrm{Otherwise}
\end{cases}
\end{equation*}
The final time is 0.2. The solution has circular shock which is constantly moving away from the centre, circular expansion fan which is moving towards the centre and contact discontinuity which separates shock and expansion waves. 
\begin{figure} [h!] 
\centering
\includegraphics[scale=0.28]{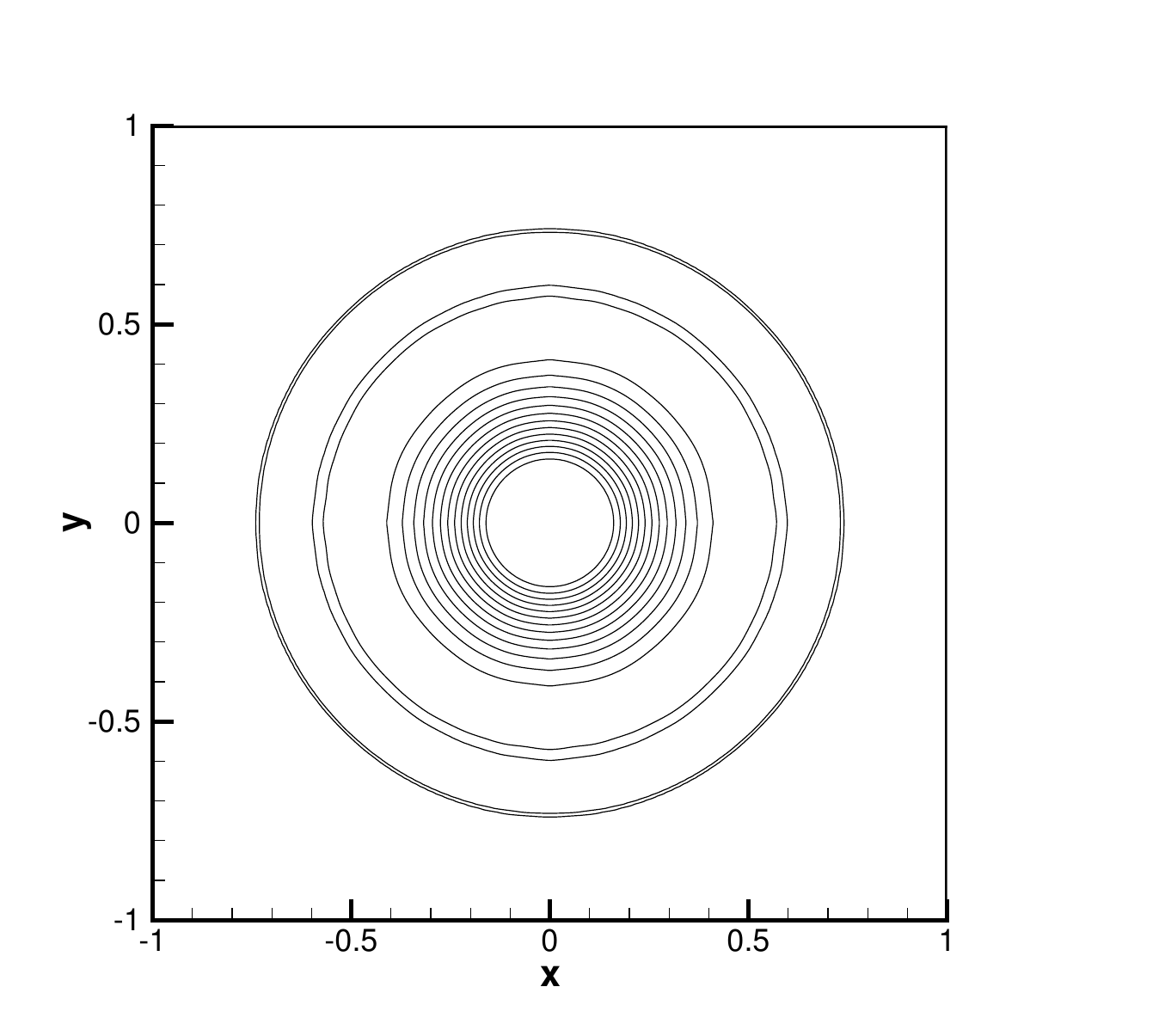}
\caption{Density contours (0.22:0.042:0.86) 
using $200 \times 200$ mesh. }
\label{fig:EPmrsu}
\end{figure}
Figure \ref{fig:EPmrsu} shows the density contours where shock, contact and expansion waves can be seen.

\subsection{Flow Over a Bump}
This is one of the difficult test case due to presence of stagnation point at the front and rear end of the bump. The bump height is 4\% of the chord length. For the numerical simulation three different flow fields are considered, namely, Mach 1.4 (supersonic flow), Mach 0.85 (transonic flow) and Mach 0.5 (subsonic flow) over the bump. 

\subsubsection{Supersonic Flow with M = 1.4 :}
\begin{figure} [h!] 
\centering
\includegraphics[trim=1cm 0.9cm 1cm 12cm, clip=true, scale=0.5, angle = 0]{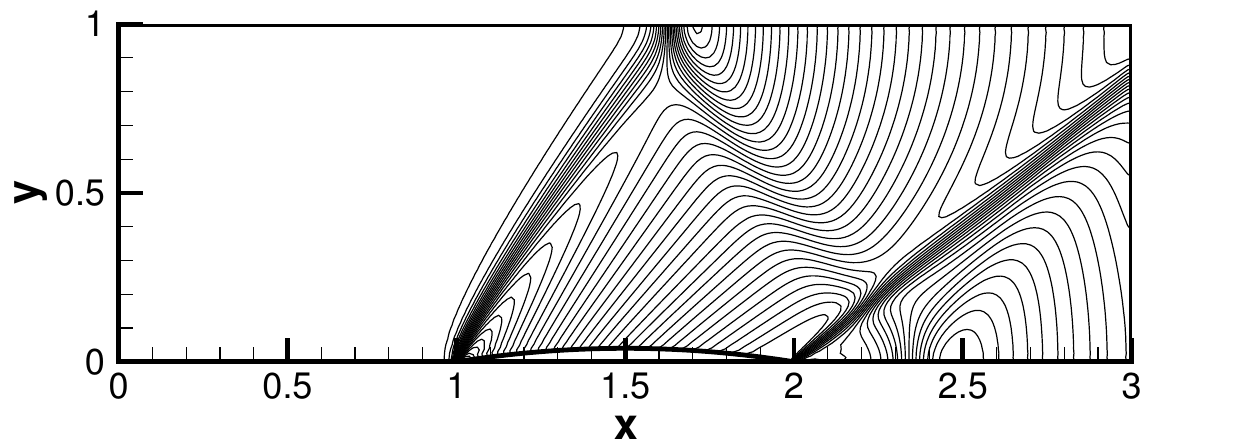}
\caption{Mach 1.4 pressure contours (0.5:0.02:1.4)  with $240 \times 80 $ quadrilateral mesh.}
\label{fig:SuperBumpMRSU}
\end{figure}

Figure ~\ref{fig:SuperBumpMRSU} shows the pressure contours for Mach 1.4.
In this test case shock appears from the leading edge of the bump which hits the top wall and then reflects back on bottom wall. The reflected shock interacts with the shock generated from the trailing edge of the bump and it reflects again from the bottom boundary. MRSU captures all the essential flow features accurately.

\subsubsection{Transonic Flow with M = 0.85 :}
\begin{figure} [h!] 
\centering
\includegraphics[trim=1cm 0.9cm 0.5cm 11cm, clip=true, scale=0.5, angle = 0]{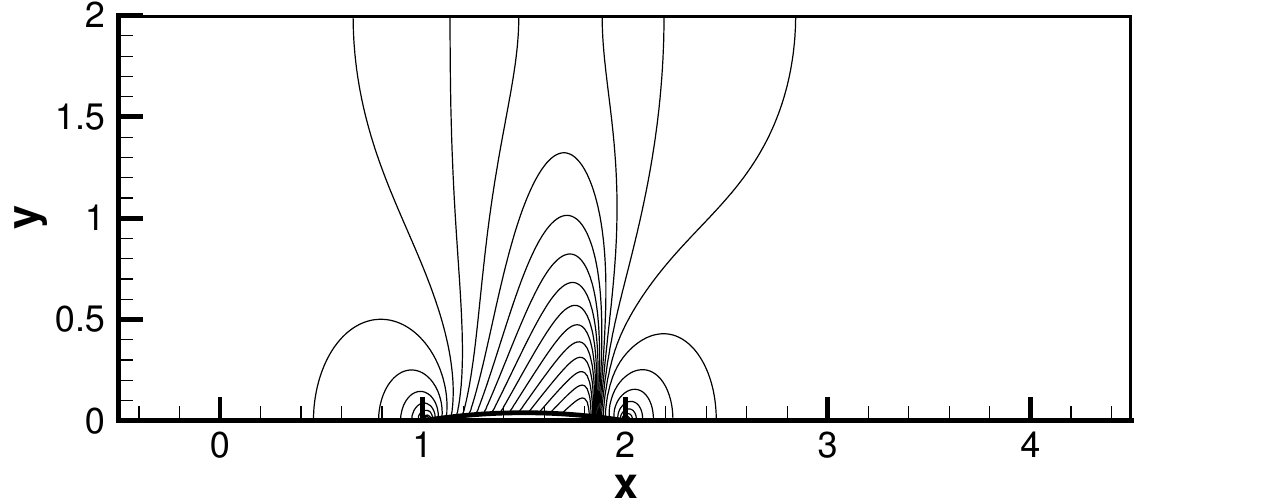}
\includegraphics[scale=0.38]{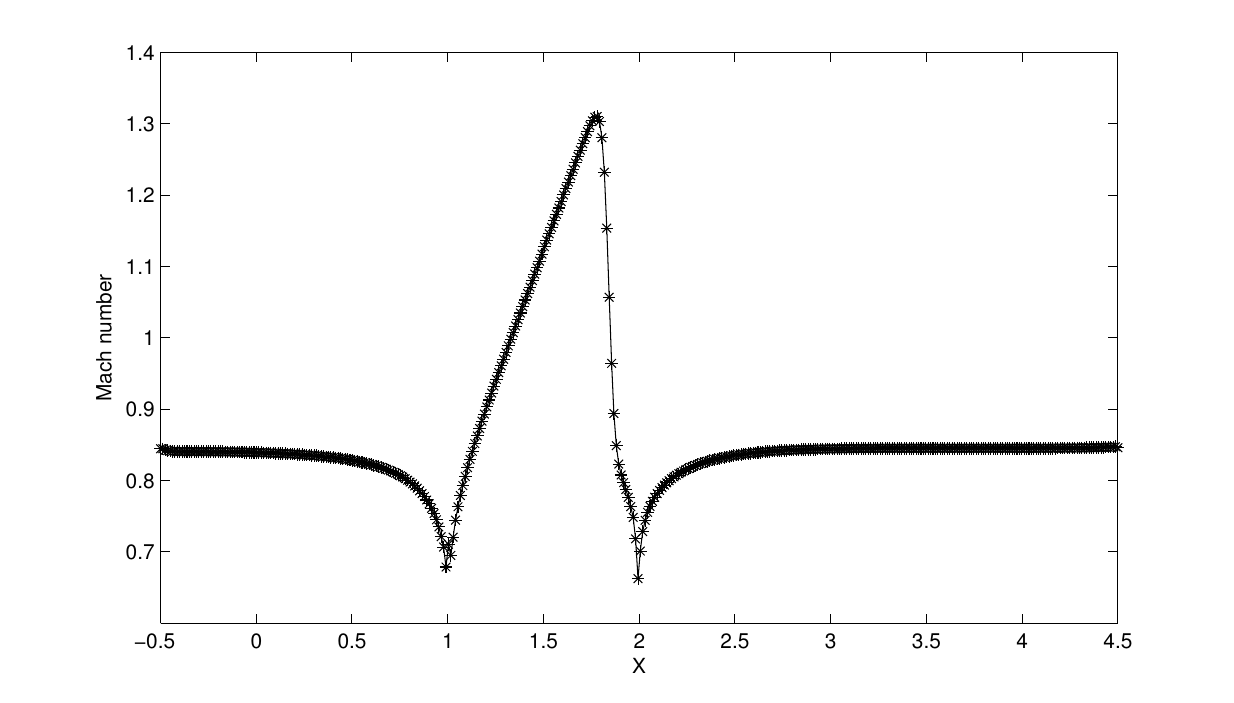}
\caption{Mach 0.85 pressure contours (0.45:0.02:0.89)  with $400 \times 160 $ quadrilateral mesh and Mach number variation along the bottom wall.}
\label{fig:TransBumpMRSU}
\end{figure}
Figure \ref{fig:TransBumpMRSU} shows the pressure contours for Mach 0.85 and the variation of Mach number along the bottom wall. 
The shock appears approximately at 86\% of the bump from the front end with upstream Mach number approximately 1.3.

\subsubsection{Subsonic Flow with M = 0.5 :}
\begin{figure} [h!] 
\centering
\includegraphics[trim=1cm 0.9cm 1cm 12cm, clip=true, scale=0.5, angle = 0]{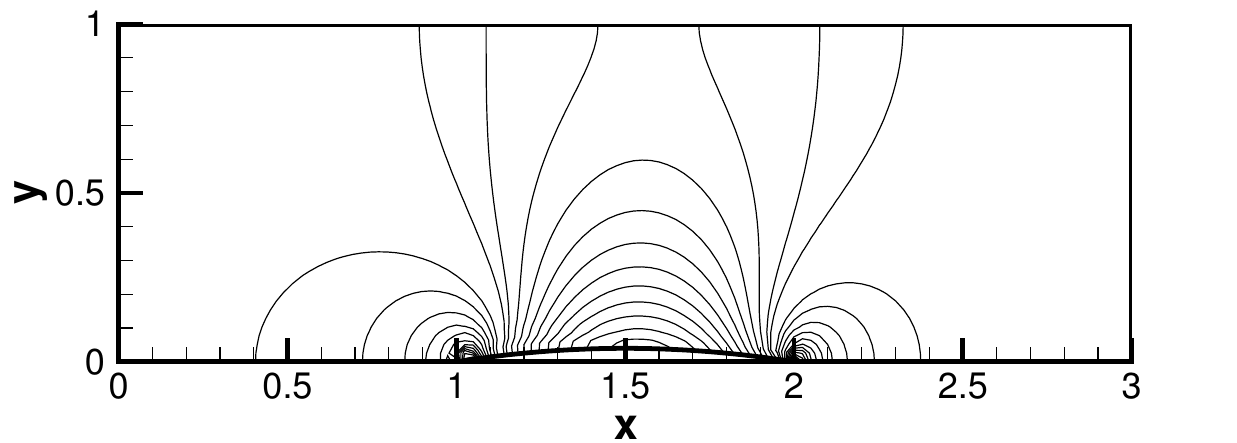}
\caption{Mach 0.5 pressure contours (0.68:0.003:0.775)  with $120 \times 40 $ quadrilateral mesh.} 
\label{fig:SubBumpMRSU}
\end{figure}
In this test case no shock wave appear. Figure \ref{fig:SubBumpMRSU} shows the pressure contours for Mach 0.5. For both transonic as well as subsonic flow over the bump, Riemann invariant based boundary conditions \cite{ChH} are used.

\subsection{Supersonic Flow Over a Reverse Bump}
This is a new test case introduced here by reversing the bump. The inlet Mach number is 1.4.
Expansion waves originate from the front end of the reverse bump. Curved surface of the reverse bump compresses the flow isentropically which generates Mach waves. These Mach waves coalesce to form an oblique shock at an angle of $50^o$ approximately with the horizontal wall. This oblique shock hits the inviscid top wall where slip boundary condition is applied. 
\begin{figure} [h!] 
\centering
\includegraphics[scale=0.65]{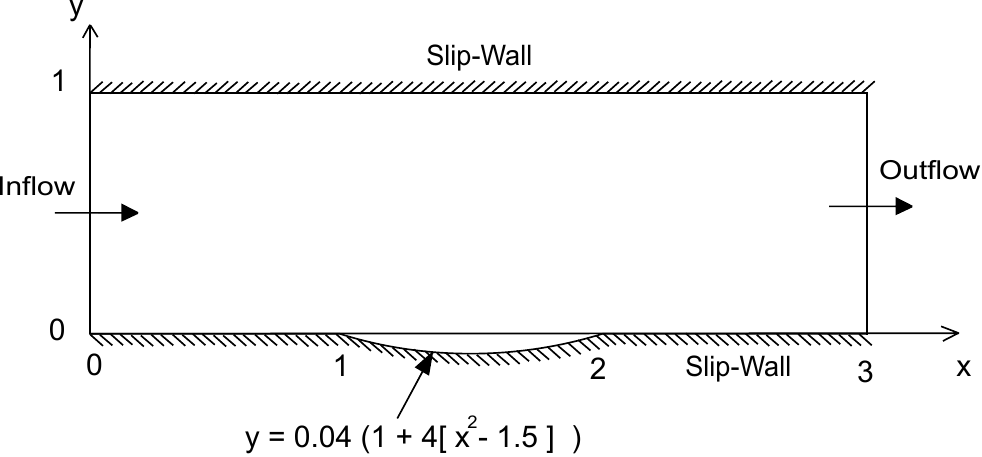}
\includegraphics[trim=1cm 0.9cm 1cm 12cm, clip=true, scale=0.5, angle = 0]{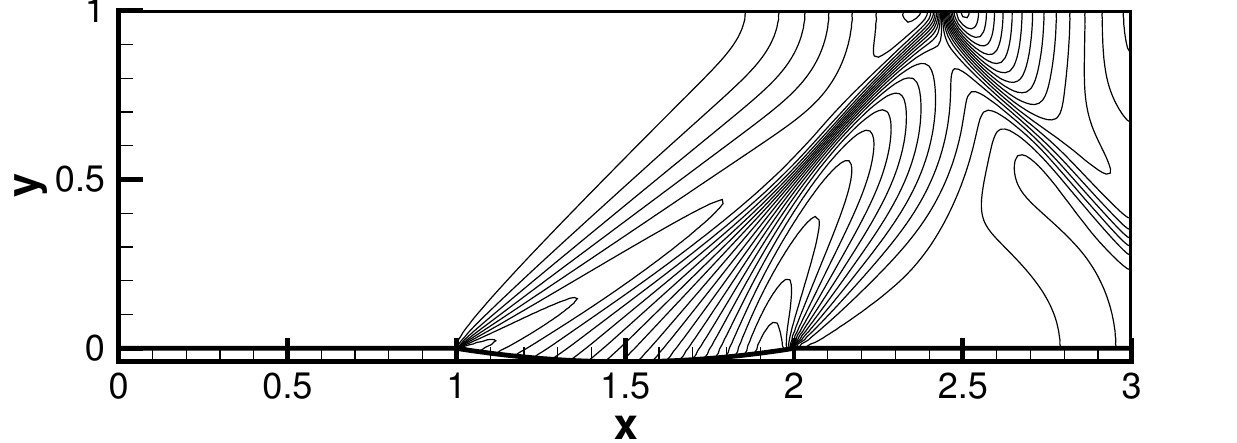}
\caption{(Top) geometry of the reverse bump and (bottom) Mach 1.4 pressure contours (0.4:0.032:1.2) using $240 \times 80 $ quadrilateral mesh.}
\label{fig:RevBumpMRSU}
\end{figure}The incident shock wave reflects from the top wall. From the rear end of the reverse bump a second expansion wave originates which interacts with the reflected shock wave. The important feature of this test case is formation of oblique shock wave through isentropic compression. Figure \ref{fig:RevBumpMRSU} shows the pressure contours on $240 \times 80 $ quadrilateral mesh.

\subsection{Two-dimensional Riemann Problems \cite{LW}}
\begin{figure} [h!] 
\centering
\includegraphics[scale=0.35]{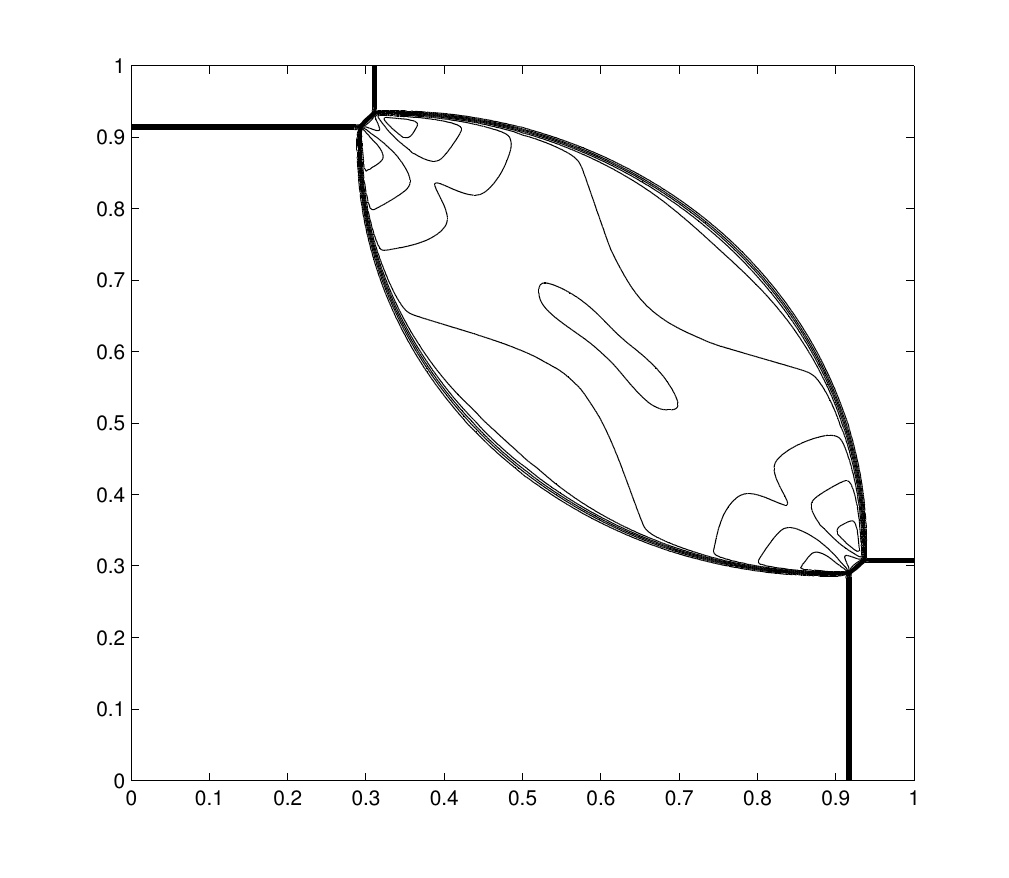}
\includegraphics[scale=0.35]{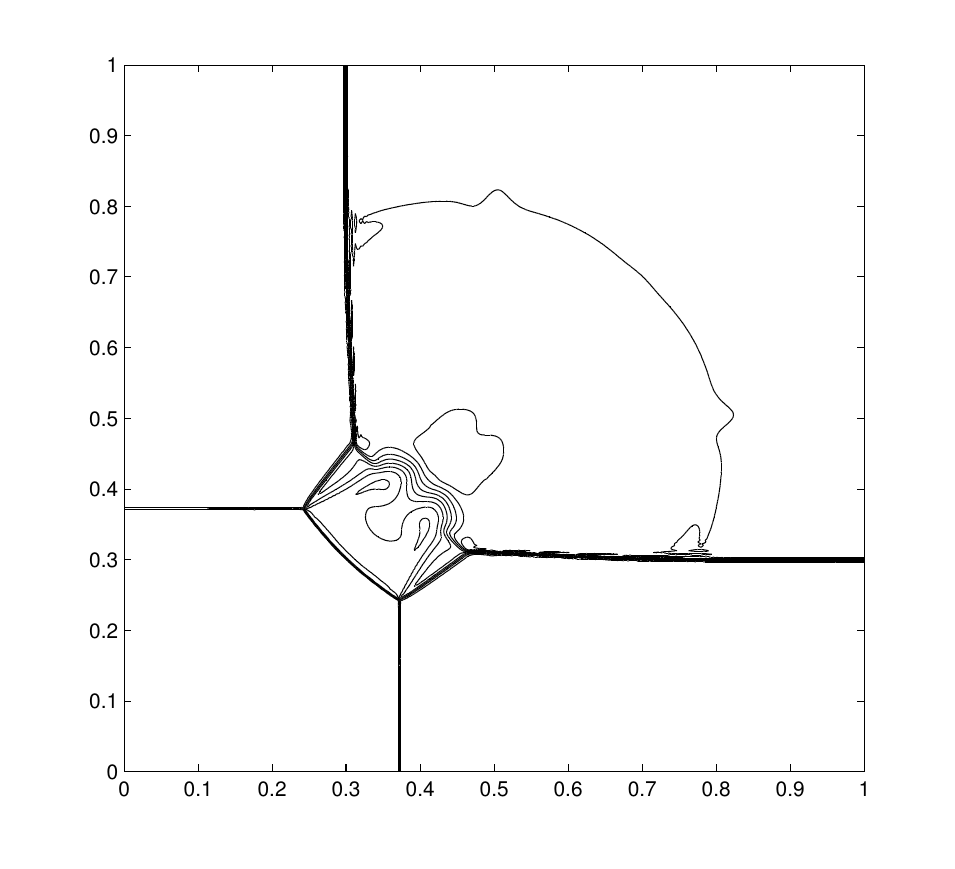}
\caption{Density contours for two-dimensional Riemann problems using $400 \times 400 $ mesh.}
\label{fig:RieProbMRSU}
\end{figure}
Two cases of two-dimensional Riemann problems involving all shock waves are considered.
Riemann problems are solved on domain $[0,\,1]^2$. This square domain is divided into four quadrants where initial constant states are defined. These problems are proposed in such a way that the solution between these quadrants have only one wave. The initial conditions are given in \cite{LW}. 
Figure \ref{fig:RieProbMRSU} shows the density contours for these cases. In both the cases, shock waves are accurately captured by the proposed scheme.


\subsection{Hydraulic Jump \cite{APGN}}
\begin{figure} [h!] 
\centering
\includegraphics[scale=0.28]{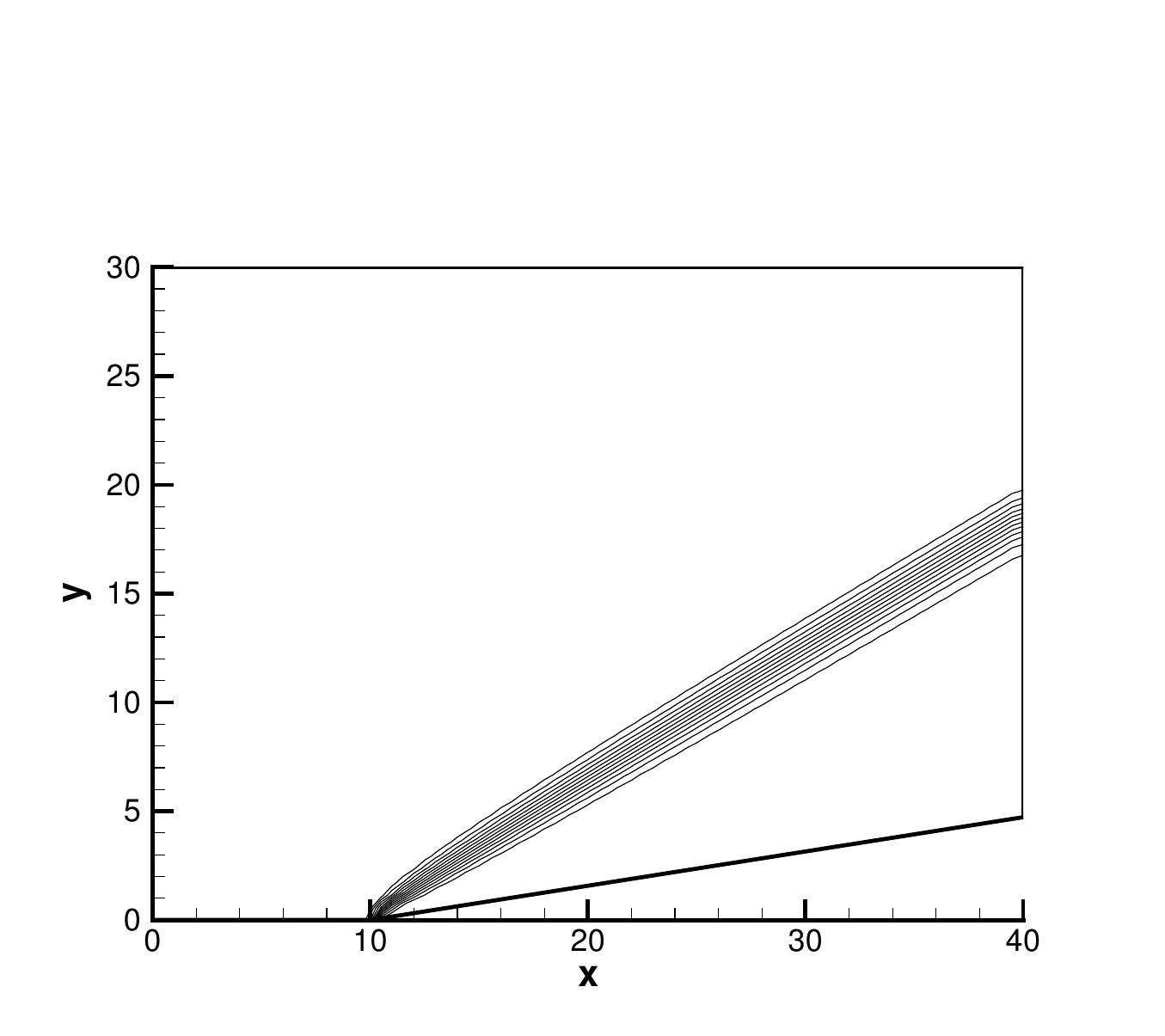}
\includegraphics[scale=0.37]{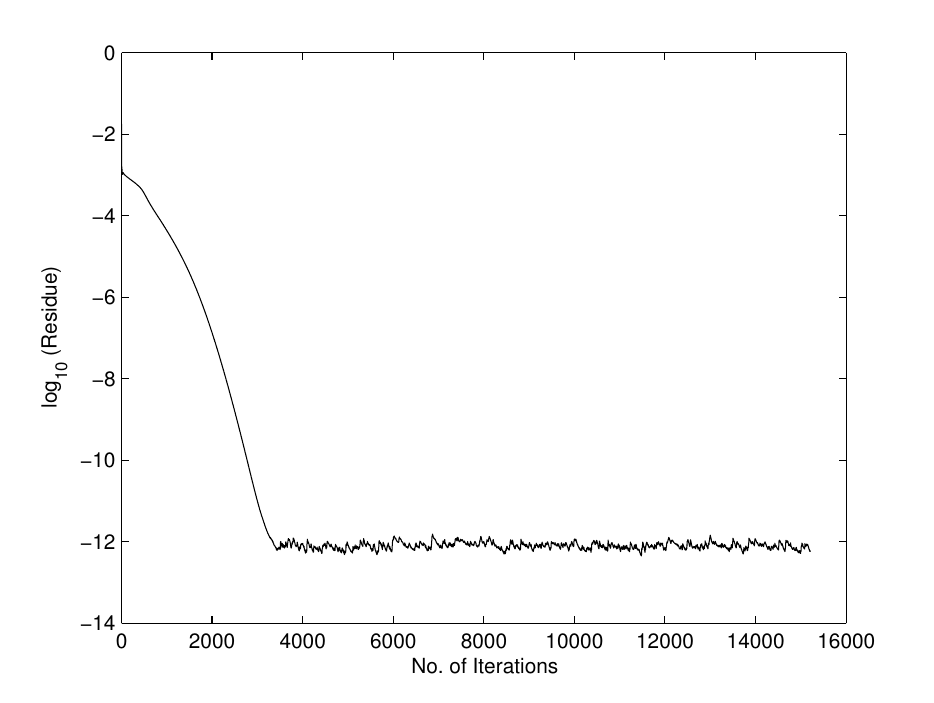}
\caption{Hydraulic jump using $80 \times 60$ quadrilateral mesh and residue plot.}
\label{fig:HjUMP}
\end{figure}
This test case gives a hydraulic jump in a convergent wall section. In this case the flow is supercritical. The computational domain is $[0,\,40] \times [0,\,30]$ and it is discretized using $80 \times 60 $ mesh.
The initial conditions are $H = 1\, m$, $u_1 = 8.57 \,m/s$ and  $u_2 = 0$. The wall angle for the convergent section is  $8.95^o$. The boundary conditions are reflective at top and bottom walls, whereas supercritical boundary conditions
are imposed on inflow and outflow boundaries. Figure \ref{fig:HjUMP} shows the contours of height $H$ along with the residue plot. The hydraulic jump is captured very well with such crude grid.

\subsection{Hydraulic Jump Interaction in a Convergent Channel}
This test case is a slight modification of previous test case where the walls are converging from both top and bottom with fixed angle $8.95^o$ which forms a convergent channel. Again, the flow is supercritical. The computational domain is $[0,\,50] \times [0,\,30]$ and it is discretized using $100 \times 60 $ mesh.
\begin{figure} [h!] 
\centering
\includegraphics[trim=1cm 0.5cm 1.5cm 7.5cm, clip=true, scale=0.33, angle = 0]{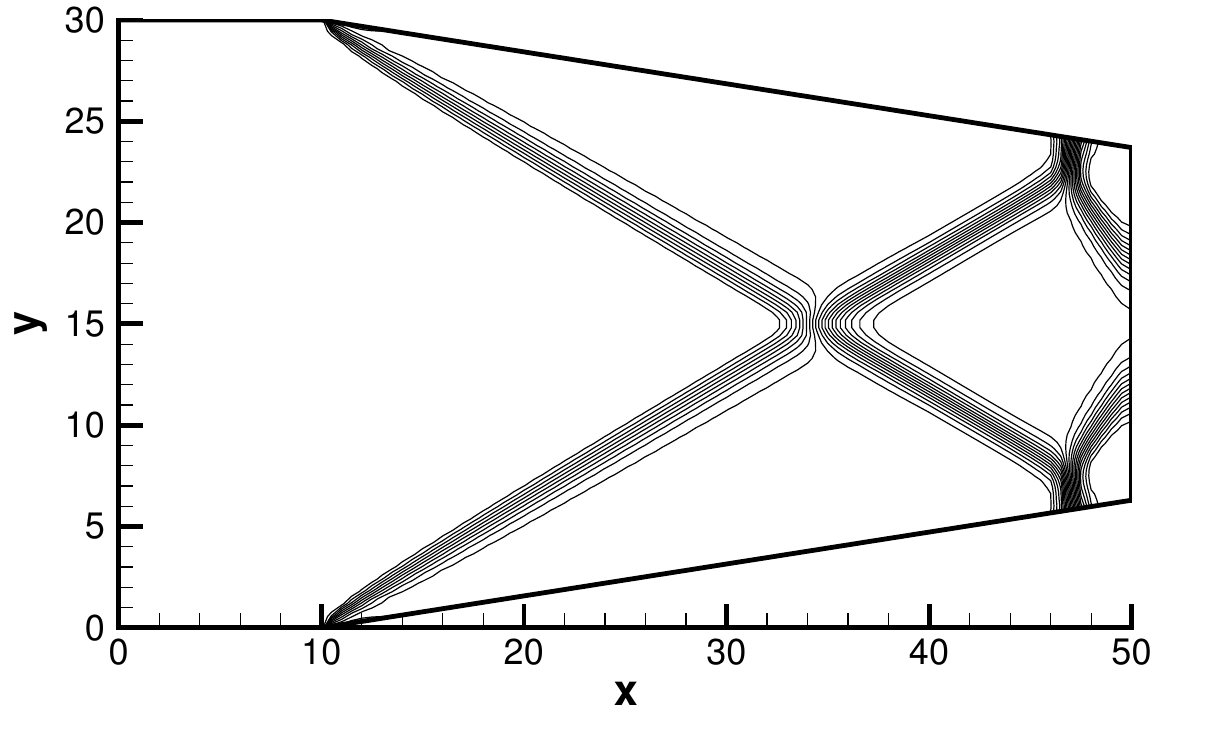}
\includegraphics[trim=5.5cm 5.5cm 7cm 5.5cm, clip=true, scale=0.48, angle = 0]{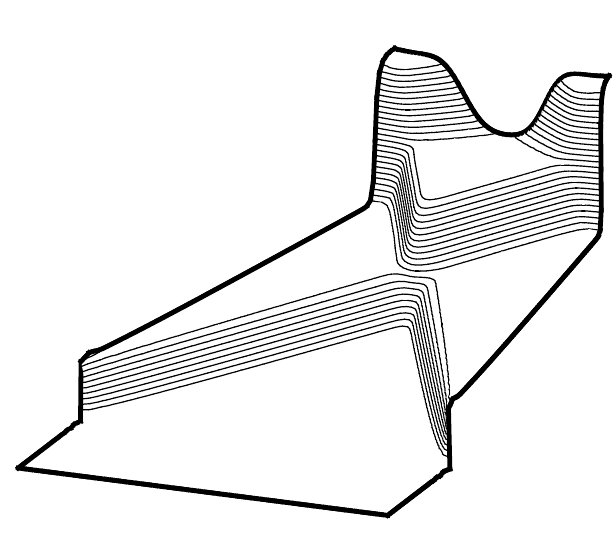}
\caption{Contour and surface plots for hydraulic jump interaction problem using $100 \times 60$ mesh.}
\label{fig:HjUMPint}
\end{figure}
The initial conditions are same as before. The boundary conditions are reflective at top and bottom walls, whereas supercritical boundary conditions are applied to the inflow and outflow boundaries. Hydraulic jump generated from top and bottom wall intersects each other and then reflects from the top and bottom boundaries. This interaction of hydraulic jump is captured well by MRSU scheme. Contour plot of height and surface plot are given in figure \ref{fig:HjUMPint}.

\subsection{Partial Dam Break Test Case \cite{APGN}}
This test case simulates the partial dam break due to sudden opening of sluice gate in a rectangular channel.
\begin{figure} [h!] 
\centering
\includegraphics[scale=0.28]{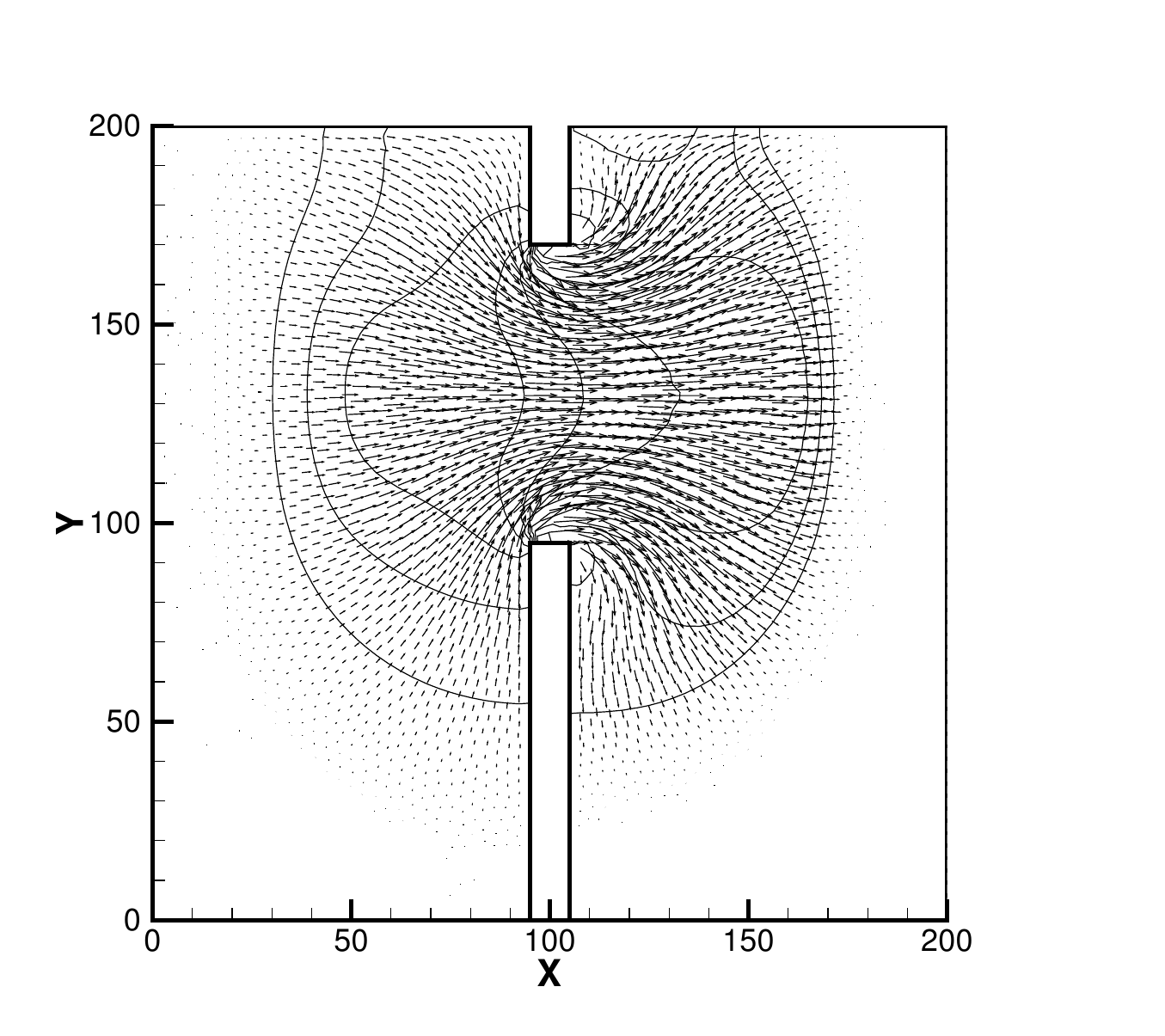}
\includegraphics[trim=0cm 0.4cm 0cm 0.9cm, clip=true, scale=0.52, angle = 0]{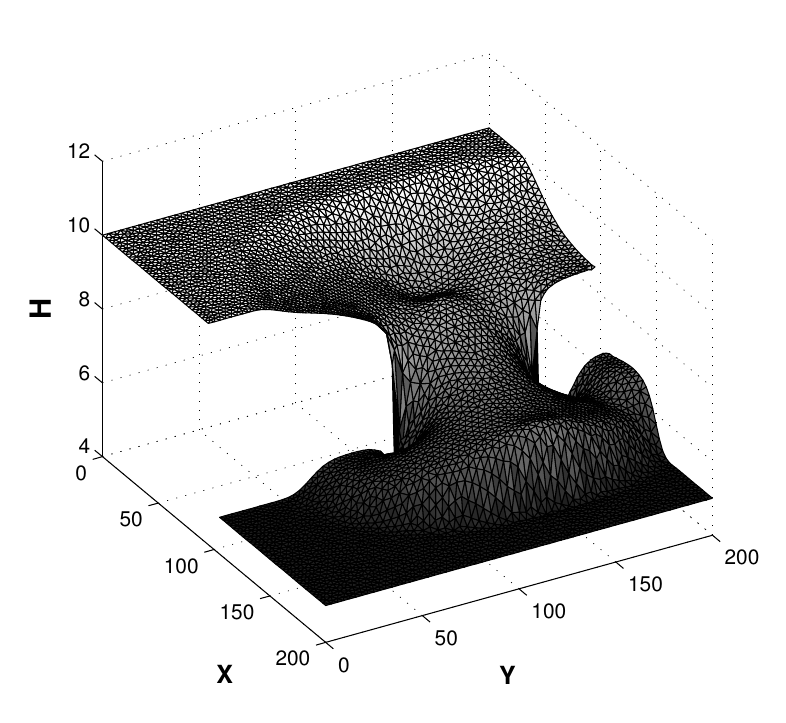}
\caption{Contour plot along with velocity vector plot (left) and surface plot of water column height (right) for partial dam break problem.}
\label{fig:PDB}
\end{figure}
The domain is $[0,\,200 \,m]^2$ and the length of sluice gate is 75m long as shown in first figure of ~\ref{fig:PDB}. The initial condition is
$$H = \begin{cases} 10\, m & \textrm{If} \,\,\,x< 100\, m \\ 5\, m & \textrm{Otherwise} \end{cases} $$
and $u_1 = u_2 = 0$.
Boundary conditions are transmissive on left as well as right boundaries, whereas reflective boundary conditions are used for remaining boundaries. Strong bore in the upstream and negative waves in the downstream directions are created due to sudden opening of sluice gate. The solution is plotted at t = 7.2 seconds. The unstructured triangular mesh is used for the simulation. The number of nodes is 4238 and the number of elements is 8114. Figure \ref{fig:PDB} shows the height contours along with velocity vector plot and surface plot respectively.

\section{Conclusions}
A new stabilized finite element method namely, Method of Relaxed Streamline Upwinding (MRSU) is proposed for hyperbolic conservation laws. The proposed scheme is based on relaxation system which replaces hyperbolic conservation laws by semi-linear system with stiff source term. Six symmetric discrete velocity models are also introduced in two dimensions which symmetrically spread foot of the characteristics in all four quadrants thereby taking information symmetrically from all directions. Main aim of this paper is to introduce a new idea based on relaxation scheme which is not only potentially interesting but also easy to implement. There are several advantages of the proposed scheme.
\begin{enumerate}
\item MRSU scheme gives exact diffusion vectors which are required for stabilization of the numerical scheme.
\item Diffusion vectors in MRSU scheme are simply conserved variable vectors with supremum eigenvalue of corresponding Jacobian matrices as a coefficient. This makes the scheme very robust. Moreover, computationally expensive Jacobian matrices are not involved in the diffusion terms. 
\item Extension of MRSU scheme from scalar to vector conservation laws is direct.
\end{enumerate}
The efficacy of the proposed scheme is shown by solving various one and two-dimensional test cases for Burgers equation (which includes both convex and non-convex flux functions), Euler equations and shallow water equations. Moreover, new test cases are proposed for all the three equations. In case of two-dimensional Burgers equation, set of test cases involving normal and oblique discontinuities is proposed along with their exact solutions. For Euler equations, the variations in Mach number are from 0.5 to 20 which covers subsonic, transonic and supersonic regimes. This shows the capability of proposed scheme for handling wide range of problems. Moreover, spectral stability analysis is carried out which gives expression of critical time step. Error analysis shows optimal convergence rate for the proposed scheme. The proposed scheme is easy to implement in the existing code of stabilized finite element methods without much modification. 

\section*{Acknowledgments}
The author would like to thank Prof. A.S.Vasudeva Murthy, Tata Institute of Fundamental Research - Centre for Applicable Mathematics, Bangalore and Prof. C.S. Jog, Department of Mechanical Engineering, Indian Institute of Science, Bangalore for giving constructive comments and suggestions which helped to improve the quality of the manuscript. The author is also grateful to the reviewers for their valuable suggestions.

\section*{Appendix A : Derivation of local Maxwellian Distribution ($\mathbf{F}$) for one-dimensional three discrete velocity model}

Let $ \mathbf{F} = \{F_1\,\,F_2\,\,F_3\}^T$ be the local Maxwellian distribution function then, the moment relations (equation \eqref{Mome111}) are satisfied as
\begin{subequations}\label{mo}
\begin{align}
 F_1 +F_2+ F_3 &= U  \label{mo1}
 \\ \lambda F_1 +\lambda F_2+\lambda F_3 &= G(U) \label{mo2}
\end{align}
\end{subequations}
$\Lambda$ matrix for three discrete velocity model is
\begin{equation*}
\Lambda = \text{diag}\{-\lambda, 0, \lambda\} 
\end{equation*}
which makes equation \eqref{mo2} as
\begin{equation}\label{mo3}
F_3- F_1 = \frac{G(U)}{\lambda}
\end{equation}
Assuming $\mathbf{F}$ as a linear combination of $U$ and $G(U)$,
\begin{equation*}
\mathbf{F} = \left\{ {\begin{array}{c}
F_1\\
F_2\\
F_3 \\
\end{array} } \right\} = \left\{ {\begin{array}{c}
c_1 U + d_1G(U)\\
c_2 U + d_2G(U)\\
c_3 U + d_3G(U)\\
\end{array} } \right\}
\end{equation*}
Using moment relations \eqref{mo1} and \eqref{mo3}, following relations are obtained
\begin{align*}
\sum_{i=1}^3 c_i = 1, \ \ \ \ \ \  \sum_{i=1}^3 d_i = 0, \nonumber
\\ (c_3-c_1)U + (d_3-d_1)G(U)= \frac{G(U)}{\lambda} \nonumber
\end{align*}
Choosing $ c_1 = c_2=c_3 = \frac{1}{3} $ and $ d_1 = \frac{-1}{2\lambda} = -d_3,\,\,
 d_2 = 0$, we get
\begin{equation*}
\mathbf{F} =   \left\{ {\begin{array}{c}
\frac{U}{3} - \frac{G}{2\lambda}\\
\frac{U}{3} \\
\frac{U}{3} + \frac{G}{2\lambda}\\
\end{array} } \right\}
\end{equation*}

\section*{Appendix B : Chapman-Enskog type Expansion of Relaxation System}
Consider one-dimensional scalar hyperbolic conservation law
\begin{equation*}
\frac{\partial U}{\partial t} + \frac{\partial G(U)}{\partial x} = 0, \,\,\, (x,t)  \in \Omega_T \subset \mathbb{R}\times \mathbb{R}_+
\end{equation*}
with appropriate initial and boundary conditions. The relaxation system for above conservation law is 
\begin{subequations}
\begin{align}
\frac{\partial U}{\partial t} + \frac{\partial W}{\partial x} & = 0  \label{RS1d11}
\\ \frac{\partial W}{\partial t} + \lambda^2 \frac{\partial U}{\partial x} &= -\frac{1}{\epsilon} (W-G(U)) , \ \ \epsilon \rightarrow 0  \label{RS1d}
\end{align}
\end{subequations}
Rewriting equation \eqref{RS1d}, we get
\begin{equation}\label{ab1}
W = G(U) - \epsilon \left( \frac{\partial W}{\partial t} + \lambda^2 \frac{\partial U}{\partial x}\right) 
\end{equation}
which is 
\begin{equation}\label{ab2}
W = G(U) + \mathcal{O}(\epsilon)
\end{equation}
Differentiating with respect to time, we get
\begin{align}
\frac{\partial W}{\partial t}  & = \frac{\partial  G(U)}{\partial t} + \mathcal{O}(\epsilon) \nonumber
\\ & = \frac{\partial  G(U)}{\partial U} \frac{\partial  U}{\partial t} + \mathcal{O}(\epsilon) \nonumber
\\ & = -\frac{\partial  G(U)}{\partial U} \frac{\partial  W}{\partial x} + \mathcal{O}(\epsilon) \nonumber
\end{align}
Using equation \eqref{ab2} we can write
\begin{align}
\frac{\partial W}{\partial t}  & =  -\frac{\partial  G(U)}{\partial U} \left \{\frac{\partial [G(U) + \mathcal{O}(\epsilon)]  }{\partial x} \right\}+ \mathcal{O}(\epsilon) \nonumber
\\ &  =  -\frac{\partial  G(U)}{\partial U} \left \{\frac{\partial G(U)  }{\partial U}\frac{\partial U  }{\partial x}  \right\}+ \mathcal{O}(\epsilon) \nonumber
\\ &  =  -\left(\frac{\partial  G(U)}{\partial U} \right)^2\frac{\partial U  }{\partial x} + \mathcal{O}(\epsilon) \nonumber
\end{align}
Substituting above equation in \eqref{ab1}, we get
\begin{align*}
W & = G(U) - \epsilon \left( \left[ -\left(\frac{\partial  G(U)}{\partial U} \right)^2\frac{\partial U  }{\partial x} + \mathcal{O}(\epsilon) \right] + \lambda^2 \frac{\partial U}{\partial x}\right) 
\\ & = G(U) - \epsilon \left[   \lambda^2 -\left(\frac{\partial  G(U)}{\partial U} \right)^2  \right] \frac{\partial U  }{\partial x}  - \mathcal{O}(\epsilon^2) 
\end{align*}
Substituting this in equation \eqref{RS1d}, we get
\begin{align*}
\frac{\partial U}{\partial t} + \frac{\partial G(U)}{\partial x} &= \epsilon \frac{\partial }{\partial x} \left( \frac{\partial U}{\partial x} [\lambda^2 - (G'(U))^2] \right)+ \mathcal{O}(\epsilon^2)
\end{align*}

\section*{Appendix C : Weak Formulation of System of Hyperbolic Conservation Laws in Conservation Form}
The standard Galerkin finite element approximation for conserved variable $U$ and the flux function $G_j(U)$ are
\begin{equation*}
U \approx U^h = \sum_{\forall i} N_i^h U_i^h, \ \ \ \  G_j \approx G_j^h = \sum_{\forall i} N_i^h (G_j)_i^h
\end{equation*}
where group formulation is used for flux function. 
Defining test and trial functions as 
\begin{align*}\mathcal{V}^h&=\{N^h~\in~\mathcal{H}^{1}( \Omega^h) \linebreak \, \textrm{and} \, \ N^h  = 0 \,\,\text{on}\,\, \Gamma_D \}
\\ \mathcal{S}^h  &=  \{  U^h \in \mathcal{H}^{1}( \Omega^h) \, \textrm{and} \, U^h = U^h_D \,\,\text{on}\,\, \Gamma_D \}  
\end{align*}
the weak formulation is, find $U^h \in \mathcal{S}^h $ such that $\forall \, N^h \in \mathcal{V}^h$
\begin{align*}
\int_{\Omega^h} N^h \cdotp   \left(\frac{\partial U^h}{\partial t}  +    \frac{\partial G_1^h}{\partial x} +\frac{\partial G_2^h}{\partial y}\right)\, d\Omega^h + \sum_{e=1}^{\textit{Nel}} \int_{\Omega_e^h} & \left[ \tau_1 A_1^h \frac{\partial N^h}{\partial x} +\tau_2 A_2^h \frac{\partial N^h}{\partial y}\right] \cdotp \left(    \frac{\partial G_1^h}{\partial x} + \frac{\partial G_2^h}{\partial y} \right)\, d\Omega_e^h 
\\ & + \sum_{e=1}^{\textit{Nel}}  \int_{\Omega_e^h}  \delta^e \left(  \frac{\partial N^h}{\partial x}\cdotp\frac{\partial U^h}{\partial x} + \frac{\partial N^h}{\partial y}\cdotp\frac{\partial U^h}{\partial y}\right) \, d\Omega_e^h  = 0
\end{align*}
where $A_i =\frac{\partial G_i}{\partial U}, \, i=1, 2$ are the flux Jacobian matrices. 
For the reference, this scheme is called as Method of Streamline Upwinding (MSU) scheme which does not use relaxation framework.  Unlike in MRSU scheme, construction of stabilization parameter $\tau$ in MSU scheme is not a trivial task. The element level $\tau$ is constructed based on the norm definition given by Tezduyar and Osawa \cite{TOnew}. This results in further increase of computational cost.

Flux Jacobian matrices for two-dimensional Euler equations are 
\begin{equation*}
A_1(U)  =   \left [ 
\begin{array}{cccc}
0&1&0 &0 \\
-u_1^2 + \frac{\gamma-1}{2} (u_1^2+u_2^2) & (3-\gamma) u_1&-(\gamma-1)u_2 &\gamma-1 \\
-u_1 u_2&u_2&u_1  &0 \\
-(\gamma E - (\gamma-1)(u_1^2+u_2^2)) u_1 &\gamma E - \frac{\gamma-1}{2}(2u_1^2 (u_1^2+u_2^2)) &-(\gamma-1)u_1 u_2& \gamma u_1\\
\end{array} \right ]  \nonumber
\end{equation*}
and
\begin{equation*}
A_2(U)  =  \left [ 
\begin{array}{cccc}
0&0&1 &0 \\
-u_1 u_2&u_1&u_2  &0 \\
-u_2^2 + \frac{\gamma-1}{2} (u_1^2+u_2^2) &-(\gamma-1)u_1 & (3-\gamma) u_2&\gamma-1 \\
-(\gamma E - (\gamma-1)(u_1^2+u_2^2)) u_2 &-(\gamma-1)u_1 u_2&\gamma E - \frac{\gamma-1}{2}(2u_2^2 (u_1^2+u_2^2)) & \gamma u_2\\
\end{array} \right ]  \nonumber
\end{equation*}
whereas for two-dimensional shallow water equations these matrices are given as
\begin{equation*}
A_1(U) = \left [ 
\begin{array}{ccc}
0&1 &0 \\
gH-u_1^2 & 2u_1 & 0\\
-u_1u_2 & u_2 & u_1\\
\end{array} \right  ] 
\end{equation*}
and
\begin{equation*}
A_2(U) = \left [ 
\begin{array}{ccc}
0&0 &1 \\
-u_1u_2 & u_2 & u_1\\
gH-u_2^2 & 0&2u_2 \\
\end{array} \right  ] 
\end{equation*}


\begin{thebibliography}{99}
\bibitem{APGN}
Alcrudo F. and P. Garcia-Navarro, A high resolution Godunov type scheme in finite volumes for 2D shallow-water equations, International Journal for Numerical Methods in Fluids 16(6), 489 - 505 (1993).

\bibitem{AKe}
Alexander Kurganov, \textit{et. al.}, Adaptive Semidiscrete Central-Upwind Schemes for Non-convex Hyperbolic Conservation Laws, SIAM J. Sci. Comput. , Vol. 29, No. 6, pp. 2381-2401, 2007.

\bibitem{PAG}
Aregba-Driollet D., Natalini R., Discrete Kinetic schemes for Multidimensional systems of conservation laws, SIAM Journal of Numerical Analysis, vol. 37, no. 6, 1973-2004, 2000.

\bibitem{CATM}
Arvanitis C. \textit{et. al.}, Adaptive finite element relaxation schemes for hyperbolic conservation laws, Mathematical Modeling and Numerical Analysis, ESAIM:M2AN, Vol. 35, No. 1, 2001, pp 17-33.

\bibitem{MB}
Banda M., Variants of relaxed schemes and two dimensional gas dynamics, Journal of Computational and Applied Mathematics, vol.175, no. 1, 41-62, 2005.


\bibitem{MBAK}
Banda M.K. \textit{et. al.}, A lattice- Boltzmann relaxation scheme for coupled convection-radiation systems, Journal of Computational Physics, 226 (2007) 1408 - 1431.



\bibitem{Kirk}
Benjamin S. Kirk and Graham F. Carey. Development and Validation of a SUPG Finite Element Scheme for the
Compressible Navier-Stokes Equations using a Modified Inviscid Flux Discretization. International Journal for Numerical
Methods in Fluids, 57(3): 265 - 293, 2008.

\bibitem{GaB}
Ben-Dor G.,  Shock wave reflection phenomena, 2nd Ed., Shock wave and high pressure phenomena, Berlin, Springer 2007.

\bibitem{BFLS}
Bereux F. and Sainsaulieu L., A Roe-type Riemann solver for hyperbolic systems with relaxation based on a time dependent wave decomposition, Numerische Mathematik, vol. 77, 143-185, 1997.

\bibitem{BGK} 
Bhatnagar P.L., E.P. Gross, M. Krook, A model for collision processes in gases I. Small amplitude processes in charged and neutral one-component systems, Phys. Rev. 94(1954) 511-525.


\bibitem{FSBillig}
Billig F.S., Shock-wave shapes around spherical and cylindrical-nosed bodies, Journal of Spacecraft and Rockets, Vol. 4, No. 6 (1967), pp. 822-823.


\bibitem{PBGU}
Bochev P.B. , Gunzburger M.D., Least-Squares Finite Element Methods. Springer, 2009.


\bibitem{BF}
Bochut F., Construction of BGK model with a family of kinetic entropies for a given system of conservation laws. Journal of Statistical Physics 1999, 95: 113-170.





\bibitem{CF}
Cavalli F., \textit{et. al.}, A family of relaxation schemes nonlinear convection diffusion problems. Communications in Computational Physics,  vol. 5, no. 2-4, 532-545, 2009.

\bibitem{LCC}
Catabriga L., Coutinho A.L.G.A., Improving convergence to steady state of implicit SUPG solution of Euler equations, Commun. Numer. Meth. Engng 2002, 18:345 - 353.

\bibitem{CRR}
Cecchi M.M., Redivo-Zaglia M. and Russo G., Extrapolation methods for hyperbolic systems with relaxation, Journal of Computational and Applied Mathematics, vol. 66, no. 1-2, 359-375, 1996.

\bibitem{ACh}
Chalabi A., Convergence of relaxation scheme for hyperbolic conservation laws with stiff source terms, Math. Comput. 68 (1999), pp. 955-970.


\bibitem{CLL}
Chen G. Q., Levermore C.D. and Liu T.P., Hyperbolic conservation laws with stiff relaxation terms and entropy, Comm. Pure Appl. Math.. 47 (1994), pp. 787-830.

\bibitem{TJC}
Chung T.J., Finite element analysis in fluid dynamics.McGraw-Hill international, 1978.


\bibitem{CFL}
Courant R., Friedrichs K. and Levy H., Uber die partielle differentialgleichungen der matematischen physik, Math. Annal. 100 (1928) 32- 74.

\bibitem{CIR}
Courant R., Isaacson E. and Rees M., On the solution of nonlinear hyperbolic differential equations by finite differences, Comm. Pure Appl. Math. 5: 243- 255, 1952.

\bibitem{JJDB}
 Donea J., Huerta A. : Finite element methods for flow problems, Wiley, 2003.

\bibitem{DHP}
Deconinck, H., Paillere, H., Struijs, R., and Roe, P. (1993). Multidimensional upwind schemes based on fluctuation-splitting for systems of conservation laws. Computational Mechanics, 11:323 - 340.


\bibitem{TJF}
 Fan H., Jin S. and Teng Z., Zero reaction limit for hyperbolic conservation laws with source terms, J. Diff. Equations 168 (2000), pp. 270-294.

\bibitem{Fle}
 Fletcher C. A. J. The Group Finite Element Formulation. Computer Methods in Applied Mechanics and Engineering,
37: 225 - 243, 1983.

\bibitem{GB1}
Graille B., Approximation of mono-dimensional hyperbolic systems: A lattice Boltzmann scheme as a relaxation method, Journal of Computational Physics, 266 (2014) 74  -  88.

\bibitem{LGg}
Gvozdeva L.G., Predvoditeleva O.A. and Fokeev V.P., Double Mach reflection of strong shock waves, Fluid Dynammics, 3 (1): 6-11, 1968.




\bibitem{JSHTW}
Hesthaven J.S., Warburton T., Nodal Discontinuous Galerkin Methods. Springer, 2008.


\bibitem{ChH}
Hirsch C., Numerical computation of internal and external flows, Vol. 2: Compuatational methods for inviscid and viscous flows, John Wiley \& Sons, 1990.

\bibitem{HR}
Hittenger J.A.F. and Roe P.L., Asymptotic analysis of the Riemann problem for the constant coefficient hyperbolic systems with relaxation, ZAMM-Journal of Applied Mathematics and Mechanics, vol. 84, no. 7, 452-471, 2004.













\bibitem{Jin}
Jin S., Xin Z., The relaxation schemes for system of conservation laws in arbitrary space dimensions, Communications in pure and applied Mathematics, vol. 48, no. 3, 235-277, 1995.

\bibitem{Kuz}
Kuzmin D., Moller M., and Turek S., High–Resolution FEM-FCT Schemes for Multidimensional Conservation Laws.
Computer Methods in Applied Mechanics and Engineering, 193:4915 - 4946, May 2004.

\bibitem{Kuz1}
Kuzmin D. and  Turek S., High–Resolution FEM-TVD Schemes Based on a Fully Multidimensional Flux Limiter. Journal of Computational Physics, 198:131 - 158, 2004.








\bibitem{CL}
Laney C.B., Computational Gasdynamics, Cambridge University Press, 1998.

\bibitem{CLD}
Lattanzio C. and Serre D., Convergence of a relaxation scheme for hyperbolic systems of conservation laws, Numer. Math., $\mathbf{88}$(2001) 121 - 134.



\bibitem{LMP}
Leveque R.J. and Pelanti M., A class of approximate Riemann solvers and their relation to relaxation schemes, Journal of Computational Physics, vol 172, 572-591, 2001.



\bibitem{HLi}
Li H. and Ben-Dor G., A shock dynamics theory based analytical solution of double Mach reflections, Shock Waves, 5 (4): 259-264, 1995.

\bibitem{CYL}
Li, X., Yu X. and Chen G., The third order relaxation schemes for hyperbolic conservation laws, J. Comp. Appl. Math., 138 (2002), pp. 93-108.

\bibitem{LW}
Liska, R. and Wendroff, B., Comparison of several difference schemes in 1D and 2D Test Problems for Euler Equations. SIAM Journal of Scientific Computation, 25:995 - 1017, 2003.



\bibitem{Mor}
Morgan K.  and Peraire J., Unstructured grid finite element methods for fluid mechanics, 1998.


\bibitem{AVS}
Murthy A. S. V., An alternative relaxation approximation to conservation laws, Journal of Computational and Applied Mathematics, 203 (2007) 437-443

\bibitem{NaR}
Natalini R., A Discrete Kinetic Approximation of Entropy Solutions to Multidimensional Scalar Conservation Laws, Journal of Differential Equations, vol. 148, 292-317, 1999.

\bibitem{Nat}
Natalini R., Recent mathematical results on hyperbolic relaxation problems, Analysis of systems of conservation laws, Pitman Research Notes in Mathematics Series, Longman, Harlow, 1998.

\bibitem{RN1}
Natalini R., Convergence to equilibrium for relaxation approximations of conservation laws, Comm. Pure Appl. Math., $\mathbf{49}$ (1996) 795 - 823.


\bibitem{RM}
Richtmyer R. D., Morton K.W., Difference methods of initial-value problems, Wiley, 1967.


\bibitem{HS}
Schroll H., High resolution relaxed upwind schemes in gas dynamics, J. Sci. Comput. 17 (2002), pp. 599-607.

\bibitem{MSeaid}
Seaid M., Improved applications of relaxation schemes for hyperbolic systems of conservation laws and convection-diffusion problems, Computational methods in applied mathematics, vol.6 (2006), No. 1, pp.56-86.

\bibitem{CWSo}
Shu C.W. and Osher S.: ‘Efficient implementation of essentially non-oscillatory shock-capturing schemes, II, Journal of Computational Physics, 83, 32–78 (1989).

\bibitem{Spekreijse} Spekreijse S., Multigrid solution of monotone second order discretizations of hyperbolic conservation laws, Mathematics of Computation, vol. 49, No. 179, pp. 135-155, 1987.  


\bibitem{TH} 
Tezduyar TE, Huges TJR, .Finite element formulation of convected dominated flows with particular emphasis on the compressible Euler equations. In: Proceedings of AIAA 21st aerospace sciences meeting. AIAA Paper 83-0125, Reno, Nevada (1982).


\bibitem{TOnew} 
Tezduyar T.E., Osawa Y., Finite element stabilization parameters computed from element matrices and vectors. \emph{Comput. Methods Appl. Mech. Eng.} $\mathbf{190}$: 411-430 (2000).


\bibitem{TS}
Tezduyar TE, Senga M., SUPG finite element computation of inviscid supersonic flows with $YZ\beta$ shock-capturing. \emph{Computers and Fluids,} 36: 147-159, 2007.







\bibitem{TOro}
Toro, E.F., Riemann Solvers and Numerical Methods for Fluid Dynamics: A practical Introduction, Springer -Verlag, 3rd edition, 2009,



\bibitem{HVi}
Viviand H., Numerical solution of two dimensional reference test cases. AGARD AR-211: Test cases for inviscid flow field methods, 1985.


\bibitem{WCPP}
Woodward P., Colella P., The numerical solution of two dimensional fluid flow with strong shocks, Journal of Compu. Phys., vol. 54, No.1, 115 - 173, 1984.


\bibitem{WQX}
Xu W.Q., Relaxation limite for piecewise smooth solutions to conservation laws, J. Diff. Eq. $\mathbf{162}$ (2000), 140 - 173.

\bibitem{Yee_Warming_Harten}
Yee H.C., Warming R.F. and  Harten A., A high-resolution numerical technique for inviscid gas-dynamics problems with weak solutions, Proceedings in eight international conference on numerical methods in fluid dynamics, Lecture notes in Physics, Vol. 170, Springer, New York/Berkin, pp. 546-552, 1982.



\bibitem{OCZT}
Zienkiewicz O.C., Taylor R.L., The Finite Element Method, Volume 3, $6^{th}$ edition, Butterworth-Heinemann, 2005.


\bibitem{ZGM}
Zhang, S.Q., Ghidaoui, M.S., Gray, W.G., Li,N.Z., A kinetic flux vector splitting scheme for shallow water flows, Advances in Water Resources 26 (2003) 635 - 647.

\bibitem{ZO}
Zoppou, C., Roberts, S., Explicit schemes for dam-break simulations, J. Hydraul. Eng. 129 (1), 11 - 34 (2003).


\end{thebibliography}

\end{document}